\documentclass[12pt]{report}
\usepackage{amsfonts,amsmath, eucal, amssymb}

\usepackage{suthesis-2e}

\usepackage{thesis_chatterjee_macros}

\begin{document}

    \title{Concentration inequalities with \\ exchangeable pairs}
    \dept{Statistics}
    \author{Sourav Chatterjee}
    \principaladviser{Persi Diaconis}
    \firstreader{Amir Dembo}
    \secondreader{David O. Siegmund}
    \submitdate{June 2005}
    \copyrightyear{2005}
    \figurespagefalse
    \tablespagefalse
 
    \beforepreface
    \prefacesection{Abstract}
    The purpose of this dissertation is to introduce a version of Stein's method of exchangeable pairs to solve problems in measure concentration. We specifically target systems of dependent random variables, since that is where the power of Stein's method is fully realized. Because the theory is quite abstract, we have tried to put in as many examples as possible. Some of the highlighted applications are as follows:~(a)~We shall find an easily verifiable condition under which a popular heuristic technique originating from physics, known as the ``mean field equations'' method, is valid. No such condition is currently known.~(b)~We shall present a way of using couplings to derive concentration inequalities. Although couplings are routinely used for proving decay of correlations, no method for using couplings to derive concentration bounds is available in the literature. This will be used to obtain~(c)~concentration inequalities with explicit constants under Dobrushin's condition of weak dependence.~(d)~We shall give a method for obtaining concentration of Haar measures using convergence rates of related random walks on groups. Using this technique and one of the numerous available results about rates of convergence of random walks, we will then prove~(e)~a quantitative version of Voiculescu's celebrated connection between random matrix theory and free probability.

    \prefacesection{Acknowledgements}
First of all, I would like to thank my parents for bringing me up the way they did, because, despite all the deficiencies in my character I am sure it is immensely better than what I could have managed on my own. They don't know any math at all, but they have the sincere respect for academics that is so pleasantly ubiquitous in the land from where I hail. It is also worthy of note, dear reader, that they love me more than anything else in the universe.

Having said that, let me add that my character and enthusiasm alone would not have made me a mathematician. There is one person who is overwhelmingly responsible for that, and that is my adviser, Persi Diaconis. I recall that Persi once told me, towards the beginning of my Ph.D.\ quest, that his job was to take me to the ``top''. Well, there is still a long way to go, but I am grateful to him for at least showing me with his own attitude what the ``top'' demands and bringing it closer to the earth. I also thank him for numerous other things, not the least for reading this thesis with an unbelievable amount of diligence and care. 

I would also like to thank professors Amir Dembo and David Siegmund for agreeing to be readers for this dissertation, and also for the numerous things I learned over the past three years from these two astute mathematicians. Amir has been particularly kind to me on many occasions, sparing a lot of his time on academic discussions and giving advice on various issues.

Charles Stein has been an inspiration in all the work that I did during my Ph.D. 

Besides that, I thank all the other professors and staff in this wonderful department who have never flinched about answering a question from a nagging graduate student --- something that is embedded in the culture of Sequoia Hall.

I am grateful to professor Arup Bose from the Indian Statistical Institute, Kolkata, for initiating my research career by helping me write my first real paper.

It is imperative that I thank all the professors who answered my various questions by e-mail; in that respect, I am particularly grateful to professor Michel Ledoux. The other members of this illustrious list include Michel Talagrand, Boguslaw Zegarlinski, Kenneth Alexander and David Aldous. I hope I have not missed anyone out!

Thanks are also due to Dadabhai (Rajib) and Boudi (Mukta), and of course, little Roompi, for taking care of me in various ways in the last three years. They add a touch of sanity to my world that is otherwise mired in academics.

Coming to friends, I have a few of them to whom I owe a bit (not money, fortunately). First of all, there is Soumik Pal, my closest confidante and friend, from whom I learned about life as much as I did about mathematics. He never fails to amaze me with his mathematical originality at regular intervals of time; for instance, I should mention that he has a beautiful alternative martingale argument for one of the basic theorems in this dissertation (Theorem \ref{hoeffding}).

Then there is Pinaki Biswas, and he may or may not be surprised to see his name in this list; but I am grateful to him, if for nothing else, then for a friendship that has lasted seventeen years now, and that is something too valuable by itself to discount.

I must also thank my three closest friends at Stanford --- George Chang, Nancy Ruonan Zhang, and Brit Turnbull. They are the ones most responsible for making my life livable in this desolate rural countryside.

Finally, there is one person from the past, an old retired clerk to whom I owe more than what I realized when he was alive. He was my childhood friend; in the long hours that we used to spend together, he introduced me to everything ranging from Euclid to atheism, from simultaneous equations to Shakespeare; in short, everything that inspired me to think and to dream. I would like to dedicate this thesis to the memory of my grandfather, the late Shri Tarapada Chatterjee.

    \afterpreface
 
\chapter{Introduction}\label{introduction}
The theory of concentration inequalities tries to answer the following question: Given a random variable $X$ taking value in some measure space $\xx$ (which is usually some high dimensional Euclidean space), and a measurable map $f:\xx \ra \rr$, what is a good explicit bound on $\pp\{|f(X)-\ee f(X)| \ge x\}$? Exact evaluation or accurate approximation is, of course, the central purpose of probability theory itself. In  situations where this is not possible, concentration inequalities aim to do the next best job.

A bound is {\it good} if it is sufficiently rapidly decreasing; gaussian bounds are usually considered satisfactory. The reasons for insisting on good bounds (as opposed to Chebychev type bounds) are theoretical as much as practical. In fact, the theoretical interest often supercedes the practical aspect, because in spite of all the activity, concentration bounds often give bad numbers when calculated numerically. Theoretically their importance stems, in large part, from the Bonferroni inequality: If we know that if a collection of events $\{A_k\}_{1\le k\le n}$ are so rare that $\pp\{A_k\} \le e^{-cn}$ for each $k$ for some fixed constant $c$, then $\pp\{\cup_k A_k\} \le ne^{-cn}$, which is again ``small'', since the $e^{-cn}$ term ``kills'' the $n$. Contrary to what someone unfamiliar with the literature might feel, this seemingly crude technique has been successful in establishing surprisingly efficient results, mainly because ``rare events are often approximately disjoint''. 

This was at the center of the earliest line of thought about controlling the suprema of empirical processes, developed mainly by David Pollard and others in the eighties (see Pollard \cite{pollard90} for an exposition). Subsequent developments in the nineties, based on Talagrand's concentration inequalities \cite{talagrand95, talagrand96b, talagrand96c} and the more recent ``entropy method'' of Ledoux \cite{ledoux01} and Massart \cite{massart00}, are more subtle. They will be discussed later in detail. All in all, concentration inequalities form the backbone of this very important branch of modern theoretical statistics and machine learning. In return, empirical process theory forms the most prominent area of application for concentration bounds.

Concentration inequalities are also used in theoretical computer science, random matrix theory and a variety of other fields, for reasons more or less the same as mentioned before. Often, applications get obscured because they were not mentioned in the abstract or the keywords, and that is a sign that concentration inequalities are rapidly attaining the status of standard tools like the Borel-Cantelli lemmas in the classical probability literature. There may soon be a time when researchers will start using Talagrand's inequalities without explicit reference in the abstract. 

The theory of concentration inequalities for functions of independent random variables, some of which was described above, has reached a high level of sophistication by now. However, concentration inequalities for functions of dependent random variables are still hard to get, the main tools being logarithmic Sobolev and transportation cost inequalities. One shortcoming of these methods is that explicit constants are very hard or almost impossible to get. We shall discuss these methods in detail in the next chapter.

The main purpose of this thesis is to construct a modification of Stein's method of exchangeable pairs, which is a well-known tool from probability theory, to derive concentration inequalities with explicit constants for functions of dependent random variables. We postpone a discussion of Stein's method until Section \ref{stein}. 

\section{Summary of thesis}
We now give a brief chapter by chapter description of this thesis. The theory is too abstract and requires too much notation to describe in this brief introductory discussion. Instead, we shall present some theorems which were obtained as applications of the theory, for the purpose of enticing the prospective reader to delve deeper.

In Chapter \ref{review}, we shall review the major existing results from the concentration inequality literature, including the early martingale techniques, Talagrand's inequalities, and modern entropy based methods. We shall also describe the basic philosophy of Stein's method in the last section of this chapter.

Chapter \ref{results} is the first part of our theory, which is the more ``abstract'' part. We shall develop our basic results in this chapter, and apply them to work out some simple examples involving dependent variables: In section \ref{curieweiss}, we shall obtain a simple and explicit concentration bound for $m-\tanh(\beta m)$, where $m$ is the magnetization in the Curie-Weiss model of ferromagnetic interaction. By a generalization of this example, we shall derive in section \ref{meanmodels} a broad condition under which the naive mean field equations --- one of the standard tools from physics that is frowned upon by mathematicians (for good reasons) --- is valid. In section \ref{least}, we shall work out an example about estimation of parameters in statistical models of dependent data. In section \ref{perm}, we shall apply our techniques to obtain tail bounds of the correct order for (a) the number of fixed points in a random permutation, and (b) the Spearman's footrule distance between a random permutation and the identity. Finally, in section~\ref{spinglass}, we shall obtain a concentration result about the Sherrington-Kirkpatrick model of spin glasses which holds at all temperatures.

In Chapter \ref{coupling} , we shall develop some advanced tools (Lemma \ref{inv1}, Lemma \ref{invlmm}) so that the theorems from Chapter \ref{results} can be easily applied to more complex problems than the ones worked out in Chapter \ref{results}. Indeed, Lemma \ref{invlmm}, in conjunction with Theorem \ref{hoeffding}, is probably the first general tool which allows one to use couplings to prove concentration inequalities. Couplings are widely used to establish rates of temporal and spatial decay of correlations; so it seems natural that there should be some result which makes them useful for concentration inequalities, too.

An application of these tools gives Theorem~\ref{dobthm}, which allows us to derive concentration inequalities for arbitrary functions in complicated models of dependent random variables, such as Ising type spin models and random proper colorings of finite graphs, under a Dobrushin type condition of weak dependence. We shall now state this theorem, after introducing some required notation. 

Let $\Omega$ be a Polish space and let $f:\Omega^n \ra \rr$ be a function satisfying a Lipschitz condition with respect to a generalized Hamming distance on $\Omega^n$: For all $x,y\in \Omega^n$, 
\begin{equation}\label{flip}
|f(x)-f(y)| \le \sum_{i=1}^n c_i \ii\{x_i \ne y_i\}
\end{equation}
where $c_1,\ldots,c_n$ are some fixed nonnegative constants. This just means that the value of the function does not change by more than $c_i$ if the $i^{\mathrm{th}}$ coordinate is altered. 

Let $X = (X_1,\ldots,X_n)$ be an $\Omega^n$-valued random variable with law $\mu$. For any $x\in \Omega^n$, let $\bar{x}^i$ denote the element of $\Omega^{n-1}$ obtained by omitting the $i^{\mathrm{th}}$ coordinate of $x$. 
For each $i\le n$ and $x\in \Omega^n$, let $\mu_i(\cdot|\bar{x}^i)$ denote the law of $X_i$ given $\bar{X}^i = \bar{x}^i$. 
Finally, let us recall that for a square matrix $A$, the $L^2$ operator norm of $A$ is defined as:
\[
\|A\|_2 := \max_{\|y\|=1} \|Ay\|.
\]
Then we have the following result from section \ref{dobrushin}:
\vskip.1in
\noindent{\bf Theorem \ref{dobthm} } 
{\it Suppose $A = (a_{ij})$ is an $n\times n$ matrix with nonnegative entries and zeros on the diagonal such that for any $i$, and any $x, y \in \Omega^n$,} 
\[
d_{TV}(\mu_i(\cdot|\bar{x}^i), \mu_i(\cdot|\bar{y}^i)) \le \sum_{j=1}^n a_{ij} \ii\{x_j\ne y_j\},
\]
{\it where $d_{TV}$ is the total variation distance on the space of probability measures on $\Omega$.
Suppose $f$ satisfies the generalized Lipschitz condition (\ref{flip}). If $\|A\|_2 < 1$, we have}
\[
\pp\{|f(X)-\ee f(X)| \ge t\} \le 2e^{-(1-\|A\|_2)t^2/\sum_i c_i^2}
\]
{\it for each $t\ge 0$.}
\vskip.1in
\noindent This is possibly the first result which gives a direct connection between the Dobrushin condition \cite{dobrushin70} from statistical physics and concentration inequalities. The log-Sobolev inequalities of Stroock \& Zegarlinski \cite{stroockzeg92a,stroockzeg92b, stroockzeg92c} do not give explicit constants, and it is also not clear whether they extend to systems beyond the lattice.

Some applications of the above theorem, to spin systems and graph colorings, will be worked out in sections \ref{spinsystems} and \ref{graphcolorings}. In particular, we shall work out explicit concentration bounds for the magnetization in the Ising model on an arbitrary graph at high temperature.

In section \ref{group} of the same chapter, we shall derive a tool (Theorem \ref{gpaction}) for obtaining the concentration of measures which are invariant under group actions. This is possibly the first result which gives an explicit connection between rates of convergence to stationarity for random walks on groups and the concentration of Haar measures. Using this result, we shall obtain the following quantitative bound related to random matrices and free probability theory, which again, is probably the first result of its kind. The limiting version of this result is quite well-known; following from a celebrated result of Voiculescu \cite{voiculescu91}, it says, roughly, that the distribution of the eigenvalues of $M+N$, where $M$ and $N$ are two high dimensional hermitian matrices of the same order, is approximately determined by the eigenvalue distributions of $M$ and $N$. 

In the following theorem, the term ``empirical distribution function of $H$'' is a commonly used random matrix jargon, which just means the probability distribution function on $\rr$ which puts mass $1/n$ on each eigenvalue of the matrix $H$.
\vskip.1in
\noindent{\bf Theorem \ref{free} }
{\it Let $\Delta_1$ and $\Delta_2$ be two $n\times n$ real diagonal matrices. Let $U$ and $V$ be independent Haar distributed random elements of $\uu_n$, the group of all unitary matrices of order $n$. Let}
\[
H = U\Delta_1 U^* + V\Delta_2 V^*,
\]
{\it and let $F_H$ be the empirical distribution function of $H$. Then, for every $x\in \rr$,
$\var(F_H(x)) \le \kappa n^{-1}\log n$.
where $\kappa$ is a universal constant \underline{not} depending on $n$, $\Delta_1$, $\Delta_2$ or $x$.
Moreover we also have the concentration inequality}
\[
\pp\{|F_H(x)-\ee (F_H(x))| \ge t\} \le 2\exp\biggl(-\frac{nt^2}{2\kappa\log n}\biggr)
\]
{\it for every $t\ge 0$, where $\kappa$ is the same as in the variance bound.}
\vskip.1in
\noindent 
We think that it will be very hard to derive such a result using available techniques, because there is no independence, and gaussianity is involved in a complicated way, and the standard concentration results for gaussian measures are all with respect to the Euclidean metric, which cannot give a universal bound like the above (bounds, if any, will involve $\Delta_1$ and $\Delta_2$, and will be inefficient, because the function $F_H$ is badly nonsmooth). 
On the other hand, our bound follows quite easily from the theory developed in this dissertation. We do not know, however, if it is of the correct order (even after discounting the $\log n$ factor).
 
Finally, let us mention some of the deficiencies of this dissertation. One shortcoming is in the range of examples. Although we have tried our best to provide as many as we can, it is probably not enough to cover all fields of interest; in particular, we have no examples from empirical process theory. Moreover, the author is not completely happy with the quality of some of the applications. For instance, the result about spin glasses seems to have a wide scope, but the author has been frustrated by attempts to exploit it further. Another unfinished aspect is that we could not find good examples for the ``unbounded differences'' theorems of section \ref{unbdd}.

It is also unfortunate that the technique in Chapter \ref{dobrushin} has only been applied to get Theorems \ref{dobthm} and \ref{gpaction}. We feel that the scope of Lemma \ref{invlmm}, which allows us to get concentration bounds using couplings, extends much beyond that.

Finally, a major incompleteness, which is a weakness of the idea of concentration itself, is the lack of lower bounds. 

We shall attempt to overcome these and other deficiencies in future work.

\chapter{Review of existing literature}\label{review}
This chapter will be devoted to the discussion of the main existing tools for proving concentration bounds. We apologize in advance for unjust omissions, if any. Most of the material in this chapter, except for the very recent developments, is taken from the wonderful monograph by Michel Ledoux \cite{ledoux01}. In fact, in many places we have kept the notation and even the language, intact. 

While we shall attempt to give a comprehensive summary of the main theoretical results, we shall, in general, refrain from discussing applications in this chapter. The main reason being, applications invariably entail some lengthy background, the discussion of which would be a digression from the central theme of this chapter. This is the same reason why we shall not discuss concentration inequalities in empirical process theory. Of course, generous references will be provided.

In the last section of this chapter, we shall discuss the basics of Stein's method.

\section{Hoeffding type inequalities}\label{hoeffineq}
The Azuma-Hoeffding martingale inequality \cite{hoeffding63, azuma67} remained the last word in concentration for a very long time, culminating in the bounded difference inequality, observed by Schechtman \cite{schechtman81} and also by Shamir \& Spencer \cite{shamirspencer87}, and its extensive popularity among combinatorialists and discrete mathematicians following the expository work of McDiarmid \cite{mcdiarmid89}. 

In its most widely used form, the Azuma-Hoeffding inequality goes as follows:
\begin{thm}\label{azuma}
\textup{[Azuma-Hoeffding inequality \cite{hoeffding63, azuma67}]}
Let $\{X_i\}_{1\le i\le n}$ be a martingale difference sequence adapted to some filtration. Suppose $c_1,\ldots,c_n$ are constants such that $|X_i|\le c_i$ almost surely for each $i$. Then for all $t \ge 0$ we have
\[
\pp\{\max_{1\le k\le n} \sum_{i=1}^k X_i \ge t\} \le \exp\biggl(-\frac{t^2}{2\sum_{i=1}^n c_i^2}\biggr).
\]
\end{thm}
Before stating the bounded difference inequality, we must state the most ``accurate'' version of the Hoeffding inequality for independent summands, which was established by Bennett \cite{bennett62}. 
\begin{thm}
\textup{[Bennett's inequality \cite{bennett62}]}
Let $Y_1,\ldots,Y_n$ be independent real-valued random variables bounded in magnitude by some constant $C$. Let $S= \sum_{i=1}^n Y_i$. Then, for every $t \ge 0$,
\[
\pp\{ S\ge \ee(S) + t\} \le \exp\biggl\{-\frac{\sigma^2}{C^2} h\biggl(\frac{Ct}{\sigma^2}\biggr)\biggl\},
\]
where $h(u) = (1+u)\log(1+u) -u$, and $\sigma^2 = \sum_{i=1}^n \ee (Y_i^2)$.
\end{thm}
Originally proved by Bernstein for Bernoulli random variables, this is also sometimes called Bernstein's  inequality.

The bounded difference inequality, stated in its simplest form, is the following powerful corollary of the Azuma-Hoeffding inequality:
\begin{thm}\label{bddiff}
\textup{[Bounded difference inequality \cite{schechtman81, shamirspencer87, mcdiarmid89}]}
Let $Z = f(X_1,\ldots,X_n)$ be a function of independent random variables $X_1,\ldots,X_n$. Let $\xp_i$ be an independent copy of $X_i$, $i=1,\ldots,n$. Suppose $c_1,\ldots,c_n$ are constants such that for each $i$,
\[
|f(X_1,\ldots,X_{i-1}, \xp_i,X_{i+1},\ldots,X_n)-f(X_1,\ldots,X_n)| \le c_i \ \text{  a.s.}
\]
Then, for any $t\ge 0$ we have
\[
\pp\{ Z-\ee(Z)\ge t\} \le \exp\biggl(-\frac{2t^2}{\sum_{i=1}^n c_i^2}\biggr).
\]
\end{thm}
The above inequality, though widely useful, does not in general convey the theoretical essence of the situation, because the bounds on the differences do not reflect the typical size of the differences in many interesting problems. A better result, from the conceptual point of view, is the following inequality, first discovered by Efron and Stein \cite{efronstein81}, which has since come to be known as the Efron-Stein inequality. The present version is due to Steele \cite{steele86}:
\begin{thm}
\textup{[Efron-Stein inequality \cite{efronstein81, steele86}]}
Keeping the notation exactly as in the previous theorem, we have
\[
\var(Z) \le \frac{1}{2}\sum_{i=1}^n \ee\bigl[(f(X_1,\ldots,X_{i-1},\xp_i,X_{i+1},\ldots,X_n) - f(X_1,\ldots,X_n))^2\bigr].
\]
\end{thm}
The Efron-Stein inequality has been used to bound variances in many complicated problems. For a very recent application, one can see the work of Reitzner \cite{reitzner03} on random polytopes. However, the inequality is not guaranteed to give a bound of the correct order. For instance, the actual variance of the number of $k$-gons in an Erd\H{o}s-R\'enyi random graph is less than the Efron-Stein bound by a factor of $k$ (Cf. \cite{blm03}, section 6). Also, it is only a variance bound and gives no useful information about the tail behavior.

All of the above results are based on martingale methods. The limitations of martingale techniques were gradually recognized in the late eighties, and people started looking for alternatives. In recent years, however, there has been some kind of a minor resurgence of interest in extending the old martingale arguments. One of the more successful efforts has been the so-called ``divide and conquer martingale'' method of Kim \& Vu \cite{kimvu04}. We refer to this paper for further references to the current literature around the martingale method.

A significant discovery was recently made by Boucheron, Lugosi \& Massart \cite{blm03}, who established an exponential version of the Efron-Stein inequality. In our opinion, this work has resulted in the completion of the quest for maximum efficiency in concentration via Hoeffding type inequalities. 

The methodology used to prove the exponential Efron-Stein is based on the ``entropy method'' introduced by Ledoux \cite{ledoux01} and Massart \cite{massart00}, and generalized by Boucheron, Lugosi \& Massart in \cite{blm00}. This powerful new  method is based on the modified log Sobolev approach, to be discussed in section \ref{entropy}.

Before stating the exponential Efron-Stein inequality, we have to introduce some notation.

Let $X_1,\ldots,X_n$ be independent random variables taking values in a measurable space $\xx$. Denote by $X^n_1$ the vector of these $n$ random variables. Let $f:\xx^n\ra \rr$ be a measurable map. Let $Z=f(X_1,\ldots,X_n)$. Let $\xp_1,\ldots,\xp_n$ denote independent copies of $X_1,\ldots,X_n$, and let $Z^{(i)} = f(X_1,\ldots,X_{i-1},\xp_i, X_{i+1},\ldots,X_n)$. 

Define the random variables $V_+$ and $V_-$ by
\[
V_+ = \ee\biggl[\sum_{i=1}^n (Z -Z^{(i)})^2\ii_{Z>Z^{(i)}}\biggl| X^n_1\biggr]
\]
and 
\[
V_-=\ee\biggl[\sum_{i=1}^n(Z-Z^{(i)})^2\ii_{Z<Z^{(i)}}\biggl| X^n_1\biggr].
\] 
The following result appears as Theorem 2 in \cite{blm03}:
\begin{thm}
\textup{[Exponential Efron-Stein inequality]}
For all $\theta >0$ and $\lambda \in (0,1/\theta)$, 
\[
\log \ee\bigl[\exp(\lambda(Z - \ee (Z))\bigr] \le \frac{\lambda\theta}{1-\lambda\theta} \log\ee\biggl[\exp\biggl(\frac{\lambda V_+}{\theta}\biggr)\biggr].
\]
On the other hand, we have for all $\theta > 0$ and $\lambda\in (0,1/\theta)$,
\[
\log \ee\bigl[\exp(-\lambda(Z - \ee (Z))\bigr] \le \frac{\lambda\theta}{1-\lambda\theta} \log\ee\biggl[\exp\biggl(\frac{\lambda V_-}{\theta}\biggr)\biggr].
\]
\end{thm}
The paper \cite{blm03} has several nice applications of this result, including applications to empirical processes and subgraph counts. In particular, the authors prove that Talagrand's famous convex distance inequality (to be discussed in the next section) is a corollary of the exponential Efron-Stein inequality. 

The more recent paper \cite{bblm05} on generalized moment inequalities in the same spirit as the exponential Efron-Stein inequality, is also worthy of note.

\section{Talagrand's concentration inequalities}\label{talagrand}
This section is devoted to the discussion of the deep investigation of Michel Talagrand, which climaxed in a series of intricate and powerful results about concentration of measure in product spaces \cite{talagrand95, talagrand96b, talagrand96c}. These have since come to be known as Talagrand's concentration inequalities. Rooted in abstract geometric formulation, Talagrand's techniques have found wide applications in fields ranging from combinatorial optimization (e.g.\ traveling salesman problem in \cite{talagrand95}) to random matrix theory \cite{gz00, alonvu02}.

Before we begin our discussion, let us introduce some basic notation: Throughout, we consider a product probability measure $\mu^n$ on a product space $\xx^n$. For any vector $c = (c_1,\ldots,c_n)\in \rr^n_+$, the generalized Hamming metric $d_c$ is defined on $\xx^n$ as
\begin{equation}\label{lips}
d_c(x,y) := \sum_{i=1}^n c_i \ii\{x_i\ne y_i\}.
\end{equation}
Thus, the bounded difference inequality says that if $f:\xx^n \ra \rr$ is a function such that $|f(x)-f(y)| \le d_c(x,y)$ for all $x$ and $y$, then $\mu^n\{f - \int fd\mu^n \ge t\} \le \exp(-t^2/2\|c\|^2)$, where $\|c\|^2 := \sum_{i=1}^n c_i^2$.

A fundamental weakness of this inequality and its early variants is the following: They only allow us to consider functions $f$ which satisfy the Lipschitz condition for some fixed $c$. It was soon realized that a lot of open questions in measure concentration were about functions which do not satisfy (\ref{lips}), but rather, obey
\begin{equation}\label{lips2}
f(x)-f(y)\le \sum_{i=1}^n c_i(x)\ii\{x_i\ne y_i\}
\end{equation}
or some variant of this, for some vector field $c:\xx^n \ra \rr^n$ with uniformly bounded norm. (The asymmetry in the above expression is often a help rather than hindrance in applications.) It is thus desirable to have a Hoeffding type bound for such functions based on $\sup\{\|c(x)\|^2: x\in \xx^n\}$. Talagrand's famed convex distance inequality is, in essence, a generalization of this idea in the abstract geometric setting. 

Given a set $A \subseteq \xx^n$, let us define the ``convex distance'' of a point $x$ to $A$ as
\[
D_A (x) := \sup_{\|c\|=1} d_c(x,A),
\]
where $d_c(x,A) = \inf_{y\in A} d_c(x,y)$. Note that $D_A(x) = d_{c(x)}(x,A)$ for some $c(x)$ depending on $x$ if the supremum is attained (which is usually the case).

The following version of Talagrand's main result is taken from Ledoux \cite{ledoux01}, Theorem 4.6:
\begin{thm}\label{convexdist}
\textup{[Talagrand's convex distance inequality \cite{talagrand95}]}
For every measurable non-empty subset $A$ of $\xx^n$, and for every product probability measure $\mu^n$ on $\xx^n$, 
\[
\int e^{D_A(x)^2/4} \mu^n(dx) \le \frac{1}{\mu^n(A)}.
\]
In particular, for every $t\ge 0$, 
\[
\mu^n\{D_A \ge t\} \le \frac{1}{\mu^n(A)} e^{-t^2/4}.
\]
\end{thm}
To connect this with concentration problems for specific functions, the usual route is the following: Given a function $f$ with median $m_f$ (that is, $\mu^n\{f\ge m_f\} \ge 1/2$ and $\mu^n\{f\le m_f\} \ge 1/2$), let $A = \{x: f(x)\le m_f\}$. Then $\mu^n(A) \ge 1/2$. Next, find a function $r(t)$ such that $f(x)-m_f > t$ implies $D_A(x) > r(t)$. It is then easy to see that
\[
\mu^n\{f-m_f > t\} \le \mu^n \{D_A > r(t)\} \le 2e^{-r(t)^2/4}.
\]
The abstract formulation allows us to go beyond Lipschitz functions, in the following way (which is also the usual method of applying the convex distance inequality): Given a function $f$ with median $m_f$, we have to find a function $r(t)$ having the following property: Whenever $f(x) \ge m_f + t$ and $f(y) \le m_f$, there exists some $c = c(x)\in \rr^n$ such that $\|c\|=1$ and $\sum_{i=1}^n c_i(x)\ii\{x_i\ne y_i\} > r(t)$. In the particular case when $r(t) = t$, we get the condition $f(x)-f(y) \le \sum_{i=1}^n c_i(x)\ii\{x_i\ne y_i\}$, which was stated as the initial motivation for this discussion. However, we may have $r(t) \ne t$ also, as demonstrated in Talagrand's bound for the concentration of the longest increasing subsequence in a random permutation (\cite{talagrand95}, section 7.1). 

As is usual in geometric measure concentration (to be discussed in section \ref{geometric}), Talagrand's inequality gives concentration around the median rather than the mean, but that is a minor inconvenience. 

A very important and striking consequence of Theorem \ref{convexdist} is the following:
\begin{cor}\textup{[Talagrand \cite{talagrand95}]}
For every product probability $\mu^n$ on $[0,1]^n$, every convex $1$-Lipschitz function $f$ on $\rr^n$, and every $t\ge 0$,
\[
\mu^n\{|f-m_f| \ge t\} \le 4e^{-t^2/4},
\]
where $m_f$ is the median of $f$ for $\mu^n$.
\end{cor}
The applications of Talagrand's convex distance inequality are too diverse to summarize in a few paragraphs. Many of the most striking examples were worked out by Talagrand himself in his landmark paper \cite{talagrand95}, including applications to the traveling salesman problem and first passage percolation. Applications to the concentration of spectral measures of random matrices were made by Guionnet \& Zeitouni \cite{gz00}. Concentration of individual eigenvalues of gaussian random matrices was established using Theorem \ref{convexdist} in a short but remarkable paper by Alon, Krivelevich \& Vu \cite{alonvu02}. 

Another important result of Talagrand is the ``control by several points'' method, which we shall not discuss here. For a description of this technique, one can look at Talagrand's original paper \cite{talagrand95} or section 4.3 in Ledoux' book \cite{ledoux01}. Talagrand used this method to prove his celebrated Bernstein type bound for empirical processes. This result is presented in a nice way as Theorem 7.4 in \cite{ledoux01}.

\section{Logarithmic Sobolev inequalities}
Logarithmic Sobolev inequalities were introduced by Gross \cite{gross76} as the infinitesimal version of hypercontractivity in quantum field theory. They soon became a tool of fundamental importance, and the last three decades have seen vigorous activity surrounding the theory and applications of log-Sobolev inequalities. It is not our purpose to go deeply into all that; we shall be only concerned with their relevance in measure concentration. For other applications in probability and statistical mechanics, one can look at the lecture notes by Guionnet and Zegarlinski \cite{gz03} which concentrate on physical applications, and the paper by Diaconis and Saloff-Coste \cite{ds96} on applications of log-Sobolev inequalities to finite Markov chains. An excellent survey of known mathematical results can be found in \cite{ane00}.

Suppose $X_0,X_1,\ldots,$ is a stationary reversible Markov chain on some space $\xx$. The {\it Dirichlet form} corresponding to this chain (or more appropriately, corresponding to the associated kernel) is a functional $\mathcal{E}$ defined on the space of all pairs $f,g$ of maps from $\xx$ into $\rr$ which satisfy $\ee(f(X_0)^2)< \infty$ and $\ee(g(X_0)^2)< \infty$. It is defined as
\[
\mathcal{E}(f,g) := \frac{1}{2}\ee\bigl[(f(X_1)-f(X_0)(g(X_1)-g(X_0))\bigr].
\]
The kernel is said to satisfy a {\it logarithmic Sobolev inequality} with constant $c$ if for all $f$ such that $\ee(f(X_0)^2)< \infty$, we have
\[
\ee\biggl(f(X_0)^2 \log \frac{f(X_0)^2}{\ee f(X_0)^2}\biggr) \le 2c \mathcal{E}(f,f).
\]
Often, probability measures have natural reversible kernels associated with them. In such cases, instead of saying that the kernel satisfies a log-Sobolev inequality, we say that the probability distribution satisfies such an inequality. For instance, a distribution on $\rr^n$ which has density $\rho(x)$ with respect to Lebesgue measure, is, under appropriate conditions, the stationary distribution of the Langevin diffusion process $(X_t)_{t\ge 0}$ with constant volatility matrix $\sqrt{2}I$ and drift $\nabla \log \rho(x)$. Instead of $X_0$ and $X_1$ we now take $X_0$ and $X_h$, where $h > 0$, but divide the right hand side by $h$ in the definition of the log-Sobolev inequality. For the Langevin diffusion under suitable conditions, if we take take $h\downarrow 0$, then $h^{-1}\mathcal{E}(f,f) \ra \ee(\|\nabla f(X_0)\|^2)$. This observation motivates the following definition: A probability measure $\mu$ on $\rr^n$ is said to satisfy a logarithmic Sobolev inequality with constant $c$ if for all locally Lipschitz~$f$, 
\[
\int f(x)^2 \log \biggl(\frac{f(x)^2}{\int f(u)^2d\mu(u)}\biggr) d\mu(x) \le 2c \int \|\nabla f(x)\|^2 d\mu(x).
\]
To most people, this is the usual definition of a log-Sobolev inequality, but the earlier form is useful for discrete problems. The left hand side is usually called the entropy of $f^2$ with respect to $\mu$, and is denoted by $\mathrm{Ent}_\mu(f^2)$. 

The single most important property of log-Sobolev inequalities is perhaps the tensorizing property: If $\mu_1, \ldots,\mu_n$ are probability measures on $\rr$ satisfying log-Sobolev inequalities with constants $c_1,\ldots,c_n$, then the product measure $\mu_1\otimes\cdots\otimes \mu_n$ on $\rr^n$ satisfies log-Sobolev inequality with constant $\max_i c_i$.

The connection between logarithmic Sobolev inequalities and concentration was made in an unpublished but now famous argument of I.\ Herbst. The following theorem, which summarizes the end result, is taken from Ledoux \cite{ledoux01}, Theorem 5.3:
\begin{thm}
\textup{[Herbst's lemma]}
Let $\mu$ be a probability measure on $\rr^n$ satisfying a log-Sobolev inequality with constant $c$.
Then, every $1$-Lipschitz function $f:\xx \ra \rr$ is integrable and for every $t\ge 0$,
$\mu\{f\ge \int fd\mu + t\} \le e^{-t^2/2c}$.
\end{thm}
Ledoux also presents a discrete version of the above result; the following occurs as Theorem 5.17 in \cite{ledoux01}:
\begin{thm}
Let $\mu$ be the stationary distribution of a reversible Markov chain $\{X_k\}_{k\ge 0}$ on a countable set $\xx$, which satisfies a logarithmic Sobolev inequality with constant $c$. Then, any $f:\xx \ra \rr$ which satisfies
\[
\frac{1}{2}\sup_{x\in \xx}\ee((f(X_1)-f(X_0))^2|X_0=x) \le 1,
\]
also satisfies $\mu\{f\ge \int fd\mu + t\} \le e^{-t^2/4c}$ for any $t\ge 0$.
\end{thm}
One class of measures for which explicit log-Sobolev constants (and hence, concentration inequalities) are easily available are high dimensional probability measures with strictly log-concave density with respect to Lebesgue measure; this is a consequence of an important work of Bakry and \'Emery \cite{bakryemery85}.
\begin{thm}
\textup{[Bakry-\'Emery criterion]}
Suppose $\mu$ is a measure on $\rr^n$ with density $e^{-U(x)}$ with respect to Lebesgue measure, where $\mathrm{Hess}\; U(x) \ge a I$ for some fixed constant $a>0$, for all $x\in \rr^n$. Then $\mu$ satisfies a log-Sobolev inequality with constant $1/a$. 
\end{thm}
The result can be supplemented by the observation that if $\mu$ satisfies a logarithmic Sobolev inequality with constant $c$, then the measure defined by $d\nu = Z^{-1}e^V d\mu$ (where $Z$ is the normalizing constant) satisfies a log-Sobolev inequality with constant $ce^{4\|V\|_\infty}$. This is the famous perturbation argument of Holley and Stroock \cite{holleystroock87}. However, this is not a very useful result in high dimensions, because $\|V\|_\infty$ blows up as the dimension increases.

It is to be noted that we do not necessarily have to go through the Herbst argument and Bakry-\'Emery's result to get the concentration inequality for strictly log-concave measures. Ledoux \cite{ledoux01},  pp.\ 39--41, has a direct proof.

However, in other cases, like spin systems on a lattice, there is no direct argument for concentration inequalities, and the route through the Herbst lemma must be used. In their important work \cite{stroockzeg92a, stroockzeg92b, stroockzeg92c}, Stroock and Zegarlinski established that the Glauber dynamics associated with Ising and other spin models satisfy log-Sobolev inequalities at high temperature, wherever Dobrushin's condition of weak dependence holds. (See Georgii \cite{georgii88}, Chapter 8, for the definition and examples for Dobrushin's condition. For a readable account of the work of Stroock and Zegarlinski, see the lecture notes \cite{gz03}.) 

This implies, by Herbst's argument, that those spin systems satisfy concentration properties at high temperature. However, since explicit log-Sobolev constants are not known, it is not possible to get explicit concentration bounds along those lines. 

This brings us to the point where we can state one of our achievements in this dissertation, which we already stated once in the Introduction: In section \ref{dobrushin}, we shall prove a concentration inequality with explicit constants for systems which satisfy Dobrushin's condition. We shall apply our result to get explicit concentration bounds for generalized Lipschitz functions of spins in Ising type models and uniformly chosen proper $k$-colorings of graphs with maximum degree $< k/2$.

For further information and references about the use of logarithmic Sobolev inequalities in the field of measure concentration, the reader is encouraged to look at the comprehensive survey \cite{ledoux99}. For further details about their applications to spin systems, as well as some easy proofs of known results, one can look at \cite{ledoux01a}.

\section{Other entropy based methods}\label{entropy}
Ever since the advent of Herbst's argument, entropy-based methods have played an increasingly important role in the concentration literature. In fact, at the time of writing this thesis, they are probably the most active tool of research in concentration inequalities. We have already discussed the logarithmic Sobolev approach to concentration as developed by Ledoux \cite{ledoux01} and Massart \cite{massart00}. In this section, we briefly describe two other important information theoretic methods, namely, the modified log-Sobolev inequalities of Ledoux and the transportation cost inequalites of Marton.
\vskip.1in

\noindent {\bf Modified log-Sobolev inequalities.} A probability measure $\mu$ on $\rr^n$ is said to satisfy a {\it modified logarithmic Sobolev inequality} if there is a function $\beta(\rho) \ge 0$ on $\rr$, such that, whenever $\|\nabla f\|_\infty \le \rho$,
\[
\mathrm{Ent}_\mu(e^f) \le \beta(\rho)\int \|\nabla f\|^2 e^f d\mu
\]
for all locally absolutely continuous $f$ such that $\int e^f d\mu <\infty$.

Modified log-Sobolev inequalities were introduced by Bobkov and Ledoux \cite{bobkovledoux97} to provide an easier alternative to Talagrand's method for proving exponential concentration. They showed that like log-Sobolev inequalities, modified log-Sobolev inequalities also tensorize in a certain sense. An important consequence of this is the following concentration result for product of measures which satisfy Poincar\'e inequalities (the exponential distribution being the prototype for such measures):
\begin{thm}
Suppose $\mu$ is a probability measure on $\rr$ satisfying a Poincar\'e inequality with constant $c$; that is, for any locally absolutely continuous $f$, $\var_\mu(f) \le c\int |f^\prime|^2 d\mu$. Then, any function $F:\rr^n \ra \rr$ satisfying
\[
\sum_{i=1}^n \biggl|\frac{\partial F}{\partial x_i}\biggr|^2 \le a^2 \ \ \text{ and } \ \ \max_{1\le i\le n}\biggl|\frac{\partial F}{\partial x_i}\biggr| \le b
\]
$\mu^n$-almost everywhere is integrable with respect to $\mu^n$, and for every $t\ge 0$,
\[
\mu^n \{ F \ge \int Fd\mu^n + t\} \le \exp\biggl(-\frac{1}{K} \min\biggl(\frac{t}{b}, \frac{t^2}{a^2}\biggr)\biggr).
\]
\end{thm}
This occurs as Corollary 5.15 in \cite{ledoux01}. For more on modified log-Sobolev inequalities, one can look at the survey \cite{ledoux99}. Some recent developments using the modified log-Sobolev approach (alternatively, the ``entropy method'') have already been discussed in section \ref{hoeffineq}. 
\vskip.1in

\noindent{\bf Transportation cost inequalities.}
Transportation cost inequalities were introduced by Marton \cite{marton96} as a version of measure concentration which works by investigating distances between measures. 

Let's begin with some familiar definitions.
The informational divergence (or Kullback-Leibler divergence, or relative entropy) of a measure $\nu$ with respect to another measure $\mu$ on the same space is defined as:
\[
D(\nu\|\mu) := \int \log \biggl(\frac{d\nu}{d\mu} \biggr)d\nu.
\]
If these measures are probability measures defined on a Polish space $(\xx,d)$, then the $L^1$ and $L^2$ {\it Wasserstein distances} between $\nu$ and $\mu$ are defined as 
\[
W_1(\nu, \mu) := \inf_\pi\ee_\pi d(Y,X) \ \text{ and } \ W_2(\nu, \mu) := \inf_\pi \bigl[\ee_\pi d(Y,X)^2\bigr]^{1/2}
\]
where $Y$ and $X$ are random variables distributed according to the laws $\nu$ and $\mu$, and the infimum is taken over all distributions $\pi$ on $\xx^2$ that have $\nu$ and $\mu$ as marginals.

A probability measure $\mu$ on $\xx$ is said to satisfy a {\it transportation cost inequality} with constant $c$, if 
\[
W_1(\nu, \mu) \le \sqrt{2c D(\nu\|\mu)} \ \ \text{ for all probability measures $\nu$ on $\xx$.}
\]
Similarly, $\mu$ satisfies a {\it quadratic transportation cost inequality} (or, as Marton prefers to call it, a distance-divergence inequality) if
\[
W_2(\nu,\mu) \le \sqrt{2c D(\nu\|\mu)} \ \ \text{ for all probability measures $\nu$ on $\xx$.}
\]
It is not difficult to deduce from either of the above conditions that for any two measurable sets $A, B\subseteq \xx$, 
\[
d(A,B) \le \sqrt{2c\log \frac{1}{\mu(A)}} + \sqrt{2c \log \frac{1}{\mu(B)}},
\]
where $d(A,B) := \inf\{d(x,y): x\in A, \ y\in B\}$. Putting $B = \{x: d(x,A) \ge t\}$ in the above, we obtain an abstract measure concentration inequality,
\[
1-\mu(A_t)\le e^{-t^2/8c},
\]
where $A_t := \{x:d(x,A)<t\}$. This interpretation of measure concentration was put forward by Milman in the late seventies. In particular, it easily implies concentration inequalities for Lipschitz functions. This will be 
discussed in some detail in section~\ref{geometric}. 

The above was the original line of argument by Marton connecting transportation cost inequalities with concentration of measure.

Quadratic transportation cost inequalities have been more successful for Euclidean spaces, mainly because they have a tensorizing property similar to log-Sobolev inequalities, which is not true for transportation cost inequalities based on the $W_1$ metric. 

There seems to be a close connection between log-Sobolev and quadratic transportation cost inequalities. Indeed, 
Otto \& Villani \cite{ottovillani00} proved that any measure which satisfies a logarithmic Sobolev inequality also satisfies a quadratic transportation cost inequality. The converse, however, is still an open question. For more information on mass transportation, we refer to the recent treatise \cite{villani03}.

One significant success of the transportation cost method has been its application to concentration for Markov chains. Originally proved by Marton \cite{marton96b} for contracting Markov chains, the results were later extended to Doeblin recurrent chains and $\Phi$-mixing processes by Samson \cite{samson00} and simultaneously by Marton in a manuscript which was not sent for publication. We shall now briefly describe the result for contracting Markov chains.

Let $\mu$ be a probability measure on $\xx = \xx_1\times \cdots \times \xx_n$ induced by a Markov chain with successive transition kernels $\Pi_i$, $i=1,\ldots,n$; that is,
\[
d\mu(x_1,\ldots,x_n) =  \Pi_1(dx_1)\Pi_2(x_1, dx_2) \cdots \Pi_n(x_{n-1}, dx_n).
\]
Assume that the Markov chain is {\it contracting}, that is, there exists $\rho < 1$ such that for each $2\le i\le n$ and $x,y \in \xx_i$,
\[
d_{TV}(\Pi_i(x,\cdot), \Pi_i(y, \cdot)) \le \rho.
\]
Then, the following analog of Theorem \ref{convexdist} holds:
\begin{thm}
\textup{[Marton \cite{marton96b}]}
For any measurable nonempty subset $A$ of $\xx$, 
\[
\int e^{(1-\rho)^2 D_A(x)^2 /4 } d\mu(x) \le \frac{1}{\mu(A)}.
\]
\end{thm}
As in Theorem \ref{convexdist}, this implies 
\[
\mu\{D_A \ge t\} \le \frac{1}{\mu(A)} e^{-(1-\rho)^2t^2/4}.
\]
Consequently, if $f:\xx\ra \rr$ is a map with median $m_f$ and $r(t)$ is a function such that whenever $f(x) > m_f +t$ and $f(y)\le m_f$, there exists a vector $c(x)\in \rr^n$ (not depending on $y$) with norm $1$ such that $f(x)-f(y)\le \sum_{i=1}^n c_i(x)\ii\{x_i\ne y_i\}$, then
\[
\mu\{f-m_f \ge t\} \le 2e^{-(1-\rho)^2r(t)^2/4}.
\]
At this point, we should also mention that there is some recent work of Houdr\'e and Tetali \cite{houdretetali01, houdretetali04} on concentration for Markov chains using a different approach, which we shall not discuss here.

Since the pioneering work of Marton, a significant number of important contributions were made by various authors. Talagrand \cite{talagrand96a} proved, among other things, that the gaussian measure on $\rr^n$ satisfies a transportation cost inequality. 
Dembo \cite{dembo97} used the transportation cost and other information theoretic methods to rederive most of Talagrand's abstract inequalities.
The papers by Bobkov \& G\"otze \cite{bobkovgotze99}, and Bobkov, Gentil \& Ledoux \cite{bgl01} are also important. We refer to Ledoux \cite{ledoux01}, Chapter~6, for details.

More recently, Marton \cite{marton03, marton04} has worked on developing transportation cost inequalities for highly dependent systems of random variables that usually occur in statistical physics. In such models, the only tractable objects are the conditional distributions of small subcollections given the rest. The results in \cite{marton03, marton04} hold under Dobrushin-Shlosman type contractivity conditions. Dobrushin type conditions of weak dependence, which originated in statistical physics, will be discussed in section~\ref{dobrushin}, where we shall also present a method based on exchangeable pairs, for directly obtaining concentration inequalities with explicit constants under similar situations.

\section{Geometric measure concentration}\label{geometric}
This section contains a brief discussion of a closely related topic, called ``concentration of measure'' in geometric analytic circles. It is, in fact, a generalization of the idea of concentration inequalities, which has been a topic of interest in geometric functional analysis and convex geometry in the last three decades. The idea rests on the following basic observation: In high dimensional spaces, certain probability measures, including products of well behaved one-dimensional measures, exhibit the ``concentration property'', which means that any set which has measure $\ge 1/2$, ``engulfs'' most of the space when slightly expanded. The idea was pioneered by Paul L\'{e}vy \cite{levy51} who observed this phenomenon for the uniform measure on high dimensional spheres; but the connection of his work with concentration inequalities remained obscure for many years until the papers by Amir \& Milman \cite{amirmilman80} and Gromov \& Milman \cite{gromovmilman83} revived mathematical interest in this very fundamental feature of high dimensional measure spaces. 

Let us now formalize the notion expressed vaguely in the last paragraph. Given a Polish space $(\xx,d)$ and a probability measure $\mu$ on $\xx$, the ``concentration function'' $\alpha_{(\xx,d,\mu)}: [0,\infty) \ra [0,1]$ is defined as:
\[
\alpha_{(\xx,d,\mu)}(t) := \sup\{1-\mu(A_t): A\subseteq \xx, \ \mu(A) \ge 1/2\},
\] 
where $A_t := \{x\in \xx: d(x,A) < t\}$, with $d(x,A):= \inf\{d(x,y):y\in A\}$, as usual. This definition is due to Amir and Milman \cite{amirmilman80}. A measure is ``concentrated'' when $\alpha$ decreases rapidly as $t$ grows. The connection with concentration inequalities for Lipschitz maps comes through the following easy observation: If $f:\xx\ra\rr$ is a Lipschitz function with Lipschitz constant $\|f\|_{\mathrm{Lip}}$, and $m_f$ is a median of $f$ (that is, $\mu\{f\ge m_f\} \ge 1/2$ and $\mu\{f\le m_f\} \ge 1/2$), then
\[
\mu\{|f-m_f|\ge t\} \le  2\alpha_{(\xx,d,\mu)}(t/\|f\|_{\mathrm{Lip}}).
\]
The advantage of using concentration functions is that they also apply to functions which are not necessarily Lipschitz, and indeed, not well-behaved in any classical sense, as is demonstrated for instance by Talagrand's treatment of the concentration problem for the longest increasing subsequence of a random permutation (\cite{talagrand95}, section 7.1).

A major success of the geometric approach to measure concentration came with the simple proof of the famous Dvoretzky theorem by Milman \cite{milman71}. To state Dvoretzky's theorem, we first need a definition: A Banach space $(E, \|\cdot\|)$ is said to contain a subspace $(1+\varepsilon)$-isomorphic to $\rr^k$ if there are vectors $v_1,\ldots,v_k$ in $E$ such that for all $t = (t_1,\ldots,t_k)\in \rr^k$, 
\[
(1-\varepsilon) \|t\| \le \biggl\|\sum_{i=1}^k t_i v_i\biggr\| \le (1+\varepsilon) \|t\|.
\] 
Dvoretzky's theorem is the following very important result:
\begin{thm}
\textup{[Dvoretzky's Theorem \cite{dvoretzky61}]}
For each $\varepsilon >0$ there exists $\eta(\varepsilon)>0$ such that every Banach space $E$ of dimension $n$ contains a subspace $(1+\varepsilon)$-isomorphic to $\rr^k$ where $k = [\eta(\varepsilon)\log n]$.
\end{thm}
This result was at the center of the vigorous activity around the local theory of Banach spaces in the decades 1970-90. 
Milman's treatment and subsequent developments show that Dvoretzky's theorem can be viewed as a manifestation of the measure concentration phenomenon in a certain sense. We refer the interested reader to the lecture notes \cite{milmanschechtman86} for extensive details. 

For further details about the functional analytic aspect (which is not relevant to this dissertation), we refer to Milman \cite{milman88} and Chapters 2 and 3 of Ledoux \cite{ledoux01}.

\section{Concentration on groups}\label{gpreview}
A group $G$ is called ``topological'' if it is endowed with a topology which makes the group operations continuous. It is a classical result that on any compact topological group $G$, there exists a unique probability measure $\mu$ which is left and right invariant; that is, if $X$ is a $G$-valued random variable following the law $\mu$, then $xX$ and $Xx$ have the same law as $X$ for every $x\in G$. It also follows that $X^{-1}$ has the same law. This measure $\mu$ is called the ``Haar measure on $G$'' (or the ``normalized Haar measure'' in some texts, but Haar measures will always be normalized for us). For a self-contained proof of the existence and uniqueness of Haar measures on compact groups, we refer to Rudin \cite{rudin73}, Theorem 5.14. 

There is not much literature on the concentration of Haar measures. Before stating whatever little we could find, let us clarify that the works of Gromov \& Milman \cite{gromovmilman83} and Pestov \cite{pestov00a, pestov00b} are about a different kind of ``concentration'' on groups, which is not related to our investigation.

One early result is due to Maurey \cite{maurey79}, who investigated the Haar measure on the group $S_n$ of all permutations of $n$ elements.
\begin{thm}\label{maurey}
\textup{[Maurey \cite{maurey79}]}
Let $\mu$ denote the uniform probability measure on $S_n$. Let $d_n$ denote the normalized Hamming metric on $S_n$: $d_n(\sigma,\tau) = \frac{1}{n}\sum_{i=1}^n \ii\{\sigma(i)\ne \tau(i)\}$. Then, for any $A\subseteq S_n$ such that $\mu(A)\ge 1/2$, and any $t\ge 0$, we have
$\mu(A_t) \ge 1 - 2e^{-nt^2/64}$, where $A_t = \{\sigma\in S_n: d_n(\sigma,A)< t\}$, as usual. Consequently, for any function $f:S_n \ra \rr$ which satisfies $|f(\sigma)-f(\tau)|\le d_n(\sigma,\tau)$ for all $\sigma,\tau\in S_n$, we have
\[
\mu\{|f-m_f|\ge t\} \le 2e^{-nt^2/64} \ \ \text{for all } t\ge 0,
\]
where $m_f$ is the median of $f$ with respect to $\mu$.
\end{thm}
Maurey's result was generalized in the lecture notes of Milman and Schechtman (\cite{milmanschechtman86}, Theorem 7.12) using the classical martingale argument to give the following theorem:
\begin{thm}
\textup{[Milman \& Schechtman \cite{milmanschechtman86}]} 
Let $G$ be a group, compact with respect to a translation invariant metric $d$. Let $G = G_0\supseteq G_1\supseteq \cdots \supseteq G_n =\{1\}$ be a decreasing sequence of closed subgroups of $G$. Let $a_k$ be the diameter of $G_{k-1}/G_k$,  $k=1,\ldots,n$. Let $\mu$ be the Haar measure and let $f:G\ra \rr$ be a function satisfying $|f(x)-f(y)| \le d(x,y)$ for all $x,y\in G$. Then for all $t\ge 0$, 
\[
\mu\{|f - {\textstyle{\int}} fd\mu|\ge t\} \le 2\exp\biggl(-\frac{t^2}{4\sum_{k=1}^n a_k^2}\biggr).
\]
Moreover, if $A\subseteq G$ is such that $\mu(A)\ge 1/2$, then for all $t > 0$, 
\[
\mu(A_t) \ge 1 - 2\exp\biggl(-\frac{t^2}{16\sum_{k=1}^n a_k^2}\biggr),
\]
where $A_t := \{x\in G: d(x,A)< t\}$. 
\end{thm}
Maurey's theorem may be easily recovered by considering the tower $S_n \supseteq S_{n-1}\supseteq \cdots \supseteq S_0 = \{1\}$. 
Although the theorem looks pretty general, it is not very clear what sorts of examples it might cover. In fact, we could not find in the literature any interesting application of this result other than the original application to rederive Maurey's theorem. Still, we decided to include it in the review because it has the looks of a powerful general result which has not yet been fully exploited. 

Next, let us discuss a result of Talagrand about the concentration of the Haar measure on $S_n$ from his seminal work \cite{talagrand95}.
For each $A \subseteq S_n$ and $\sigma \in S_n$, define
\[
U_A(\sigma) := \{s\in\{0,1\}^n: \text{ for some } \tau \in A, \ \ii\{\tau(i)\ne \sigma(i)\} \le s_i \text{ for each } i\}.
\]
Let $V_A(\sigma)$ be the convex hull of $U_A(\sigma)$ in $[0,1]^n$, and let
\[
f(A,\sigma) := \inf\{\|s\|^2: s\in V_A(\sigma)\}.
\]
Then, Talagrand (\cite{talagrand95}, Theorem 5.1) has the following result:
\begin{thm}
\textup{[Talagrand \cite{talagrand95}]} 
For every $A \subseteq S_n$, we have
\[
\int_{S_n} e^{f(A,\sigma)/16} d\mu(\sigma) \le \frac{1}{\mu(A)},
\]
where $\mu$ is the uniform probability on $S_n$. 
\end{thm}
As the reader might have observed, it is not clear what is going on. Talagrand does not provide any example for the above theorem in his paper, although applications should be quite similar to those for product measures, because the setup is more or less the same. A workable corollary of this theorem was devised by McDiarmid \cite{mcdiarmid02}. 

Finally, let us state a result from Gromov \& Milman \cite{gromovmilman83} about the concentration of the Haar measure on $SO_n$ --- the group of $n\times n$ orthgonal matrices with determinant~$1$. 
\begin{thm}
\textup{[Gromov \& Milman \cite{gromovmilman83}]}
Consider the group $SO_n$ of $n\times n$ orthogonal matrices with determinant $1$. Let $d$ be the Hilbert-Schmidt metric on $SO_n$, defined as $d(A,B) = [\sum_{i,j}(a_{ij}-b_{ij})^2]^{1/2}$, where $A=(a_{ij})$ and $B=(b_{ij})$ are any two elements of $SO_n$. Let $\mu$ be the Haar measure on $SO_n$. Then, for any set $A\subseteq SO_n$ such that $\mu(A)\ge 1/2$, and for any $t\ge 0$, we have
\[
\mu(A_t) \ge 1 - \sqrt{\frac{\pi}{8}}e^{-(n-1)t^2/8}.
\]
In particular, if $f:SO_n \ra \rr$ is a function such that $|f(A)-f(B)|\le d(A,B)$ for all $A,B$, then $\mu\{|f-m_f|\ge t\} \le \sqrt{\frac{\pi}{2}} e^{-(n-1)t^2/8}$ for each $t\ge 0$, where $m_f$ is the median of $f$ with respect to $\mu$.
\end{thm}
This result was used by Voiculescu \cite{voiculescu91} in his seminal work connecting free probability theory to random matrices. However, it is not of much use for a quantitative analysis of the eigenvalue spectrum. In section \ref{freeprob}, we shall develop a new method, using the {\it rank metric}, defined as $d(M,N) = \mathrm{rank}(M-N)$, for analyzing the concentration of unitary matrices which will be suitable for studying  
spectral distributions of related random matrices which arise in free probability theory.

\section{Stein's method of exchangeable pairs}\label{stein}
Stein's method was introduced by Charles Stein \cite{stein72} in the context of normal approximation for sums of dependent random variables. Stein's version of his method, best described as the ``method of exchangeable pairs'', attained maturity in his later work \cite{stein86}. A reasonably large literature has developed around the subject, but it has almost exclusively developed as a method of proving distributional convergence with error bounds. Stein's attempts at getting large deviations in \cite{stein86} did not, unfortunately, prove fruitful. The main purpose of this dissertation is to outline a simple way of deriving concentration inequalities using the method of exchangeable pairs, and applying it to problems involving dependent variables.

We shall now briefly describe Stein's method for distributional approximation. Suppose we want to show that a random variable $X$ taking values in some space $\xx$ has approximately the same distribution as some other random variable $Z$. The procedure involves four steps:
\begin{enumerate}
\item Identify a ``characterizing operator'' $T$ for $Z$, which has the defining property that for any function $g$ belonging to a fixed large class of functions, $\ee Tg (Z) = 0$. For instance, if $\xx = \rr$ and $Z$ is a standard gaussian random variable, then $Tg(x) := g^\prime(x)-xg(x)$ is a characterizing operator, acting on all locally absolutely continuous $g$ with moderate growth at infinity.  
\item Construct a random variable $\xp$ such that $(X,\xp)$ is an exchangeable pair.
\item Find an operator $\alpha$ such that for any suitable $h:\xx \ra \rr$, $\alpha h$ is an antisymmetric function (that is, $\alpha h(x,y)\equiv -\alpha h(y,x)$) and
\[
|\ee(\alpha h(X,\xp)|X = x) - Th(x)| \le \varepsilon_h,
\]
 where $\varepsilon_h$
is a small error depending only on $h$. 
\item Take a function $g$ and find $h$ such that $Th(x) = g(x) - \ee g(Z)$. By antisymmetry of $\alpha h$, it follows that $\ee(\alpha h(X,\xp))=0$. Combining with the previous step, we get $|\ee g(X) - \ee g(Z)|\le \varepsilon_h$. 
\end{enumerate}
Note that the operator $S$ defined as
\[
Sh(x) := \ee(\alpha h(X,\xp)|X=x)
\]
is a characterizing operator for the distribution of $X$. Thus, the basic principle of Stein's method is to prove the closeness of the distributions of $X$ and $Z$ by showing that the characterizing operators are close.

Stated differently, this can be viewed as the study of stationary distributions of reversible Markov chains using their generators, but differences exist. The essential difference is that Stein's method involves only the pair $(X,\xp)$ instead of the whole chain. The restriction of attention makes a lot of conceptual and practical difference.

There are other variants of Stein's method, most notably the dependency graph approach popularised by Arratia, Goldstein and Gordon \cite{agg90}, and the size-biased and zero-biased couplings of Barbour, Holst \& Janson \cite{bhj92}, but we shall not discuss those. 

Stein's method has been successfully used to prove convergence to gaussian and Poisson distributions in various situations involving dependent random variables. (Poisson approximation by Stein's method was introduced by Chen \cite{chen75} and became popular after the publication of \cite{agg89, agg90}.) It is not our purpose here to go deeply into the regular versions of Stein's method. For further references and exposition, we refer to the recent monograph \cite{stein04}. For applications of the method of exchangeable pairs and other versions of Stein's method to Poisson approximation, one can look at the survey paper by Chatterjee, Diaconis \& Meckes \cite{cdm05}.

\chapter{Theory and examples: Part I}\label{results}
In this chapter, we shall concentrate on the ``abstract'' part of  our theory. The more directly applicable part will be presented in the next chapter. 

Each theorem will be followed by a demonstrative example, which will always be the easiest application, involving sums of independent random variables. Let us now briefly mention the examples that we shall work out in this chapter.

In section \ref{curieweiss}, we shall obtain a simple and explicit concentration bound for $m-\tanh(\beta m)$, where $m$ is the magnetization in the Curie-Weiss model of ferromagnetic interaction (which will be defined in that section for the reader who is not familiar with ferromagnetic models). By a generalization of this example, we shall derive in section \ref{meanmodels} a broad condition under which a mathematically dubious but often useful technique from physics, known as {\it the naive mean field equations}, is valid. In section \ref{least}, we shall work out an example about estimation in models of dependent data (like Markov Random Field models) from mathematical statistics. In section \ref{perm}, we shall apply our techniques to obtain tail bounds of the correct order for (a) the number of fixed points in a random permutation, and (b) the Spearman's footrule distance between a random permutation and the identity. Finally, in section \ref{spinglass}, we shall obtain a temperature-free concentration result about the Sherrington-Kirkpatrick model of spin glasses.

Throughout the remainder of this thesis, unless otherwise specified, the following notation will remain fixed and understood without mention:
\begin{itemize}
\item $\xx$ is a Polish space and $X$ is a random variable taking values in $\xx$.
\item $f:\xx \ra \rr$ is a measurable map. The object of interest is $f(X)$, and we assume, without loss of generality, that $\ee f(X) = 0$. (We shall abandon this last assumption on certain occasions.)
\item $\xp$ is another random variable defined on the same probability space as $X$, such that $(X,\xp)$ is an exchangeable pair.
\item $F:\xx^2 \ra \rr$ is a measurable map such that $F(X,\xp) = - F(\xp,X)$ almost surely (i.e.\ $F$ is {\it antisymmetric}) and $\ee(F(X,\xp)|X) = f(X)$. For now, we shall assume that both $F$ and $f$ are known. In Chapter \ref{coupling}, we shall investigate a method for obtaining $F$ from a given $f$ using couplings. 
\item We associate a function $v$ with $f$, which will serve as a seemingly crude but useful ``stochastic bound on the size of $f(X)^2$'':
\begin{equation}\label{ddef}
v(x):= \frac{1}{2}\ee\bigl(|(f(X)-f(\xp))F(X,\xp)|\bigl| X = x\bigr).
\end{equation}
\end{itemize}
The general principle of Stein's method of exchangeable pairs is to study the behavior of $f(X)$ using $F(X,\xp)$. Exchangeable pairs often have natural constructions, and the theory usually gives us ways to infer about $f(X)$ if we know things about $F(X,\xp)$; however, getting information about $F$ is usually the difficult part in practice. As mentioned before, we shall deal with that in Chapter \ref{coupling}.

In the next two sections we are going to present two basic theoretical results. We shall work out our first serious example in section \ref{curieweiss}.

\section{A variance formula}
In this section, we establish a formula for the variance of $f(X)$. The importance of this simple formula lies in the fact that in a sense, it contains the essence of our theory; further embellishments are technical. 

We begin with the following fundamental lemma, which  is the foundation of all further investigation:
\begin{lmm}\label{fundamental}
For any measurable $h:\xx\ra\rr$
such that $\ee|h(X)F(X,\xp)|< \infty$, we have
\[
\ee(h(X)f(X)) = \frac{1}{2}\ee((h(X)-h(\xp))F(X,\xp)).
\]
\end{lmm}
{\bf Proof.} Note that $\ee(h(X)f(X)) = \ee(h(X)F(X,\xp))$. 
Using the exchangeability of $X$ and $\xp$, and the antisymmetric
nature of $F$, we have
\[
\ee(h(X)F(X,\xp)) = \ee(h(\xp)F(\xp,X)) = -\ee(h(\xp)F(X,\xp)).
\]
Thus,
\[
\ee(h(X)F(X,\xp)) = \frac{1}{2}\ee((h(X)-h(\xp))F(X,\xp)).
\]
This completes the proof. \hfill $\Box$\\

An immediate consequence of the above lemma is the following identity, which may be viewed as a weak law of large numbers for exchangeable  pairs:
\begin{thm}\label{weaklaw}
With the same notation as above, we have 
\[
\var(f(X)) =
\frac{1}{2}\ee((f(X)-f(\xp))F(X,\xp))
\]
whenever $\ee(f(X)^2)<\infty$.
\end{thm}
{\bf Proof.} Recall that $\ee(f(X)) = 0$ by assumption. Thus, we can directly use Lemma~\ref{fundamental} with $h = f$. \hfill $\Box$
\medskip

\noindent {\bf Remark.} At this point, we should mention that Reinert \cite{reinert95} has a technique for proving weak law type results for process valued functions using Stein's method, which she used to prove weak convergence of certain types of empirical processes. However, Reinert takes the usual Stein's method approach of treating the concentration problem as a problem of distributional approximation, where the target distribution is a point mass. Proceeding along this route cannot give exponential tail bounds.

As a second remark, note that from the above formula, we easily get $\ee(f(X)^2)\le \ee(v(X))$, where $v$ is as defined in (\ref{ddef}).  
\medskip

\noindent {\bf Example.} To quickly see how this works, let $X = \sum_{i=1}^n Y_i$, where $Y_i$'s are independent square integrable random variables. Let $\mu_i =\ee(Y_i)$ and $\sigma_i^2 = \var(Y_i)$. An exchangeable pair is created by choosing a coordinate $I$ uniformly at random from $\{1,\ldots,n\}$, and defining
\[
\xp = \sum_{j\ne I} Y_j + Y^\prime_I,
\]
where $Y^\prime_1,\ldots,Y^\prime_n$ are independent copies of $Y_1,\ldots,Y_n$. Let $F(x,y) = n(x-y)$. Then
\begin{align*}
\ee(F(X,\xp)|Y_1,\ldots,Y_n) &= \frac{1}{n}\sum_{i=1}^n \ee(n(Y_i - Y^\prime_i)|Y_1,\ldots,Y_n) = \sum_{i=1}^n (Y_i - \mu_i).
\end{align*}
Thus, we have $f(x) = x - \sum \mu_i$. The theorem now gives
\[
\var(X) = \frac{n}{2}\ee((X-\xp)^2) = \frac{1}{2}\sum_{i=1}^n \ee(Y_i - Y^\prime_i)^2 = \sum_{i=1}^n \sigma_i^2.
\]
Note that this example was meant to be just a basic illustration. It requires independence of the $Y_i$'s, while for the variance formula to hold, we just need them to be uncorrelated. Substantial examples will be provided in sections~\ref{curieweiss} and~\ref{meanmodels}.

\section{An exponential inequality}
Before working on further examples, let us describe our first concentration result, which is an exchangeable pairs version of the classical Hoeffding inequality (Theorem \ref{azuma}). Also, from now on we shall assume that
\begin{equation}\label{integrability}
\ee(e^{\theta f(X)}|F(X, \xp)|) < \infty \  \text{ for each } \theta \in \rr
\end{equation}
in all our theorems, with the exception of Theorem \ref{moments}, which is designed for situations where the above condition  does not hold. The validity of this assumption will be trivial to check in most of our examples.
\begin{thm}\label{hoeffding}
If $C$ is a constant such that $|v(X)|\le C$ almost surely (where $v(X)$ is defined in (\ref{ddef})), then we have $\ee(e^{\theta f(X)}) \le e^{C\theta^2/2}$ for all $\theta$, and consequently, 
$\pp\{|f(X)| \ge t\} \le 2e^{-t^2/2C}$ 
for each $t\ge 0$.
\end{thm}
{\bf Remark.} A version of this result about reversible Markov kernels was observed by Schmuckenschl\"ager \cite{schmuck98}. (Note that an exchangeable pair is formally equivalent to a reversible kernel.) However, Schmuckenschl\"ager restricts himself to a very special class of kernels (those with ``positive curvature'' in the sense of Bakry and \'Emery \cite{bakryemery85}) and the resulting technique becomes almost unusable in most practical problems. Indeed, the paper \cite{schmuck98}, though imaginative in many ways, has few concrete applications. 
\medskip

\noindent {\bf Example.}
To see that Theorem \ref{hoeffding} is a generalization of Hoeffding's inequality, consider $X= \sum_{i=1}^n Y_i$, where $Y_i$'s are now independent random variables with $\ee(Y_i)=\mu_i$ and $|Y_i-\mu_i|\le c_i$ almost surely, where $\mu_i$'s and $c_i$'s are finite constants. Creating the exchangeable pair as in the example following Theorem \ref{weaklaw}, and putting $f(X)=X-\ee(X)$, we get
\begin{align*}
v(X) &= \frac{1}{2n}\sum_{i=1}^n \ee(n(Y^\prime_i- Y_i)^2|X) \\
&= \frac{1}{2} \sum_{i=1}^n (\ee(Y_i-\mu_i)^2 + \ee((Y_i-\mu_i)^2|X)) \le \sum_{i=1}^n c_i^2, 
\end{align*}
and so our theorem gives the usual Hoeffding bound (Theorem \ref{azuma} in Chapter \ref{review}) in this case.\\
\\
{\bf Proof of Theorem \ref{hoeffding}.} Let $m(\theta) := \ee(e^{\theta f(X)})$ be the moment generating function of $f(X)$. We can differentiate $m(\theta)$ and move the derivative inside the expectation because of the integrability assumption (\ref{integrability}).
Thus, by Lemma \ref{fundamental}, we have
\begin{align*}
m^\prime(\theta) &= \ee(e^{\theta f(X)} f(X)) \\
&= \frac{1}{2} \ee((e^{\theta f(X)} - e^{\theta f(\xp)})F(X,\xp)).
\end{align*}
Now note that for any $x,y\in \rr$, 
\begin{align}\label{expin}
\biggl|\frac{e^x-e^y}{x-y}\biggr| &=\int_0^1 e^{tx + (1-t)y}dt \nonumber\\
&\le \int_0^1 (te^x+(1-t)e^y)dt \ \ \ \text{ (by the convexity of $u\mapsto e^u$)}\nonumber\\
&= \frac{1}{2}(e^x + e^y).
\end{align}
Using this inequality, and the exchangeablity of $X$ and $\xp$, we get
\begin{align*}
|m^\prime(\theta)| &\le \frac{|\theta|}{4}\ee((e^{\theta f(X)} + e^{\theta f(\xp)})|(f(X)-f(\xp)) F(X,\xp)|) \\
&\le \frac{|\theta|}{2} \ee(e^{\theta f(X)} v(X) + e^{\theta f(\xp)} v(\xp))\\
&= |\theta| \ee(e^{\theta f(X)} v(X)) \le C|\theta| m(\theta).
\end{align*}
For $\theta \ge 0$, this means $\frac{d}{d\theta}\log m(\theta) \le C\theta$. Since $m(0)=1$, we have $\log m(\theta)\le C\theta^2/2$.
Given $t\ge 0$, we can choose $\theta = t/C$ and get
\[
\pp\{f(X)\ge t\} \le e^{-\theta t + \log m(\theta)} \le e^{-t^2/2C}.
\] 
Similarly, $\pp\{f(X)\le -t\} \le e^{-t^2/2C}$. \hfill $\Box$

\section{Example: Curie-Weiss model}\label{curieweiss}
In this section, we shall work out concentration bounds for the magnetization in the Curie-Weiss model of ferromagnetic interaction. This is arguably the simplest statistical mechanical model of spin systems. For a detailed mathematical treatment of this model, we refer to the book \cite{ellis85} by Richard Ellis. In the next section, we shall extend our technique for the Curie-Weiss model to a more general class of models with quadratic interaction. 

The state space is $\{-1,1\}^n$, the space of all possible spin configurations of $n$ particles. A typical configuration will be denoted by $\sigma = (\sigma_1,\ldots,\sigma_n)$. 
The Curie-Weiss model at inverse temperature $\beta$ and a fixed external magnetic field $h$ prescribes a joint probability density (the Gibbs measure) for these spins by the formula $p_{\beta, h}(\sigma) = Z_{\beta, h}^{-1} e^{-\beta H_h(\sigma)}$ where $Z_{\beta, h}$ is the normalizing constant, and 
\[
H_h(\sigma) := -\frac{1}{n}\sum_{i<j} \sigma_i\sigma_j - h \sum_{i=1}^n \sigma_i
\]
is the Hamiltonian for the system. We shall henceforth assume that $\beta$ and $h$ are fixed and omit the subscripts. One quantity of interest is the magnetization, defined as
\[
m = m(\sigma) := \frac{1}{n}\sum_{i=1}^n \sigma_i.
\]
It is common knowledge in the physics circles (and proved rigorously in \cite{ellis85}) that if $n$ is large and $\sigma$ is drawn from the Gibbs measure, the random variable $m(\sigma)$ is concentrated in a neighborhood of the set of solutions of the equation 
\[
x = \tanh(\beta x + \beta h).
\]
The equation has a unique root for small values of $\beta$ (the ``high temperature phase'') and multiple solutions for $\beta$ above a critical range (the ``low temperature phase''). For example, when $h=0$, $\beta_c =1$ is the critical value. 

In this context, we can use Theorems \ref{weaklaw} and \ref{hoeffding} quite easily to prove the following result:
\begin{prop}\label{curie}
The magnetization $m$ in the Curie-Weiss model satisfies 
\[
\ee\bigl(m - \tanh(\beta m + \beta h)\bigr)^2 \le \frac{2+2\beta}{n} + \frac{\beta^2}{n^2}.
\]
Moreover, we also have the concentration bound
\[
\pp\bigl\{\bigl|m -  \tanh(\beta m + \beta h)\bigr| \ge \frac{\beta}{n} + t\bigr\} \le 2e^{-nt^2/(4+4\beta)}.
\]
Here $\ee$ and $\pp$ denote expectation and probability under the Gibbs measure at inverse temperature $\beta$ and external field $h$.
\end{prop}
{\bf Remark.} Although the Curie-Weiss model is considered to be the simplest model of ferromagnetic interaction, we haven't encountered any result in the literature which gives bounds like the above. In fact, it is unlikely that there exists a short and simple bare-hands argument which gives such bounds.

The classical way to solve the Curie-Weiss model uses ideas from large deviations theory, and involves a variational problem in the limit (see Ellis \cite{ellis85}, Chapter IV). The magnetization is obtained by the usual technique of differentiating the free energy function. Although the limiting result is straightforward, getting finite sample tail bounds will not be easy along this route. 
\medskip

\noindent {\bf Proof.} Suppose $\sigma$ is drawn from the Gibbs distribution. We construct $\sss$ by taking a step in the Gibbs sampler as follows: Choose a coordinate $I$ uniformly at random, and replace the $I^{\mathrm{th}}$ coordinate of $\sigma$ by an element drawn from the conditional distribution of the $I^{\mathrm{th}}$ coordinate given the rest. 
It is well-known and easy to prove that $(\sigma,\sss)$ is an exchangeable pair. Let 
\[
F(\sigma, \sss) := \sum_{i=1}^n (\sigma_i - \sss_i).
\]
Now define
\[
m_i(\sigma) := \frac{1}{n}\sum_{j \le n, j\ne i} \sigma_j, \ \ i=1,\ldots,n.
\]
Since the Hamiltonian is a simple explicit function, the conditional distribution of the $i^{\mathrm{th}}$ coordinate given the rest is easy to obtain. An easy computation gives
$\ee(\sigma_i|\{\sigma_j, j\ne i\}) = \tanh(\beta m_i+\beta h)$. 
Thus, we have
\begin{align*}
f(\sigma) = \ee(F(\sigma,\sss)|\sigma) &= \frac{1}{n}\sum_{i=1}^n (\sigma_i -\ee(\sigma_i|\{\sigma_j, j\ne i\}))\\
&= m - \frac{1}{n}\sum_{i=1}^n \tanh(\beta m_i + \beta h).
\end{align*}
Now note that $|F(\sigma,\sss)| \le 2$, because $\sigma$ and $\sss$ differ at only one coordinate. Also, since the map $x\mapsto \tanh x$ is $1$-Lipschitz, we have 
\begin{align*}
|f(\sigma)-f(\sss)| &\le |m(\sigma)-m(\sss)| +\frac{\beta}{n}\sum_{i=1}^n |m_i(\sigma) - m_i(\sss)| \le \frac{2(1+\beta)}{n}.
\end{align*}
Thus, by Theorem \ref{weaklaw} we have
\[
\var\biggl(m - \frac{1}{n}\sum_{i=1}^n \tanh(\beta m_i + \beta h)\biggr) \le \frac{2(1+\beta)}{n}
\]
and from Theorem \ref{hoeffding},
\[
\pp\bigl\{\bigl|m - \frac{1}{n}\sum_{i=1}^n \tanh(\beta m_i + \beta h)\bigr| \ge t\bigr\} \le 2e^{-nt^2/(4+4\beta)}.
\]
Finally note that for each $i$, by the Lipschitz nature of the $\tanh$ function, we get
\[
|\tanh(\beta m_i+\beta h) - \tanh(\beta m+\beta h)| \le \beta|m_i-m| \le \frac{\beta}{n}.
\]
This completes the proof of the theorem. \hfill $\Box$

\section{Validity of naive mean field equations}\label{meanmodels}
The proof of Proposition \ref{curie} shows that there is no reason why the argument will not generalize to the Hamiltonian
\begin{equation}\label{quadmodel}
H(\sigma) := -\sum_{i<j}J_{ij} \sigma_i\sigma_j - \sum_{i=1}^n h_i\sigma_i
\end{equation} 
where $J = (J_{ij})_{1\le i,j\le n}$ is a symmetric interaction matrix with zeros on the diagonal, and $h_1,\ldots,h_n$ are fixed real numbers. In this section, we are going to carry through an extension of our treatment of the Curie-Weiss model to find conditions under which a certain widely used (but mathematically dubitable) technique from physics --- the naive mean field equations --- is valid.

Throughout this section, we shall use the standard physical notation $\langle\cdot\rangle$ to denote the expectation with respect to the Gibbs measure. 

Now, it is an easy fact that $\ee_\beta(\sigma_i|\{\sigma_j, j\ne i\}) = \tanh(\beta \sum_{j=1}^n J_{ij}\sigma_j + h_i)$. This gives what are called the {\it Callen equations} (see e.g.\ Chapter 3 of \cite{parisi88}):
\begin{equation}\label{callen}
\langle \sigma_i\rangle = \biggl \langle \tanh\biggl( \beta \sum_{j=1}^n J_{ij}\sigma_j + \beta h_i\biggr) \biggr\rangle, \ \ \ i=1,\ldots,n.
\end{equation}
These equations are somehow just a restatement of the model, and not useful for anything deeper. The {\it naive mean field equations} are a modification of the Callen equations. The general physical intuition is that, when the spins are sufficiently ``uncorrelated'', the expectation on the right hand side of the Callen equations can be taken inside the $\tanh$. This gives
\begin{equation}\label{meanfield}
``\langle \sigma_i\rangle = \tanh\biggl( \beta \sum_{j=1}^n J_{ij}\langle\sigma_j\rangle + \beta h_i\biggr) , \ \ \ i=1,\ldots,n\text{''}.
\end{equation}
The reason why this is written within quotes is that no reasonable general condition is known for the validity of these equations. Physicists know that these are {\it not} valid in models with short range interactions like the Ising model, but generally work when the interactions are of very long range (like the Curie-Weiss, for instance). Rigorous texts like \cite{ellis85} prefer to avoid talking about mean field equations. Still, this technique is used with reasonable measures of success in a variety of fields, including computer science, image processing, neural netwroks and other computational sciences. A system of $n$ equations may seem strange at first sight, but the fact is that they often reduce to a much smaller system in practice (for instance, in the Curie-Weiss model, they reduce to just one equation). For a recent survey of applications in areas outside of physics, one can look at \cite{oppersaad01}. For the  physical perspective, a classical text is Parisi's book~\cite{parisi88}.

We shall now show that, although it is difficult to obtain general conditions under which the average spins $\{\langle\sigma_i\rangle, \ i=1,\ldots,n\}$ satisfy the mean field equations, the {\it conditional averages} of the spins, defined as
\begin{equation}\label{condx}
\langle\sigma_i\rangle^- := \ee(\sigma_i|\{\sigma_j, j\ne i\}) = \tanh(\beta\sum_j J_{ij} \sigma_j + \beta h_i), \ \ i =1,\ldots,n
\end{equation}
satisfy (\ref{meanfield}) with high probability whenever the Hilbert-Schmidt norm of $J$ (denoted by $\|J\|_{HS}$) is not too large. The precise condition is
\begin{equation}\label{long}
\|J\|_{HS} := \bigl[\sum_{i,j} J_{ij}^2\bigr]^{1/2} \ll n^{1/3}.
\end{equation}
We shall also show that the mean field equations (\ref{meanfield}) themselves are approximately valid at sufficiently high temperature in the class of models satisfying~(\ref{long}). The results are summarized in the following theorem:
\begin{thm}\label{meanthm} 
Fix $\beta \ge 0$.
Let $\rho = \|J\|_{HS} = [\sum_{i,j}J_{ij}^2]^{1/2}$. Let $\langle \sigma_i\rangle^-$, $i=1,\ldots,n$ be defined as in (\ref{condx}). For $1\le i\le n$, define
\begin{align*}
\tilde{\varepsilon}_i &:= \langle\sigma_i\rangle^- - \tanh\biggl(\beta\sum_{j=1}^n J_{ij}\langle \sigma_j\rangle^- + \beta h_i\biggr), \ \ \text{ and} \\
\varepsilon_i &:= \langle\sigma_i\rangle -  \tanh\biggl(\beta\sum_{j=1}^n J_{ij}\langle \sigma_j\rangle + \beta h_i\biggr).
\end{align*}
Then we have the bound $\langle\tilde{\varepsilon}_i^2\rangle \le 2\beta^2(1+\beta\rho)\sum_{j=1}^n J_{ij}^2$ for each $i$, and hence
\begin{equation}\label{tep}
\frac{1}{n}\sum_{i=1}^n \langle\tilde{\varepsilon}_i^2 \rangle\le \frac{1}{n}2(1 + \beta \rho)\beta^2\rho^2.
\end{equation}
We also have the tail bounds 
\begin{equation}\label{tail}
\langle |\tilde{\varepsilon}_i| \ge t\rangle \le 2\exp\biggl(-\frac{t^2}{4\beta^2(1+\beta \rho) \sum_{j=1}^n J_{ij}^2}\biggr)
\end{equation}
for $1\le i \le n$. Moreover, if $\beta < 1/\rho$, we have
\begin{equation}\label{ep}
\frac{1}{n}\sum_{i=1}^n \varepsilon_i^2 \le \frac{2(1+\beta\rho)\beta^2\rho^2}{n(1-\beta\rho)^2}.
\end{equation}
\end{thm}
{\bf Discussion and examples.} 
Before we come to a critical discussion, let us, for a moment, gloat over the fact that it is difficult to imagine that such bounds (even the second moment bounds) can be easily obtained by bare-hands arguments, or by any available technique, for that matter!

The first part of the theorem says that under (\ref{long}), the conditional averages $\{\langle \sigma_i\rangle^-, \ i=1,\ldots,n\}$, instead of the unconditional ones, satisfy the mean field equations (\ref{meanfield}). In other words, we have a set of modified mean field equations
\begin{equation}\label{modmean}
``\langle \sigma_i\rangle^- = \tanh\biggl( \beta \sum_{j=1}^n J_{ij}\langle\sigma_j\rangle^- + \beta h_i\biggr) , \ \ \ i=1,\ldots,n\text{''},
\end{equation}
which are approximately valid with high probability. 
In particular, if each $\langle \sigma_j \rangle^-$ is concentrated around its mean, then we can recover the mean field equations from the information in (\ref{tep}). One condition under which this happens is the ``high temperature conditon'' $\beta <1/\rho$, which gives the second part of the theorem. But presumably the concentration of $\langle \sigma_j\rangle^-$ can happen under other conditions also.

The condition (\ref{long}) can be interpreted as some sort of a condition of long range interaction, as will be made clear by the following example.
Consider a graph $G = (V,E)$ on $V=\{1,\ldots,n\}$ with constant degree $r$. The (normalized) Ising Hamiltonian on $\{-1,1\}^V$ is defined as
\[
H(\sigma) = - \frac{1}{r}\sum_{(i,j)\in E} \sigma_i \sigma_j - h\sum_{i=1}^n \sigma_i,
\]
where $h$ is the external field. Thus, $J_{ij} = 1/r$ if $(i,j)\in E$, and $=0$ otherwise. Since the graph has constant degree $r$, a simple computation gives $\|J\|_{HS} = \sqrt{n/r}$. Thus, the modified mean field equations (\ref{modmean}) hold whenever $r \gg n^{1/3}$, which is a long range condition.

For another example, consider the following generalization of the Curie-Weiss model: Suppose each $J_{ij}$ is bounded in magnitude by $1/n$. Then $\|J\|_{HS}\le 1$, and therefore the modified mean field equations hold. Moreover, the second part of the theorem shows that in this class of models, if $\beta < 1$, then the mean field equations~(\ref{meanfield}) are valid.
\vskip.2in
\noindent {\bf Proof of Theorem \ref{meanthm}.} Throughout this proof, we shall use the fact that $\tanh$ is a Lipschitz function without explicit mention. Also, as a natural extension of the notation in the previous section, we shall write
\[
m_i = m_i(\sigma) := \sum_{j=1}^n J_{ij}\sigma_j + h_i.
\]
Note that $\langle \sigma_i\rangle^- = \tanh(\beta m_i)$. 

Now fix $i$, $1\le i\le n$. Construct $(\sigma, \sss)$ as in the proof of Proposition \ref{curie}, by choosing $\sigma$ from the Gibbs measure and taking a step in the Gibbs sampler to get $\sss$ as follows: Choose a coordinate $I$ uniformly at random, and replace the $I^{\mathrm{th}}$ coordinate of $\sigma$ by a random sample from the conditional distribution of $\sigma_I$ given $\{\sigma_j, j\ne I\}$. Define
$F(\sigma,\sss) = n\sum_{j=1}^n J_{ij}(\sigma_j-\sss_j) = nJ_{iI}(\sigma_I-\sss_I)$. 
Then, as before, 
\begin{equation}\label{thisf}
f(\sigma) = \sum_{j=1}^n J_{ij}(\sigma_j - \tanh(\beta m_j)) = m_i - \sum_{j=1}^n J_{ij} \tanh(\beta m_j) - h_i.
\end{equation}
Now, if $\sigma$ and $\sss$ differ at site $I$, then
\begin{align*}
|f(\sigma)-f(\sss)| &\le 2|J_{iI}| + \sum_{j=1}^n |J_{ij}(\tanh (\beta m_j(\sigma))- \tanh(\beta m_j(\sss)))|\\
&\le 2|J_{iI}| + \sum_{j=1}^n |J_{ij}| |\beta (m_j(\sigma)- m_j(\sss))|\\
&= 2|J_{iI}| + 2\beta \sum_{j=1}^n |J_{ij}J_{jI}|.
\end{align*}
Also, $|F(\sigma, \sigma^I)|\le 2n|J_{iI}|$. Thus, we have
\begin{align}\label{vbd}
v(\sigma) &= \frac{1}{2}\ee\bigl(|(f(\sigma)-f(\sss))F(\sigma, \sss)|\bigl|\sigma\bigr) \nonumber \\
&\le \frac{1}{2}\sum_{k=1}^n (2|J_{ik}| + 
2\beta\sum_{j=1}^n |J_{ij}J_{jk}|) 2|J_{ik}|\nonumber \\
&\le 2\sum_{k=1}^n J_{ik}^2 + 2\beta \bigl[ \sum_{j,k}J_{ij}^2J_{ik}^2\bigr]^{1/2} \bigl[\sum_{j,k} J_{jk}^2\bigr]^{1/2} \nonumber \\
&= 2(1+\beta\rho)\sum_{j=1}^n J_{ij}^2.
\end{align}
Thus, by Theorem \ref{weaklaw} and (\ref{thisf}), we get
\begin{align}\label{varf}
\bigl\langle \bigl(m_i - \sum_{j=1}^n J_{ij} \tanh(\beta m_j) - h_i\bigr)^2\bigr\rangle  &= \langle f(\sigma)^2\rangle \le 2(1+\beta \rho )\sum_{j=1}^n J_{ij}^2.
\end{align}
To complete the proof of the first set of inequalities in Theorem \ref{meanthm}, observe that
\begin{align*}
\langle\tilde{\varepsilon}_i^2\rangle &= \biggl\langle\biggl(\tanh(\beta m_i) -  \tanh\biggl(\beta\sum_{j=1}^n J_{ij}\tanh(\beta m_j) + \beta h_i\biggr)\biggr)^2\biggr\rangle\\
&\le \beta^2 \bigl\langle \bigl(m_i - \sum_{j=1}^n J_{ij} \tanh(\beta m_j) - h_i\bigr)^2\bigr\rangle,
\end{align*}
and combine with (\ref{varf}). The tail bound on $\tilde{\varepsilon}_i$ follows from the same analysis, using Theorem \ref{hoeffding}. 

Next, for each $i$ define the function $g_i:\rr^n \ra \rr$ as $g_i(x) = \sum_{j=1}^n J_{ij}\tanh(\beta x_j) + h_i$. Define $g:\rr^n \ra \rr^n$ as $g(x) := (g_1(x),\ldots,g_n(x))$.
Then note that for any $x,y\in \rr^n$,
\begin{align*}
\|g(x)-g(y)\|^2 &= \sum_{i=1}^n \biggl(\sum_{j=1}^n J_{ij}(\tanh(\beta x_j)-\tanh(\beta y_j))\biggr)^2 \\
&\le \beta^2\sum_{i=1}^n \biggl(\sum_{j=1}^n |J_{ij}(x_j-y_j)|\biggr)^2 \\
&\le \beta^2 \rho^2 \|x-y\|^2.
\end{align*}
Thus, if $\beta <1/\rho$, then $g$ is a contraction map with respect to the Euclidean metric. Therefore by  the well-known theorem about contraction maps (see, e.g.\ Theorem 9.23 in Rudin \cite{rudin76}), $g$ has a unique fixed point, which we shall call $x^*$. Now note that for any $x\in \rr^n$,
\begin{align*}
\|x-x^*\| &\le \|(x-g(x))-(x^*-g(x^*))\| + \|g(x)-g(x^*)\| \\
&\le \|x-g(x)\| + \beta \rho\|x-x^*\|.
\end{align*}
Thus for any $x$, we have $\|x-x^*\|\le (1-\beta\rho)^{-1}\|x-g(x)\|$. Applying to $x=(m_1,\ldots m_n)$, and remembering that $\langle (m_i - a)^2\rangle$ is minimized at $a = \langle m_i\rangle$, we get
\begin{align*}
\sum_{i=1}^n \langle(m_i-\langle m_i\rangle)^2\rangle &\le \sum_{i=1}^n\langle (m_i-x_i^*)^2\rangle \\
&\le \frac{1}{(1-\beta\rho)^2} \sum_{i=1}^n \bigl\langle \bigl(m_i - \sum_{j=1}^n J_{ij} \tanh (\beta m_j) - h_i\bigr)^2\bigr\rangle \\
&\le \frac{2(1+\beta \rho )\rho^2}{(1-\beta\rho)^2} \ \ \ \text{ by inequality (\ref{varf}).}
\end{align*}
Finally, observe that
\begin{align*}
\sum_{i=1}^n\varepsilon_i^2 &= \sum_{i=1}^n \bigl(\langle \tanh (\beta m_i) \rangle - \tanh(\beta \langle m_i\rangle)\bigr)^2\\
&\le \beta^2\sum_{i=1}^n \bigl\langle (m_i - \langle m_i\rangle )^2\bigr\rangle.
\end{align*}
Combined with the previous step, this completes the proof. \hfill $\Box$

\section{Example: Least squares for the Ising model}\label{least}
Consider an undirected graph $G = (V,E)$, where $V = \{1,\ldots,n\}$. Let $r$ be the maximum degree of $G$. Let $X = (X_1,\ldots,X_n)$ be a random element of $\{-1,1\}^G$. The {\it Ising model} assigns a probability distribution for $X$ according to the formula
\begin{equation}\label{ising}
\pp_\theta\{X= x\} = Z(\theta)^{-1} e^{-\theta \sum_{(i,j)\in E } x_ix_j}
\end{equation}
where $\theta$ is an unknown parameter (the ``inverse temperature'', usually denoted by $\beta$ in the physics literature), and $Z(\theta)$ is the normalizing constant. The parameter space is $\Omega := [0,\infty)$. 
The natural statistical problem in this model is the following: How to make inference about $\theta$ from a single realization of $X$? 
This is one of most elementary (yet analytically almost intractable) models of dependent discrete data.

The classical maximum likelihood approach for this problem has been discussed in detail by Pickard \cite{pickard87}. The main difficulty with directly computing the maximum likelihood estimator (MLE) is that the normalizing constant has no explicit form except in very special cases, and there is no polynomial time algorithm for exact numerical evaluation.

In a well-cited paper, Geyer \& Thompson \cite{geyerthompson92} devised a feasible Monte Carlo technique for computing the MLE in models like the above. One of the examples considered in that paper involves an instance of the autologistic model of Besag \cite{besag72}, which is a generalization of the Ising model. The autologistic model of binary data assumes that the conditional distribution of each $X_i$ given the rest can be modeled as a logistic regression, that is,
\[
\log\frac{\pp\{X_i=1|(X_j,j\ne i)\}}{1-\pp\{X_i=1|(X_j,j\ne i)\}} = \alpha_i + \sum_{j\le n, j\ne i} \beta_{ij} X_j,
\]
where the $\alpha_i$'s and the $\beta_{ij}$'s are known functions of a collection of unknown parameters. A simple verification shows that the Ising model described by (\ref{ising}) is a special case of the autologistic model.

Monte Carlo algorithms for computing the MLE in such  problems are widely used nowadays, but lingering doubts remain about the rates of convergence, specially in models with high degrees of dependence like the ones described above. 

A method which does not require simulations, but fell out of favor due to efficiency issues after the advent of Monte Carlo techniques, is the maximum pseudolikelihood approach introduced by Besag \cite{besag75}. Besag's pseudolikelihood function is defined as
\[
f(\theta|X) :=\prod_{i=1}^n f_i(\theta|X),
\]
where $f_i(\theta|X)$ is the conditional density of $X_i$ given $(X_j,j\ne i)$ under $\theta$. Now, if $X$ is a gaussian vector such that $\var_\theta (X_i|(X_j,j\ne i))$  is a constant independent of $i$ and $\theta$, the pseudolikelihood problem reduces to minimizing
\begin{equation}\label{stheta}
S(\theta) := \frac{1}{n}\sum_{i=1}^n \bigl(X_i - \ee_\theta(X_i|(X_j,j\ne i))\bigr)^2,
\end{equation}
which may be called a ``conditional least squares problem''. We shall now show that somewhat surprisingly, minimizing $S(\theta)$ to estimate $\theta$ may be a good idea even outside the gaussian framework. In particular, we shall show by way of example later in this section that it works in the Ising model, and the method can even be extended to construct nonasymptotic confidence intervals for $\theta$. One of the interesting consequences is the following:
\begin{prop}\label{isingprop}
In the Ising model on a graph with maximum degree $r$ as defined in (\ref{ising}), we have the bound
\[
\ee_\theta(S(\theta) - \inf_\psi S(\psi)) \le C\sqrt{\frac{r\log n}{n}} \ \ \text{ for every $\theta \ge 0$},
\]
where $C$ is a numerical constant. 
\end{prop}
This shows that the conditional least squares inference for $\theta$ should work in this model when $r\ll n/\log n$. However, there is a caveat: Although our results will show that $S(\theta)$ is minimized near the true value of $\theta$ under pretty general conditions, they say nothing about the sharpness of the minimum of $S(\theta)$. If $S$ is too flat, our results will be of no use.

In the remainder of this section, we shall adopt the notation $\bar{x}^i$ for the vector $(x_1,\ldots,x_{i-1},  x_{i+1},\ldots,x_n)$, where $x$ is the original vector $(x_1,\ldots,x_n)$. 
Let $X=(X_1,\ldots,X_n)$ be a random vector in $[-1,1]^n$, whose distribution is parametrized by a parameter $\theta$ belonging to some parameter space $\Omega$. For $1\le i\le n$, let
\[
\mu_i(\theta, \bar{x}^i) := \ee_\theta(X_i|\bar{X}^i = \bar{x}^i),
\]
For each $\theta$, let  $\{a_{ij}(\theta)\}_{1\le i,j\le n}$ be an array of nonnegative real numbers with zeros on the diagonal such that for any $x, y,i$ and $\theta$ we have
\[
|\mu_i(\theta, \bar{x}^i) - \mu_i(\theta,\bar{y}^i)| \le \sum_{j=1}^n a_{ij}(\theta)\ii\{x_j\ne y_j\},
\]
and let
\begin{equation}\label{mdef}
M := \frac{1}{n}\sup_{\theta \in \Omega} \sum_{i,j} a_{ij}(\theta).
\end{equation}
Lemma \ref{lslmm} below will show that when $M \ll n$, which holds in many models (including Ising type models), 
$S(\psi) - S(\theta)$ is ``approximately nonnegative with high probability'' whenever $\theta$ is the true value of the parameter and $\psi$ is any other value. This will be shown by decomposing $S(\psi)-S(\theta)$ as $A(\psi,\theta)+B(\psi,\theta)$, where $A \ge 0$ and $B\approx 0$ under $\pp_\theta$. The explicit expressions for $A$ and $B$ are as follows:
\[
A(\psi,\theta) = 
\frac{1}{n}\sum_{i=1}^n (\mu_i(\theta,\bar{X}^i) - \mu_i(\psi, \bar{X}^i))^2
\]
and 
\[
B(\psi,\theta)=  \frac{2}{n} \sum_{i=1}^n (\mu_i(\theta, \bar{X}^i) - \mu_i(\psi, \bar{X}^i))
(X_i - \mu_i(\theta, \bar{X}^i)).
\]
It is easy to check that indeed $S(\psi)-S(\theta) = A(\psi,\theta) + B(\psi,\theta)$. It is also apparent that $A(\psi,\theta)\ge 0$. The following lemma, the proof of which depends heavily on our techniques, shows that $B(\psi,\theta)\approx 0$ with high probability under $\theta$ whenever $M\ll n$. It is not completely honest, though, because the constants are bad. 
\begin{lmm}\label{lslmm}
Suppose $X$ takes values in $[-1,1]^n$ and let $M$, $A$ and $B$ be defined as above. Then for any $\psi, \theta\in \Omega$ and $t\ge 0$, we have 
\[
\pp_{\theta}\{|B(\psi,\theta)|\ge t\} \le 2e^{-nt^2/(96M + 32)}.
\] 
\end{lmm}
{\bf Proof.} Fix $\psi, \theta\in\Omega$. We produce $\xp$, as usual, by taking a step in the Gibbs sampler: A coordinate $I$ is chosen uniformly at random, and $X_I$ is replace by $\xp_I$ drawn from the conditional distribution (under $\pp_\theta$) of the $I^{\mathrm{th}}$ coordinate given $\bar{X}^I$. 
Next, for each $i$, let
\[
g_i(\bar{x}^i) := \mu_i(\theta, \bar{x}^i) - \mu_i(\psi, \bar{x}^i),
\]
and define
\[
F(x,y) = \sum_{i=1}^n (g_i(\bar{x}^i) + g_i(\bar{y}^i))(x_i - y_i).
\]
Clearly, $F$ is antisymmetric. Note that
\begin{align*}
f(X) := \ee_{\theta}(F(X,\xp)|X) &= 2\ee_{\theta}(g_I(\bar{X}^I)(X_I - \xp_I)|X) \\
&= \frac{2}{n}\sum_{i=1}^n g_i( \bar{X}^i)(X_i - \mu_i(\theta, \bar{X}^i)). 
\end{align*}
Thus, in our notation, $f(X) = B(\psi,\theta)$. From the given conditions it easily follows that if $x$ and $y$ are elements of $[-1,1]^n$ which differ only in the $i^{\mathrm{th}}$ coordinate, then
\[
|f(x)-f(y)| \le \frac{4}{n}\biggl(\sum_{j=1}^n (a_{ji}(\psi) + 2a_{ji}(\theta)) + 1\biggr).
\]
Also, quite clearly, $|F(X,\xp)|\le 8$. Combining, we have
\begin{align*}
v(X) &= \frac{1}{2}\ee_{\theta}\bigl(|(f(X)-f(\xp))F(X,\xp)|\bigl| X\bigr) \\
&\le \frac{16}{n^2}\sum_{i=1}^n \sum_{j=1}^n (a_{ji}(\psi) +2a_{ji}(\theta)) + \frac{16}{n}.
\end{align*}
Invoking Theorem \ref{hoeffding} completes the proof. \hfill $\Box$\\
\\

The above lemma is only a moral justification for using the ``conditional least squares approach'' for estimating parameters in models of dependent data. It does not show that the true value of $\theta$ is an approximate {\it global} minimizer of $S(\theta)$. For that, and also for constructing confidence regions, we need tail bounds on $S(\theta)-\inf_\psi S(\psi)$. For instance, a $(1-\alpha)$-level {\it confidence region} for the true value of $\theta$ can be defined as $\{\theta: S(\theta)-\inf_{\psi\in \Omega} S(\psi) \le t_\alpha\}$, where $t_\alpha$ is chosen such that for any $\theta$, 
\[
\pp_{\theta}\{S(\theta) - \inf_\psi S(\psi) > t_\alpha\} \le \alpha.
\]
For doing all that, we need to introduce some more notation. Define a pseudometric $d$ on $\Omega$ as follows:
\[
d(\theta, \theta^\prime) := \sup_x \frac{1}{n}\sum_{i=1}^n |\mu_i(\theta, \bar{x}^i) - \mu_i(\theta^\prime, \bar{x}^i)|.
\]
For every $\varepsilon > 0$, let $N_d(\varepsilon)$ denote the minimum number of closed $\varepsilon$-balls (w.r.t.\ the metric $d$) required to cover $\Omega$. We have the following result:
\begin{thm}\label{lse}
Suppose $X$ takes values in $[-1,1]^n$. Let $S(\theta)$ and $M$ be defined as in (\ref{stheta}) and (\ref{mdef}), and let $N_d$ be the covering number defined above. Then, for any $\varepsilon >0$ and any $\theta \in \Omega$, we have
\[
\pp_{\theta}\{S(\theta)- \inf_{\psi\in \Omega} S(\psi) \ge 4\varepsilon + t\} \le 2N_d(\varepsilon) e^{-nt^2/(96M + 32)}
\]
for all $t\ge 0$. Consequently, we also have
\[
\ee_{\theta}(S(\theta) - \inf_{\psi\in \Omega} S(\psi))\le \inf_{\varepsilon > 0}\biggl\{4\varepsilon + 2\sqrt{\frac{(96 M+32)\log (2N_d(\varepsilon))}{n}}\biggr\}.
\]
\end{thm}
\noindent {\bf Application to the Ising model.} For the Ising model described at the beginning of the section, it is easy to verify that
$\mu_i(\theta, \bar{x}^i) = \tanh \bigl(\theta\sum_{j\in N(i)} x_j\bigr)$, 
where $N(i)$ is the neighborhood of $i$ in $G$. 
Now, since $\tanh x  \in [-1,1]$ for all $x\in \rr$, we have
\[
|\mu_i(\theta, \bar{x}^i) - \mu_i(\theta,\bar{y}^i)| \le 2\sum_{j\in N(i)} \ii\{x_j\ne y_j\}.
\]
Thus, we can take $a_{ij}(\theta) = 2\ii\{(i,j)\in E\}$ irrespective of $\theta$. With this choice of $a_{ij}$'s, we have $M = 2r$, where recall that $r$ is the maximum degree of $G$. 

Now note that $\sum_{j\in N(i)} x_j \in \{0,\pm 1, \ldots,\pm r\}$, and $\tanh$ is an odd function. Therefore for any $\theta, \theta^\prime$ and $x$,
\begin{align*}
|\mu_i(\theta, \bar{x}^i) - \mu_i(\theta^\prime, \bar{x}^i)| &\le \sup_{s\in \{0,1,\ldots,r\}} |\tanh(\theta s) - \tanh(\theta^\prime s)|. 
\end{align*}
Since $\tanh$ is a Lipschitz function, therefore the right hand side is bounded by $|\theta-\theta^\prime| r$. On the other hand, since $\tanh$ is an increasing function bounded by $1$, therefore the right hand side is also bounded by $1-\tanh (\min\{\theta,\theta^\prime\})\le \exp(-\min\{\theta,\theta^\prime\})$. Combining, we have
\begin{align*}
d(\theta,\theta^\prime) &= \sup_x\frac{1}{n}\sum_{i=1}^n |\mu_i(\theta, \bar{x}^i) - \mu_i(\theta^\prime, \bar{x}^i)|\le \min\{|\theta-\theta^\prime|r, e^{-\min\{\theta,\theta^\prime\}}\}.
\end{align*}
Now fix any $\varepsilon\in (0,1)$. Let $L = -\log \varepsilon$. Equipartition the interval $[0,L]$ into $[Lr/2\varepsilon]+1$ subintervals of length $\le 2\varepsilon/r$ each. The above bound shows that the end points of these subintervals, together with $L$, form an $\varepsilon$-net with respect to the pseudometric $d$ on $[0,\infty)$. Thus, 
\[
N_d(\varepsilon) \le \frac{r |\log \varepsilon|}{2\varepsilon} + 2.
\]
It is now easy to apply Theorem \ref{lse}. The second bound directs us to the optimal choice of $\varepsilon$, which we take to be $\varepsilon = \sqrt{r/n}$. Theorem \ref{lse} now gives
\[
\pp_{\theta}\{S(\theta)-\inf_\psi S(\psi) \ge 4\sqrt{\frac{r}{n}} + t\} \le 2\sqrt{rn}(\log n)e^{-nt^2/(96r+32)}, 
\]
and 
\[
\ee_{\theta} (S(\theta) - \inf_\psi S(\psi)) \le C\sqrt{\frac{r\log n}{n}},
\]
where $C$ is a computable numerical constant. This proves Proposition \ref{isingprop}.\\
\\
{\bf Proof of Theorem \ref{lse}.} Fix $\varepsilon >0$. Let $k = N_d(\varepsilon)$, and let $\theta_1,\ldots,\theta_k$ be the centers of a collection of $\varepsilon$-balls which cover $\Omega$. Then for every $\theta$ there exists $i$ such that $|S(\theta)-S(\theta_i)| \le 4d(\theta, \theta_i) \le 4\varepsilon$. Thus,
\[
\min_{1\le i \le k} S(\theta_i) - \inf_\theta S(\theta)  \le 4\varepsilon.
\]
Therefore,
\begin{align*}
&\pp_{\theta}\{S(\theta) - \inf_\psi S(\psi) \ge 4\varepsilon + t\} \\
&\le \pp_{\theta}\{|S(\theta)-\inf_{1\le i\le k} S(\theta_i)| \ge t\} 
\le \sum_{i=1}^k \pp_{\theta}\{|S(\theta)-S(\theta_i)| \ge t\}.
\end{align*}
The first bound follows from this and Lemma \ref{lslmm}. Now note that for any $b >0$ and $a> \sqrt{e}$, we have 
\begin{align*}
\int_0^\infty (ae^{-bt^2} \wedge 1) dt &= \sqrt{b^{-1} \log a} + \int_{\sqrt{b^{-1}\log a}}^\infty ae^{-bt^2} dt \\
&\le \sqrt{b^{-1}\log a} + \int_{\sqrt{b^{-1}\log a}}^\infty \frac{t}{\sqrt{b^{-1}\log a}} ae^{-bt^2} dt\\
&= \sqrt{b^{-1}\log a} + \frac{1}{2b\sqrt{b^{-1}\log a}} \\
&\le 2\sqrt{b^{-1}\log a} \ \ \text{(since $a \ge \sqrt{e} \implies 2\sqrt{\log a} \ge 1/\sqrt{\log a}$).}
\end{align*}
Thus, 
\begin{align*}
\ee_{\theta}(S(\theta)-\inf_\psi S(\psi)) &\le 4\varepsilon + \int_0^\infty (2N_d(\varepsilon)e^{-nt^2/(96M+32)} \wedge 1) dt\\
&\le 4\varepsilon + 2\sqrt{\frac{(96M+32)\log (2N_d(\varepsilon))}{n}}.
\end{align*}
This completes the proof of Theorem \ref{lse}. \hfill $\Box$

\section{An inequality for self-bounded functions}
Theorem \ref{hoeffding} can be extended in several ways. The following extension is analogous to existing results for the so called ``self-bounded'' functions. The terminology was introduced by Boucheron, Lugosi \& Massart in \cite{blm00}. Recall the notation of these authors that was described in section \ref{hoeffineq} of Chapter \ref{review}. They call a function ``self-bounded'' if $V_+$ and $V_-$ can be bounded by some function of $Z$, which is usually a linear function. Functions of independent random variables which satisfy the self bounding property appear in reasonable amounts in the literature; examples include suprema of nonnegative empirical processes and conditional Rademacher averages. Further details and references about concentration inequalities for self-bounded functions of independent random variables are available in \cite{blm00, bblm05} and \cite{rio01}.

We shall say a function $f$ is self-bounded (w.r.t.\ $X$) if $v(X) \le Bf(X) + C$ for some constants $B$ and $C$, where $v$ is defined in (\ref{ddef}). 

In Section \ref{perm}, we shall provide some applications of this result to obtain Bernstein type concentration for functionals related to random permutations (matching problem and Spearman's footrule).
\begin{thm}\label{ext1}
Continuing with the notation introduced at the beginning of the Chapter, suppose $B,C$ are finite positive constants such that $v(x) \le Bf(x) + C$ for each $x$. Then $\pp\{|f(X)| \ge t\} \le 2e^{-t^2/(2C + 2Bt)}$ for any $t\ge 0$.
\end{thm}
{\bf Example.} As usual, consider $X = \sum_{i=1}^n Y_i$, where $Y_i$'s are now independent random variables taking value in $[0,1]$. Let $\mu_i = \ee(Y_i)$. We shall use Theorem \ref{ext1} to prove a version of Bernstein's inequality in this setup. As before, let $f(X)=X - \ee (X)$. Let $\xp$ be constructed as in the example following Theorem \ref{weaklaw}. From previous computation, we know that 
\[
v(X) = \frac{1}{2} \sum_{i=1}^n \ee((Y_i^\prime-Y_i)^2|X) = \frac{1}{2}\sum_{i=1}^n (\ee Y_i^2 - 2\mu_i \ee(Y_i|X) + \ee(Y_i^2|X)).
\]
Using $Y_i^2\le Y_i$ (because $0\le Y_i\le 1$), we get
\[
v(X) \le \frac{1}{2}\sum_{i=1}^n (\ee(Y_i) + \ee(Y_i|X)) = \frac{1}{2}(\ee(X) + X).
\]
Thus, taking $B = 1/2$ and $C=\ee(X)$, we get
\[
\pp\{|X-\ee(X)| \ge t\} \le 2e^{-t^2/(2\ee(X)+t)}.
\]
Note, for instance, if $\mu_i = 1/2$ for all $i$, then $E(X)=n/2$, and this bound is essentially equivalent to the Hoeffding bound. The increase in efficiency is apparent only when $\mu_i$'s are small.\\
\\
{\bf Proof of Theorem \ref{ext1}.} Proceeding exactly as in the proof of Theorem \ref{hoeffding}, we get
\begin{align*}
|m^\prime(\theta)| &\le |\theta|\ee(e^{\theta f(X)} v(X)) \\
&\le |\theta|\ee(e^{\theta f(X)} (Bf(X) + C)) \\
&= B|\theta| m^\prime(\theta) + C|\theta| m(\theta).
\end{align*}
Since $m$ is a convex function and $m^\prime(0) 
= \ee(f(X)) = 0$, therefore $m^\prime(\theta)$ always has the same sign as $\theta$. Thus, for $0\le \theta< 1/B$, the above inequality translates into 
\[
\frac{d}{d\theta} \log m(\theta) \le \frac{C\theta}{1-B\theta}.
\]
Using this and recalling that $m(0)=1$, we have
\begin{align*}
\log m(\theta) &\le \int_0^\theta \frac{Cu}{1-Bu}du \le \frac{C\theta^2}{2(1 - B\theta)}.
\end{align*}
Putting $\theta = t/(C+Bt)$, we get
\begin{align*}
\pp\{f(X) \ge t\} &\le \exp(-\theta t + \log m(\theta)) \le e^{-t^2/(2C + 2Bt)}.
\end{align*}
The lower tail can be done similarly (though the argument is not exactly symmetric in this case). \hfill $\Box$

\section{Example: Matching problem and Spearman's footrule}\label{perm}
Suppose $\pi$ is chosen uniformly at random from the set of all permutations of $\{1,\ldots,n\}$. It is a classical probabilistic fact (the matching problem) that as $n\ra \infty$, the distribution of the number of fixed points of $\pi$ converges weakly to the Poisson distribution with mean $1$. Error bounds can also be obtained by various methods. The question of tail bounds for this random variable, with explicit constants, is a reasonable problem to look at. Although this question can presumably be tackled by Talagrand's technique (\cite{talagrand95}, section 5; also discussed in section \ref{gpreview} here) for random permutations, it is a good test case for our theory. 

Another interesting problem related to random permutations is the behavior of $\sum_{i=1}^n|i-\pi(i)|$. This statistic, known as Spearman's footrule (see, e.g. \cite{kendall90}), arises in the nonparametric theory of statistics. 

In the rest of this section, we shall work out concentration inequalities for a generalized version of these statistics, originating in an early work of Hoeffding \cite{hoeffding51}.

Let $\{a_{ij}\}$ be an $n\times n$ array of real numbers, assumed to be in $[0,1]$ for our purposes. Let
$\pi$ be a random (uniform) permutation of $\{1,\ldots,n\}$, and let $X=\sum_{i=1}^n a_{i\pi(i)}$. This class of random variables was first studied by Hoeffding \cite{hoeffding51}, who proved that they are approximately normally distributed under certain conditions. The following proposition gives concentration inequalities in this setup: 
\begin{prop}\label{hoeff3}
Let $\{a_{ij}\}_{1\le i,j\le n}$ be a collection of numbers from $[0,1]$. Let $X = \sum_{i=1}^n a_{i\pi(i)}$, wherfe $\pi$ is drawn from the uniform distribution over the set of all permutations of $\{1,\ldots,n\}$. Then for any $t\ge 0$, $\pp\{|X -\ee(X)| \ge t\} \le
2e^{-t^2/(4\ee(X) + 2t)}$.
\end{prop}
{\bf Remarks.} Note that the bound does not explicitly depend on $n$. This is a consequence of the assumption that the $a_{ij}$'s are bounded. A version of this theorem sans the boundedness assumption can also be proved using our techniques, but then the bound will involve $n$. The classical concentration result of Maurey \cite{maurey79}, stated as Theorem \ref{maurey} in the review section \ref{gpreview} cannot give a Bernstein type bound like the above. Talagrand's theorem (\cite{talagrand95}, Theorem 5.1) might, but it is not clear from the abstract form whether it really does. McDiarmid's corollary \cite{mcdiarmid02} of Talagrand's result certainly cannot give a result like the above.
\medskip

\noindent {\bf Example 1.} (Matching problem.) Taking $a_{ij} = \ii\{i=j\}$, we get $X$ to be the number of fixed points of $\pi$. Since $\ee(X)=1$, we get the exponential tail bound $\pp\{|X-1|\ge t\} \le 2e^{-t^2/(4+2t)}$,
which does not depend on $n$.
\medskip

\noindent{\bf Example 2.} (Spearman's footrule.) Sometimes, the boundedness of the $a_{ij}$'s can be overcome by just dividing by the maximum. We give one such example now. 

In nonparametric statistics, a standard measure of distance between two permutations $\pi$ and $\sigma$ is the Spearman's footrule, defined as
\[
\rho(\pi,\sigma):= \sum_{i=1}^n |\pi(i)-\sigma(i)|.
\]
A standard reference for the uses of Spearman's footrule in nonparametric statistics is the book \cite{kendall90} by Kendall and Gibbons. From a statistical point of view, it is of interest to know the distribution of the Spearman's footrule distance between the identity and a random uniform permutation. Diaconis and Graham \cite{diaconisgraham77} proved the following theorem:
\begin{thm}
\textup{[Diaconis and Graham \cite{diaconisgraham77}]}
Let $\rho(\pi,\sigma) = \sum_{i=1}^n |\pi(i)-\sigma(i)|$. If $\pi$ is chosen uniformly, then
\begin{align*}
\ee(\rho) &= \frac{1}{3}(n^2-1),\\
\var(\rho) &= \frac{1}{45}(n+1)(2n^2+7), \ \text{ and }\\
\pp\biggl\{ \frac{\rho-\ee(\rho)}{\sqrt{\var(\rho)}}\le t\biggr\} &\ra \frac{1}{\sqrt{2\pi}}\int_{-\infty}^t e^{-x^2/2} dx
\end{align*}
as $n \ra \infty$, for each $t\in \rr$.
\end{thm}
We shall use Proposition \ref{hoeff3} to get finite sample tail bounds of the correct order. 
Without loss of generality, we can take $\sigma$ to be the identity permutation. Let 
\[
a_{ij} = \frac{|i-j|}{n}, \ \ \ 1\le i,j \le n.
\] 
Then $0\le a_{ij}\le 1$. Let $X=\sum_{i=1}^n a_{i\pi(i)} = \rho/n$. Then by Proposition \ref{hoeff3},
$\pp\{|X-\ee(X)|\ge u\} \le e^{-u^2/(4\ee(X)+2u)}$. Now put $u = t\sqrt{\var(X)}$, and observe that 
\[
\frac{X-\ee(X)}{\sqrt{\var(X)}} = \frac{\rho-\ee(\rho)}{\sqrt{\var(\rho)}}.
\]
By the Diaconis-Graham computation, $\ee(X)= (n^2-1)/3n$ and  $\var(X) = (n+1)(2+7n^{-2})/45$. Combining, we get the following tail bound:
\begin{prop}
Let $\rho$ be the Spearman's footrule distance between the identity and a uniformly chosen permutation of $\{1,\ldots,n\}$. For any $t\ge 0$, we have
\[
\pp\biggl\{\biggl|\frac{\rho-\ee(\rho)}{\sqrt{\var{\rho}}}\biggr|\ge t\biggr\} \le 2e^{-t^2/(\alpha_n +\beta_n t)},
\]
where
\begin{align*}
\alpha_n &= \frac{60(n-1)}{n(2+7n^{-2})} \ra 30 \ \text{ as } n\ra \infty, \ \text{ and }\\
\beta_n &= 2\sqrt{45}(n+1)^{-1/2}(2 + 7n^{-2})^{-1/2} \ra 0 \ \text{ as } n \ra \infty.
\end{align*}
\end{prop}
Note that our tail bound is approximately a gaussian bound (and hence, of the correct order) for large $n$. However, the constants are poor. The reason is that Proposition~\ref{hoeff3} is designed for random variables that exhibit Poissionian behavior. Our $X$ in this problem does not belong to that class, particularly because its variance is smaller than its mean by a factor significantly smaller than unity. Thus, although Proposition \ref{hoeff3} can be successfully applied to this problem to get concentration bounds of the correct order, it is not an ideal example.
\\
\\
\noindent {\bf Proof of Proposition \ref{hoeff3}.} Construct $\xp$ as follows: Choose $I,J$ uniformly and independently at random from $\{1,\ldots, n\}$. Let $\pi^\prime = \pi\circ (I,J)$, where $(I,J)$ denotes the transposition of $I$ and $J$. It can be easily verfied that $(\pi,\pi^\prime)$ is an exchangeable pair. Hence if we let 
\[
\xp := \sum_{i=1}^n a_{i\pi^\prime(i)},
\]
then $(X,\xp)$ is also an exchangeable pair. Now note that
\begin{align*}
\frac{1}{2}\ee(n(X - \xp)|\pi) &= \frac{n}{2}\ee(a_{I\pi(I)} + a_{J\pi(J)} - a_{I\pi(J)} - a_{J\pi(I)}|\pi) \\
&= \frac{1}{n} \sum_{i,j} a_{i\pi(i)} - \frac{1}{n}\sum_{i,j} a_{i\pi(j)} \\
&= X-\ee(X).
\end{align*}
Thus, we can take $f(x)= x-\ee(X)$ and $F(x,y) = \frac{1}{2}n(x-y)$. Now note that since $0\le a_{ij}\le 1$ for
all $i$ and $j$, we have
\begin{align*}
v(X) &= \frac{n}{4}\ee((X-\xp)^2|\pi) \\
&=\frac{1}{4n}\sum_{i,j} (a_{i\pi(i)} +
a_{j\pi(j)}-a_{i\pi(j)}-a_{j\pi(i)})^2 \\
&\le \frac{1}{2n} \sum_{i,j}(a_{i\pi(i)} +
a_{j\pi(j)}+a_{i\pi(j)}+a_{j\pi(i)}) \\
&= X + \ee(X) = f(X) + 2\ee(X).
\end{align*}
Applying Theorem \ref{ext1} with $B = 1$ and $C= 2\ee(X)$ completes
the proof. \hfill $\Box$

\section{Inequalities for unbounded differences}\label{unbdd}
In this section, we present a some results which are applicable when $v(X)$ is unbounded. As usual, we shall supplement with trivial examples. A nontrivial application will be worked out in section \ref{group} of the next chapter.

The first result of this section is an extension of Theorem \ref{hoeffding} which only requires reasonable bounds on the moment generating function of $v(X)$, assuming that it exists. To be more specific, note that if $v(X)$ is bounded, then
\[
\lim_{L\ra \infty} L^{-1} \log \ee(e^{Lv(X)}) = \|v(X)\|_\infty.
\] 
Thus, the number $r(L)$ defined as
\[
r(L):=L^{-1}\log \ee(e^{Lv(X)})
\]
may serve as a surrogate for an actual bound on $v(X)$, in situations where such a bound does not exist, or is not representative of the true size of $v(X)$. The following theorem allows us the flexibility of using $r(L)$ with appropriate choice of $L$:
\begin{thm}\label{bernstein}
Suppose $r(L)$ is defined as above. Fix any $L > 0$. Then we have $\pp\{|f(X)|\ge t\} \le 2e^{-t^2/(2r(L) + 4tL^{-1/2})}$ for any $t \ge 0$.
\end{thm}
{\bf Remarks.} The idea is to choose $L$ so that $L^{-1/2}\ll r(L)$. In particular, observe that if $v(X)\le C$ almost surely for some constant $C$, we can take $L\ra \infty$ and get the bound in Theorem \ref{hoeffding}. 
\medskip

\noindent {\bf Example.} Again, let $X = \sum_{i=1}^n Y_i$, where $Y_i$'s are independent zero mean random variables. However, we shall now drop the boundedness assumption, and assume only that $Y_i$'s have gaussian tails: In other words, assume that there exists $\theta > 0$ such that $\ee(e^{\theta Y_i^2}) \le K_t < \infty$ for each $i$ for some fixed constant $K_\theta$. Choosing $L=2\theta$ and applying Jensen's inequality, we have
\begin{align*}
r(L) &\le (2\theta)^{-1} \log \ee(e^{\theta\sum_{i=1}^n (\ee(Y_i^2)+Y_i^2)})\\
&\le (2\theta)^{-1}\sum_{i=1}^n 2 \log \ee(e^{\theta Y_i^2})\le \theta^{-1}n \log K_\theta.
\end{align*}
The above result now gives the bound
\[
\pp\{|X| \ge t\} \le 2\exp\biggl(-\frac{t^2}{2\theta^{-1} n\log K_\theta + 4t (2\theta)^{-1/2}}\biggr).
\]
Note how this reduces to the Hoeffding bound when the $Y_i$'s are bounded.\\
\\
{\bf Proof of Theorem \ref{bernstein}.} Let $u(X) = e^{\theta f(X)}/m(\theta)$, where $m(\theta)=\ee(e^{\theta f(X)})$, as usual. Applying Jensen's inequality, we have for $\theta \ge 0$
\begin{align*}
m^\prime(\theta) &\le \theta \ee(e^{\theta f(X)} v(X)) \\
&= L^{-1}\theta m(\theta)\ee\biggl[u(X)\biggl(\log \frac{e^{Lv(X)}}{u(X)} + \log u(X)\biggr)\biggr]\\
&\le L^{-1} \theta m(\theta) \log \ee(e^{Lv(X)}) + L^{-1}\theta \ee(e^{\theta f(X)} \log u(X)).
\end{align*}
Recall that $m^\prime(0)=\ee(f(X))=0$, $m(0)=1$, and $m$ is a convex function. Hence $m(\theta)\ge 1$ for all $\theta$. Consequently, $\log u(X) \le \theta f(X)$. Using this in the above bound, we have
\[
m^\prime(\theta) \le r(L)\theta m(\theta) + L^{-1} \theta^2 m^\prime(\theta).
\]
Written differently, this gives
\[
\frac{d}{d\theta}\log m(\theta) \le \frac{r(L)\theta}{1-L^{-1}\theta^2}
\]
for $0\le \theta < L^{1/2}$. Now take any $t \ge 0$, and let $\theta$ be the positive root of the equation 
\[
\frac{r(L)\theta}{1-L^{-1}\theta^2} = t.
\]
Explicitly, we have
\begin{equation}\label{theta}
\theta = \frac{Lr(L)}{2t}\biggl(\sqrt{1+\frac{4t^2}{r(L)^2L}} -1\biggr).
\end{equation}
Using the inequality $\sqrt{1+\alpha} \le 1+\sqrt{\alpha}$, it is easy to see that $\theta < L^{1/2}$. Another application of the same inequality gives
\begin{align*}
\sqrt{1+\alpha} - 1 - \frac{\alpha}{2(1+\sqrt{\alpha})} &\ge \sqrt{1+\alpha} - 1 - \frac{\alpha}{2\sqrt{1+\alpha}}\\
&=
\frac{(\sqrt{1+\alpha} - 1)^2}{2\sqrt{1+\alpha}} \ge 0.
\end{align*}
Using this inequality in (\ref{theta}), we get
\[
\theta \ge \frac{t}{r(L) + 2tL^{-1/2}}.
\]
Now, with the $\theta$ defined in (\ref{theta}) we have
\[
\log m(\theta) \le \int_0^\theta \frac{r(L)u}{1-L^{-1}u^2} du \le \frac{r(L)\theta^2}{2(1-L^{-1}\theta^2)} = \frac{\theta t}{2}.
\]
Thus, by Chebychev's inequality we get
\begin{align*}
\pp\{f(X) \ge t\} &\le e^{-\theta t + \log m(\theta)} \le e^{-\theta t/2}.
\end{align*}
Using the lower bound on $\theta$ derived above, we get the desired expression.
The lower tail bounds are obtained by symmetry.\hfill $\Box$\\

Our next result is the exchangeable pairs version of the Burkholder-Davis-Gundy inequality \cite{burkholder73, burkholder88, burkholder89} from classical probability. In its simplest form, this inequality says that for a martingale difference sequence $\{X_i\}_{1\le i\le n}$ adapted to some filtration $\{\mf_i\}_{1\le i\le n}$, we have for any $p\ge 1$ the inequality
\[
\ee\bigl|\sum_{i=1}^nX_i\bigr|^{2p} \le (2p-1)^p \ee\bigl|\sum_{i=1}^n \ee(X_i^2|\mf_{i-1})\bigr|^p.
\]
Note that the inequality is trivially an equality for $p=1$. Essentially, this inequality gives a $p^{\mathrm{th}}$ moment expression of the notion that a martingale which is the sum of a homogeneous sequence of differences grows like $n^{1/2}$. (Analogously, the Hoeffding inequality gives a tail estimate expression of this notion when the differences are bounded.) This is clear from the observation that the right hand side above is bounded by $n^{p-1} \sum_{i=1}^n  \ee|X_i|^{2p} = O(n^p)$ by Jensen's inequality, and so if the $X_i$'s are of comparable sizes then $\ee|\sum X_i|^{2p} = O(n^p)$, which ``shows '' that $\sum X_i = O(n^{1/2})$ upto $p^{\mathrm{th}}$ moment accuracy. 

Recently, Boucheron, Bousquet, Lugosi \& Massart \cite{bblm05} have derived a useful version of this inequality for general functions of independent random variables (rather than just sums). Moment bounds are often useful when we can take $p$ to be large (growing with $n$), because the resulting Chebychev inequalities give surprisingly efficient tail bounds. In fact, it is an easy fact that for any suitable random variable $X$ and any $t\ge 0$,
\[
\pp\{|X|\ge t\} \le \inf_{p\ge 1} \frac{\ee|X|^p}{t^p} \le \inf_{\theta\ge 0} \frac{\ee(e^{\theta|X|})}{e^{\theta t}},
\]
which shows that optimized Chebychev bounds are better than optimized Chernoff bounds. For applications of the moment inequalities to obtain efficient tail bounds for a variety of complicated functions of independent random variables, we refer to~\cite{bblm05}. 

We shall now present the exchangeable pairs version of the Burkholder-Davis-Gundy inequality.
In the following, $\|\cdot\|_p$ will denote the $L^p$ norm of random variables; that is, for a random variables $Y$, $\|Y\|_p := (\ee |Y|^p)^{1/p}$. 
\begin{thm}\label{moments}
For any positive integer $p$, we have
\[
\|f(X)\|_{2p}^2 \le (2p-1) \|v(X)\|_p.
\]
\end{thm}
{\bf Example.} 
To see why this can be called the Burkholder-Davis-Gundy inequality for exchangeable pairs, consider, as usual, sums of independent random variables. Let $Y_1,\ldots,Y_n$ be independent mean zero random variables with finite $p^{\mathrm{th}}$ moment, where $p$ is a fixed positive integer. Let $X=\sum_{i=1}^n Y_i$. Construct $\xp$ as in the previous sections, and recall that $f(X)=X$, $F(X,\xp)=n(X-\xp)$, and
\[
v(X) = \frac{1}{2}\sum_{i=1}^n (\ee(Y_i^2) + \ee(Y_i^2|X)).
\]
By the above theorem and simple applications of Minkowski and Jensen inequalities, we get
\[
\|X\|_{2p}^2 \le (2p-1)\bigl\|\sum_{i=1}^n Y_i^2\bigr\|_p,
\]
which is the classical inequality for sums of independent random variables. \\
\\
{\bf Proof of Theorem \ref{moments}.} By Lemma \ref{fundamental}, we have
\[
\ee(f(X)^{2p}) = \frac{1}{2} \ee((f(X)^{2p-1} - f(\xp)^{2p-1})F(X,\xp)).
\]
By the inequality 
\[
|x^{2p-1} - y^{2p-1}| \le \frac{2p-1}{2}(x^{2p-2} + y^{2p-2})|x-y|
\]
which follows easily from a convexity argument very similar to (\ref{expin}), we have
\begin{align*}
\ee(f(X)^{2p}) &\le (2p-1)\ee(f(X)^{2p-2} v(X))
\end{align*}
By H\"{o}lder's inequality, we get
\[
\ee(f(X)^{2p}) \le (2p-1)(\ee(f(X)^{2p}))^{(p-1)/p}(\ee(v(X)^p))^{1/p}.
\]
The proof is completed by transferring $\ee(f(X)^{2p})^{(p-1)/p}$ to the other side.\hfill $\Box$

\section{A refinement of the exponential inequality}\label{refine}
The following is a refinement of Theorem \ref{hoeffding} which allows us, in particular, to derive concentration bounds in the Sherrington-Kirkpatrick model of spin glasses in the next section. 
\begin{thm}\label{ext2}
Define 
\begin{align*}
v_1(x) &:= \frac{1}{2}\ee\bigl((f(X)-f(\xp))F(X,\xp)\bigl| X = x\bigr) \ \ \text{ and} \\
v_2(x) &:= \frac{1}{4}\ee\bigl((f(X)-f(\xp))^2 |F(X,\xp)|\bigl| X = x\bigr).
\end{align*}
Suppose $C$ and
$\varepsilon$ are constants such that $|v_1(X)| \le C$ and $v_2(X) \le
\varepsilon$ a.s. Then $\pp\{|f(X)| \ge t\} \le
2e^{-Ct^2/(2C^2 + 8\varepsilon t)}$ for any $t\ge 0$.
\end{thm}
{\bf Remark.} Note the difference between $v$ and $v_1$: the latter has no absolute value inside. This  can be considered as a ``second order refinement'' of the original inequality in Theorem \ref{hoeffding}.
\vskip.1in

\noindent{\bf Proof.} Recall the identity 
\[
m^\prime(\theta) = \frac{1}{2} \ee\bigl( (e^{\theta f(X)} - e^{\theta
  f(\xp)}) F(X,\xp)\bigr)
\]
that was derived using Lemma \ref{fundamental} in the proof of Theorem \ref{hoeffding}. Instead of simply using the inequality $|e^x-e^y|\le \frac{1}{2}(e^x + e^y)|x-y|$, we now take a finer recourse: For any $x\le y$, note that
\begin{align*}
0 &\le \frac{1}{2}(e^y + e^x)(y-x) - (e^y - e^x) \\
&= \int_x^y \biggl( \frac{e^y - e^x}{y-x}(u-x) + e^x - e^u\biggr) du \\
&\le \int_x^y \biggl(\frac{1}{2}(e^y + e^x)(u-x) + 0\biggr) du \\
&= \frac{1}{4}(e^y + e^x)(y-x)^2.
\end{align*}
Thus, we have 
\[
e^x - e^y = \frac{1}{2}(e^x + e^y)(x-y) + \delta(x,y),
\]
where
\[
|\delta(x,y)| \le
\frac{1}{4} (e^x + e^y)(x-y)^2.
\]
Combining this with the observation that $(x-y)F(x,y)= (y-x)F(y,x)$, is it not difficult to deduce that for any $\theta \ge 0$ we have
\begin{align*}
m^\prime(\theta) &\le \theta \ee(e^{\theta f(X)}v_1(X)) +
\theta^2 \ee(e^{\theta f(X)} v_2(X)) \\
&\le (C\theta + \varepsilon \theta^2) m(\theta).
\end{align*}
This gives, for any $\theta\ge 0$,
\[
\log m(\theta) \le \frac{C\theta^2}{2} + \frac{\varepsilon \theta^3}{3}.
\]
Now fix any $t\ge 0$. Let $\theta$ be the positive solution of the equation
$\varepsilon \theta^2 + C\theta = t$. Explicitly,
\[
\theta = \frac{C}{2\varepsilon}\biggl(\sqrt{1+\frac{4\varepsilon t}{C^2}} -
1\biggr). 
\]
Then
\begin{align*}
-\theta t + \log m(\theta) &\le -\theta(\varepsilon \theta^2 + C\theta) + \frac{C\theta^2}{2} + \frac{\varepsilon \theta^3}{3} \\
&= -\frac{C\theta^2}{2} - \frac{2\varepsilon \theta^3}{3} \le
-\frac{C\theta^2}{2}. 
\end{align*}
Now note that for any $\alpha \ge 0$, 
\[
\sqrt{1+\alpha} - 1 - \frac{\alpha}{2\sqrt{1+\alpha}} =
\frac{(\sqrt{1+\alpha} - 1)^2}{2\sqrt{1+\alpha}} \ge 0.
\]
Thus,
\[
\theta \ge \frac{C}{2\varepsilon} \biggl(\frac{\frac{4\varepsilon t}{C^2}}{2\sqrt{
    1+\frac{4\varepsilon t}{C^2}}}\biggr) = \frac{t}{\sqrt{C^2
    + 4\varepsilon t}}.
\]
Combining the steps, we get
\begin{align*}
\pp\{f(X) \ge t\} &\le e^{-\theta t + \log m(\theta)} \le e^{-C\theta^2/2}
\le e^{-Ct^2/(2C^2 + 8\varepsilon t)}. 
\end{align*}
The lower tail bound follows by symmetry. \hfill $\Box$

\section{Application to spin glasses}\label{spinglass}
The Sherrington-Kirkpatrick (S-K) model of spin glasses considers the Hamiltonian
\[
H(\sigma) = - n^{-1/2}\sum_{1\le i < j\le n} g_{ij} \sigma_i \sigma_j -  h\sum_{i=1}^n \sigma_i,
\]
where $g_{ij}$'s are a fixed realization of a collection of i.i.d.\ standard gaussian random variables. For recent advances in the rigorous analysis of this model, we refer to Chapter 2 of the book \cite{talagrand03} by Michel Talagrand, and also his proof \cite{talagrand05} of the Parisi formula for the limiting free energy in this model. For the physical perspective, one should look at the book \cite{mpv87} by M\'ezard, Parisi and Virasoro.

In spite of the fact that the Parisi formula has been proved, very little is explicitly known about the low temperature phase of the S-K model at the time of writing this thesis. In fact, as Talagrand admits (\cite{talagrand03}, p.\ 182), ``we do not even know where to start''! 

The mean field equations discussed in section \ref{meanmodels} are no longer valid in the S-K model, even at high temperatures. Instead, physicists think that a modification of the naive mean field equations (\ref{meanfield}), the so-called TAP equations, hold for the S-K model. These are as follows:
\begin{equation}\label{TAP}
``\langle \sigma_i\rangle  =  \tanh  \biggl(\frac{\beta}{\sqrt{n}} \sum_{j\le n, j\ne i} g_{ij} \langle \sigma_j\rangle + h - \beta^2(1-q)\langle \sigma_i\rangle\biggr),\ \ i=1,\ldots,n\text{''}
\end{equation}
where $q$ solves $q = \ee \tanh^2(\beta Z\sqrt{q} + h)$, $Z$ being a standard gaussian random variable. The validity of this conjecture has been established  rigorously for sufficiently small $\beta$ by Talagrand (\cite{talagrand03}, Theorem 2.4.20). 

In this section, we shall prove that another set of mean field type equations are valid, irrespective of the temperature, in a class of models which will encompass the S-K model. For instance, we shall be able to show that in this class of models, the magnetization $m(\sigma) = \frac{1}{n}\sum_{i=1}^n \sigma_i$ satisfies
\begin{equation}\label{skmag}
m(\sigma) \approx \frac{1}{n}\sum_{i=1}^n \tanh(\beta m_i(\sigma)) \ \ \text{with high probability,}
\end{equation}
where $m_i(\sigma) := n^{-1/2}\sum_{j\le n, j\ne i } g_{ij}\sigma_j + h$ is the ``local field'' at site $i$. It will follow from the same argument that whenever $k \ll \sqrt{n}$, we also have
\begin{equation}\label{skover}
\frac{1}{n}\sum_{i=1}^n \prod_{r=1}^k\sigma_i^r \approx \frac{1}{n}\sum_{i=1}^n \prod_{r=1}^k\tanh(\beta m_i(\sigma^r)),
\end{equation}
where $\sigma^1,\ldots,\sigma^k$ are i.i.d.\ from the Gibbs measure. The left hand side, called the ``overlap'' of the $k$ configurations, is a particularly interesting quantity in the S-K model when $k=2$. 

However, we must admit that we have found no use for all this information till now. Still, it may be interesting just because the low temperature phase of the S-K model is a highly intractable object.

Our results will be valid whenever the $L^2$ operator norms of the matrices $J$ and $J_2 := (J_{ij}^2)$ are bounded. Recall that the $L^2$ operator norm of a matrix $A$ is defined as
\[
\|A\| := \max_{\|x\|=1} \|Ax\|.
\]
Alternatively, it is the square root of the maximum eigenvalue of $A^TA$. In the S-K model, it is well-known that the norm of $J$ is bounded in probability. In fact, it is known that with $J^{(n)} = (n^{-1/2}g_{ij})_{1\le i,j\le n}$, where $\{g_{ij}\}_{1\le i\le j\le n}$ are i.i.d.\ standard gaussian random variables and $g_{ji}=g_{ij}$, we have
\[
\|J^{(n)}\| \ra 2 \ \ \text{in probability.}
\]
A proof of this result, as well as tail bounds, can be found in the survey article by Bai \cite{bai99}, section 2.2.
It is also easy to see that the norm of $J_2$ is bounded in probability in the S-K model, because 
\[
\|J_2\| \le \bigl[ n\max_{1\le i\le n} \sum_{j=1}^n J_{ij}^4\bigr]^{1/2} = \biggl[\max_{1\le i\le n} \sum_{j=1}^n \frac{g_{ij}^4}{n}\biggr]^{1/2} \ra \ee(g_{11}^4) \ \ \text{in probability,}
\]
since for each $i$, $\frac{1}{n}\sum_{j=1}^n g_{ij}^4$ is concentrated around $\ee(g_{11}^4)$ and the tails fall off sharply enough for the above to hold (a fact that can be easily proved using the Burkholder-Davis-Gundy inequality).
The key result of this section is the following:
\begin{thm}\label{skmodel}
Let $J$ be the interaction matrix in a model described by the Hamiltonian (\ref{quadmodel}). Let $J_2 = (J_{ij}^2)_{1\le i,j\le n}$. Take any $\alpha=(\alpha_1,\ldots,\alpha_n)\in \rr^n$. Fix $\beta \ge 0$,  and let 
\begin{align*}
C &= C(\alpha, \beta) := 2(1+\beta\|J\|+\beta^2 \|J_2\|)\|\alpha\|^2, \ \ \text{and}\\
\varepsilon &= \varepsilon(\alpha, \beta) := 4(\max_{1\le i\le n}|\alpha_i|)\bigl[ 1+ (\beta\|J\|+\beta^2\|J_2\|)^2\bigr] \|\alpha\|^2.
\end{align*}
For each $i$, let $m_i = m_i(\sigma) := \sum_{j=1}^n J_{ij}\sigma_j + h_i$. Let 
\[
Y = \sum_{i=1}^n \alpha_i(\sigma_i - \tanh(\beta m_i)).
\] 
If $\sigma$ is drawn from the Gibbs measure at inverse temperature $\beta$, then $\ee(Y)=0$, $\var(Y)\le C$, and for any $t\ge 0$, $\pp\{|Y|\ge t\} \le 2e^{-Ct^2/(2C^2 + 8\varepsilon t)}$.
\end{thm}
{\bf Applications.} Taking $\alpha_i = 1/n$ for each $i$, we get $\|\alpha\|^2 = 1/n$. Fix $\beta \ge 0$, and let $K = \beta\|J\| + \beta^2\|J_2\|$. As observed before, $K= O(1)$ in the S-K model. We have $C = (2+2K)/n$ and $\varepsilon = 4(1+K^2)/n^2$. Thus, if we let $m(\sigma) := \frac{1}{n}\sum_{i=1}^n \sigma_i$ be the magnetization of a configuration drawn from the Gibbs measure at inverse temperature $\beta$, and let
\[
Y = m(\sigma) - \frac{1}{n}\sum_{i=1}^n \tanh(\beta m_i(\sigma)),
\]
where, as usual, $m_i(\sigma) = \sum_{j=1}^n J_{ij} \sigma_j+ h_i$ is the local field at $i$, then $\ee(Y^2) \le (2+2K)/n$ and 
\begin{equation}\label{overlap}
\pp\{|Y|\ge t\} \le 2e^{-nt^2/(a+bt)},
\end{equation}
where $a = 4 + 4K$ and $b = 16(1+K^2)/(1+K)$ are free of $n$. This shows (\ref{skmag}).

The relation (\ref{skover}) for the $k^{\mathrm{th}}$-order overlaps can be treated exactly in the same way, by successively conditioning on $(\sigma^1,\ldots,\sigma^{r-1}, \sigma^{r+1},\ldots,\sigma^k)$, $r=1,\ldots,k$, and getting
\[
\frac{1}{n}\sum_{i=1}^n \prod_{s=1}^{r-1} \tanh(\beta m_i(\sigma^{s-1}))\prod_{s=r}^k \sigma_i^s
\approx \frac{1}{n}\sum_{i=1}^n \prod_{s=1}^r \tanh(\beta m_i(\sigma^{s-1})) \prod_{s=r+1}^k \sigma_i^s
\]
at the $r^{\mathrm{th}}$ stage. Roughly, this can be carried out till $k = o(\sqrt{n})$, since the errors accrued at each stage are like $n^{-1/2}$.
\vskip.3in

\noindent {\bf Proof of Theorem \ref{skmodel}.}
As in the proof of Proposition \ref{curie}, we construct an exchangeable pair by taking a step in the Gibbs sampler chain as follows: First, draw $\sigma$ from the Gibbs measure. Next, choose a coordinate $I$ uniformly at random, and replace the $I^{\mathrm{th}}$ coordinate of $\sigma$ by a sample drawn from the conditional distribution of the $\sigma_I$ given $\{\sigma_j, j\ne I\}$. Let 
\[
F(\sigma, \sss) := n \sum_{i=1}^n \alpha_i(\sigma_i - \sss_i) = n(\sigma_I-\sss_I).
\]
Then $F$ is antisymmetric and
\begin{align*}
f(\sigma) = \ee(F(\sigma, \sss)|\sigma) &= \sum_{i=1}^n \alpha_i(\sigma_i - \ee(\sigma_i|\{\sigma_j, j\ne i\}) \\
&= \sum_{i=1}^n \alpha_i(\sigma_i - \tanh(\beta m_i(\sigma))).
\end{align*}
For any $\sigma \in \{-1,1\}^n$ and $1\le j\le n$, let
\[
\sigma^{(j)} := (\sigma_1,\ldots, \sigma_{i-1}, -\sigma_i, \sigma_{i+1}, \ldots, \sigma_n).
\]
Define $h:\rr \ra\rr$ as $h(x) := \tanh(\beta x)$, and for each $i, j$, let
\[
b_{ij} = b_{ij}(\sigma) := h(m_i(\sigma)) - h(m_i(\sigma^{(j)})),
\]
where $m_i(\sigma) = \sum_{j=1}^n J_{ij}\sigma_j + h_i$, as defined in the statement of the theorem.
Now note that
\begin{equation}\label{fF}
f(\sigma) - f(\sigma^{(j)}) = 2\alpha_j\sigma_j -
\sum_{i=1}^n \alpha_i b_{ij}, \ \text{ and } \ F(\sigma,
\sigma^{(j)}) = 2n\alpha_j\sigma_j. 
\end{equation}
Let $p_j = p_j(\sigma)=\pp\{\sigma_j^\prime = -\sigma_j \mid \sigma, I=j\}$. Combining everything, we have
\begin{equation}\label{d1exp}
v_1(\sigma) = \frac{1}{2}\ee\bigl((f(\sigma) - f(\sss))F(\sigma, \sss)\bigl|\sigma\bigr) =
2\sum_{j=1}^n \alpha_j^2p_j - \sum_{i,j}
\alpha_ib_{ij}\alpha_j \sigma_j p_j.
\end{equation}
It is easy to verify that $\|h^{\prime\prime}\|_\infty \le \beta^2$.
Therefore,
\begin{align*}
&\bigl|h(m_i(\sigma)) - h(m_i(\sigma^{(j)})) - (m_i(\sigma) -
m_i(\sigma^{(j)})) h^\prime(m_i(\sigma))\bigr| \\
&\le
\frac{\beta^2}{2}(m_i(\sigma) - m_i(\sigma^{(j)}))^2.
\end{align*}
Now let $c_i = c_i(\sigma) := h^\prime(m_i(\sigma))$, and note that
$m_i(\sigma) - m_i(\sigma^{(j)}) = 2J_{ij}\sigma_j$. Combining, we
have
\begin{equation}\label{btoa}
|b_{ij} - 2J_{ij}\sigma_j c_i| \le 2\beta^2J_{ij}^2.
\end{equation}
Finally, note that $|c_i|\le \beta$. Using all this information, we get, for any $x,y \in \rr^n$, 
\begin{align}\label{norm}
\biggl|\sum_{i,j}  x_iy_jb_{ij}\biggr| &\le \biggl|\sum_{i,j}x_iy_j(2J_{ij}\sigma_jc_i)\biggr| +
\biggl|\sum_{i,j}x_iy_j (b_{ij}-2J_{ij}\sigma_jc_i)\biggr| \nonumber \\
&\le 2 \|J\|\bigl(\sum_{i=1}^n (x_ic_i)^2\bigr)^{1/2} \bigl(\sum_{j=1}^n (y_j\sigma_j)^2\bigr)^{1/2}
+ \sum_{i,j}|x_iy_j| 2\beta^2J_{ij}^2 \nonumber \\
&\le (2\beta\|J\| +2\beta^2\|J_2\|)\|x\|\|y\|.
\end{align}
Let $K= 2\beta\|J\| +2\beta^2\|J_2\|$. Using the above inequality in equation (\ref{d1exp}) with $x_i = \alpha_i$ and $y_j=\alpha_j\sigma_jp_j$, we get
\[
|v_1(\sigma)| \le (2+K)\sum_{i=1}^n \alpha_i^2.
\]
The first assertion of this theorem now follows by Theorem
\ref{weaklaw}. Next, we use (\ref{fF}) once again (and the inequality $(a+b)^2 \le 2a^2 + 2b^2$) to get 
\begin{align}\label{v2}
v_2(\sigma) &= \frac{1}{4}\ee\bigl((f(\sigma) - f(\sss))^2|F(\sigma,\sss)| \bigl|
\sigma\bigr) \nonumber \\
&\le 4\sum_{i=1}^n |\alpha_i|^3 + \sum_{j=1}^n
|\alpha_j|\bigl(\sum_{i=1}^n \alpha_i b_{ij}\bigr)^2 \nonumber \\
&\le (\max_i|\alpha_i|)\biggl[4\sum_{i=1}^n \alpha_i^2 + \sum_{j=1}^n\bigl(\sum_{i=1}^n \alpha_i b_{ij}\bigr)^2\biggr].
\end{align}
Now, put $x_i = \alpha_i$ and $y_j = \sum_{k=1}^n \alpha_k b_{kj}$, $i,j = 1,\ldots,n$. Then
\begin{align*}
&\|y\|^2 = \sum_{j=1}^n\bigl(\sum_{i=1}^n \alpha_i b_{ij}\bigr)^2 = \sum_{i,j} x_i y_j b_{ij}\\
&\le K \|x\|\|y\| \ \ \text{ by inequality (\ref{norm}).}
\end{align*}
Thus, $\|y\| \le K\|x\|$, and therefore by (\ref{v2}) we have
\[
v_2(\sigma) \le 2(\max_i|\alpha_i|) ( 1+ K^2) \sum_{i=1}^n \alpha_i^2.
\] 
The proof is now completed by applying Theorem \ref{ext2}.  \hfill $\Box$\\

\chapter{Theory and examples: Part II}\label{coupling}
We begin this chapter with Stein's observation \cite{stein86} that an exchangeable pair $(X,\xp)$ automatically defines a reversible Markov kernel $P$ as
\begin{equation}\label{kernel}
Pf(x):=\ee(f(\xp)|X= x),
\end{equation}
where $f$ is any function such that $\ee|f(X)|<\infty$. All other notation will be the same as what was defined at the beginning of Chapter \ref{results}.

The main purpose of the subsequent discussion is to connect the concentration properties of the distribution of $X$ with the rate of convergence to stationarity of a Markov chain following the kernel $P$. In particular, information about the rate of decay of $P^k f(x) - P^k f(y)$, where $(x,y)$ is any point in the support of $(X,\xp)$, can give us a bound on $F(X,\xp)$ in a way that we are going to describe in the following pages. This is a direction that has not been systematically explored in the concentration literature, to the best of our knowledge. We shall also describe a coupling technique for getting a handle on the antisymmetric function $F$ when $f$ is known. 

The techniques so developed, will be applied in section \ref{dobrushin} to get concentration under a famous ``weak dependence'' condition from statistical physics. Examples from classical spin systems at high temperature and random proper graph colorings are worked out in the subsequent sections.

In section \ref{group}, we shall use the results of this section in another direction to obtain concentration of Haar measures on compact groups based on convergence rates of random walks. The result will be used in section \ref{freeprob} to derive a quantitative version of Voiculescu's connection \cite{voiculescu91} between random matrices and free probability theory.

\section{Explicit construction of $F(X,\xp)$}
Our construction of $F$ will involve 
the {\it Poisson equation} associated with the kernel $P$, which can be written as follows:
\begin{equation}\label{poisson}
g-Pg = f.
\end{equation}
Here $f$ is a given function, and the objective is to solve for $g$. The Poisson equation is an object of deep mathematical significance. Its importance in the theory of Markov chains was realized after the work of Paul Meyer \cite{meyer71} and the subsequent investigation of Neveu~\cite{neveu72}. The contributions by Nummelin \cite{nummelin91} are also significant. A classical textbook reference is the book by Revuz \cite{revuz84}. Poisson's equation has been used in the probability literature by too many authors to mention in this short space; for a recent survey of the literature about Poisson's equation in the context of discrete Markov chains, we refer to Makowski and Shwartz \cite{makowskishwartz02}.

The way we can think of using the solution to Poisson's equation in our problem is the following:
If we let $F(x,y) = g(x)-g(y)$, where $g$ is a solution to (\ref{poisson}), then $F$ is antisymmetric and 
\[
\ee(F(X,\xp)|X) = g(X)-\ee(g(\xp)|X) = g(X)-Pg(X)= f(X).
\]
Generally, the solution to Poisson's equation is given by
\[
g = \sum_{k=0}^\infty P^k f,
\]
but problems with convergence are not uncommon. The solution need not be unique, either. Various conditions have been proposed for the existence and uniqueness of solutions of Poisson's equation in various situations. For a summary of such results, we refer to \cite{makowskishwartz02}.

We shall not need the weakest conditions for the existence of solutions to Poisson's equation in any of our applications. The following lemma formalizes our construction of $F$ under generous assumptions, which are satisfied for geometrically ergodic Markov chains, and all of our examples belong to that class.
\begin{lmm}\label{inv1}
Let $f:\xx \ra \rr$ be a measurable function such that $\ee f(X)=0$. Suppose there is a finite constant $L$ such that 
\begin{equation}\label{sum1}
\sum_{k=0}^\infty|P^k f(x)-P^k f(y)|\le L \ \ \text{for every $x$ and $y$.}
\end{equation} 
Then the function 
\[
F(x,y) := \sum_{k=0}^\infty (P^k f(x)-P^k f(y))
\]
satisfies $F(X,\xp)=-F(\xp,X)$ and $\ee(F(X,\xp)|X)=f(X)$.
\end{lmm}
{\bf Proof.} The convergence of the series defining $F$ follows from the summability assumption (\ref{sum1}). Also, $F$ is clearly antisymmetric. Now, 
\[
\ee(P^k f(X) - P^k f(\xp)|X) = P^k f(X) - P^{k+1}f(X).
\]
Thus, for any $N$, we have
\[
\sum_{k=0}^N \ee(P^kf(X)-P^k f(\xp)|X) = f(X) - P^{N+1}f(X).
\] 
The condition (\ref{sum1}) ensures that the above partial sums converge everywhere. Consequently, the sequence $\{P^{N+1}f(X)\}_{N\ge 0}$ also converges everywhere. Again, condition (\ref{sum1}) implies that the limit is a constant, since for any $x$ and $y$, $P^k f(x) - P^k f(y) \ra 0$ as $k \ra \infty$. Since $\ee F(X,\xp) =0 = \ee f(X)$, this constant can only be zero. \hfill $\Box$
\vskip.3in

Although Lemma \ref{inv1} gives an explicit expression for $F$, it is  somewhat inconvenient for practical purposes. We shall now give a coupling version of Lemma \ref{inv1}, that will be easier to work with in practice.

Let $\{X_k\}_{k\ge 0}$ and $\{\xp_k\}_{k\ge 0}$ be two chains from the kernel defined by $(X,\xp)$, with arbitrary initial values, and coupled according to some coupling scheme which satisfies the following property:
\vskip.1in
\begin{itemize}
\item[\bf P] For every initial value  $(x,y)$, and every $k$, the marginal distribution of $X_k$ depends only on $x$ and the marginal distribution of $\xp_k$ depends only on $y$.
\end{itemize}
\vskip.1in
Under this assumption, we have the following lemma:
\begin{lmm}\label{invlmm}
Suppose the chains $\{X_k\}$ and $\{\xp_k\}$ satisfy the property \textup{\bf P} described above. 
Let $f:\xx \ra \rr$ be a function such that $\ee f(X) = 0$. Suppose there exists a finite constant $L$ such that for every $(x,y)\in \xx^2$, 
\begin{equation}\label{summable}
\sum_{k=0}^\infty |\ee(f(X_k)-f(\xp_k)|X_0=x,\xp_0=y)| \le L.
\end{equation}
Then, the function $F$ defined as
\[
F(x,y) := \sum_{k=0}^\infty \ee(f(X_k)-f(\xp_k)|X_0 = x, \xp_0 = y)
\]
satisfies $F(X,\xp) = -F(\xp,X)$ and $\ee(F(X,\xp)|X) = f(X)$.
\end{lmm} 
{\bf Remark.} In practice, we will start our chains with $X_0=X$ and $\xp_0=\xp$ for directly obtaining bounds on $F(X,\xp)$ in the process of verifying (\ref{summable}).
\vskip.1in

\noindent{\bf Proof.} Property {\bf P} implies that
\[
\ee(f(X_k)|X_0=x,\xp_0=y) = \ee(f(X_k)|X_0=x) = P^k f(x),
\]
where $P$ is the kernel defined in (\ref{kernel}). Similarly, $\ee(f(\xp_k)|X_0=x,\xp_0=y) = P^k f(y)$. The rest is a rewriting of the conclusions of Lemma \ref{inv1}. \hfill $\Box$
\vskip.2in

\noindent{\bf Example.} Let $\xx = \rr^n$, and let $X= (X_{(1)},\ldots,X_{(n)})$ be a vector with independent components. (The subscripts are put in brackets because we will be using $\{X_k\}_{k\ge 0}$ to denote a Markov chain on $\rr^n$.) As a preliminary example, and also as a prelude to section \ref{dobrushin}, we shall now work out the concentration of $f(X)$, where $f$ is a function satisfying
\[
|f(x)-f(y)|\le \sum_{i=1}^n c_i\ii\{x_i\ne y_i\}.
\] 
In other words, if $x$ and $y$ differ only at coordinate $i$, then $|f(x)-f(y)|\le c_i$. Our technique will give an easy way to recover a version of the bounded difference inequality (Theorem \ref{bddiff}) without using martingale or information theoretic results. We shall only give a brief description of the steps in our solution, because the details will be worked out under a more general setting in the next section.

We produce $\xp$ by choosing a coordinate $I$ uniformly at random, and replacing the $I^{\mathrm{th}}$ coordinate of $X$ by $X^*_{(I)}$, where the vector $(X_{(1)}^*,\ldots,X_{(n)}^*)$ has the same distribution as $X$ and is independent of $X$. It is easy to see that $(X,\xp)$ is an exchangeable pair. 

The coupling is done in the natural way: At every step, choose the same coordinate $I$ for both chains, and replace the $I^{\mathrm{th}}$ coordinate of both by the same realization of $X_{(I)}^*$. Since the number of coordinates at which $X_k$ and $\xp_k$ differ cannot increase with $k$, it is clear (by the coupon collector's problem) that the coupling time for the chains has expected value bounded by $n\log n$. This proves condition (\ref{summable}). 

Now suppose we start with $(X_0,\xp_0)=(x,y)$. If $x$ and $y$ differ only at coordinate $i$, then $|f(X_k)-f(\xp_k)| \le c_i$ for every $k$, and the expected coupling time is $n$. Thus, by the representation of $F$ in Lemma \ref{invlmm}, we get
$|F(x,y)|\le n c_i$. Consequently, we have
\begin{align*}
v(X) &= \frac{1}{2}\ee\bigl(|(f(X)-f(\xp))F(X,\xp)|\bigl|X\bigr) \\
&\le \frac{1}{2n}\sum_{i=1}^n (c_i)(nc_i) = \frac{1}{2}\sum_{i=1}^n c_i^2.
\end{align*}
Theorem \ref{hoeffding} now gives us a version of the well-known bounded difference inequality (Theorem \ref{bddiff}), albeit with a missing $2$ in the exponent.

\section{Concentration under weak dependence}\label{dobrushin}
Let $\xx = \Omega^n$, where $\Omega$ is a Polish space. Consequently, our random variable $X$ will now have $n$ coordinates, which we shall assume to be {\it weakly dependent} in the sense of Dobrushin, a familiar notion from statistical mechanics which we shall define later in this section.

We shall also assume that our function $f$ satisfies a Lipschitz condition with respect to a generalized Hamming distance on $\xx$:
\begin{equation}\label{lipsch}
|f(x)-f(y)| \le \sum_{i=1}^n c_i \ii\{x_i \ne y_i\}
\end{equation}
for some fixed constants $c_1,\ldots,c_n$. As mentioned before, this just means that the value of the function does not change by more than $c_i$ if the $i^{\mathrm{th}}$ coordinate is altered. 

The generalization of the Hamming metric is useful because it allows us to get concentration for lower dimensional marginals. For instance, if we want the concentration of a function of $(X_1, \ldots, X_k)$, where $k\le n$, we can just consider every function of $(X_1,\ldots,X_k)$ as a function of $(X_1,\ldots,X_n)$ with $c_{k+1}=\cdots = c_n = 0$.

In the following, we shall be using the notation $\bar{x}^i$ to  denote the element of $\Omega^{n-1}$ obtained by omitting the $i^{\mathrm{th}}$ coordinate of the vector $x\in \Omega^n$. 

We shall let $\mu$ denote the law of $X$, and for each $i$ and $x$, $\mu_i(\cdot|\bar{x}^i)$ will stand for the law of the $i^{\mathrm{th}}$ coordinate of $X$ given that $\bar{X}^i = \bar{x}^i$. 

At this point, let us also mention that we shall  usually denote the $i^{\mathrm{th}}$ coordinate of a vector $x$ by $x_i$. Unfortunately, we are also using the notation $X_k$ for the $k^{\mathrm{th}}$ element of a Markov chain. We hope that the context will clarify any ambiguity. In particular, to denote the $i^{\mathrm{th}}$ coordinate of $X_k$ we shall use the notation $X_{k,i}$. 

Finally, let us recall that for a square matrix $A$, the $L^2$ operator norm of $A$ is defined as:
\[
\|A\|_2 := \max_{\|y\|=1} \|Ay\|.
\]
While we are at it, let us also recall that the $L^\infty$ and $L^1$ operator norms can be expressed as 
\[
\|A\|_\infty := \max_{1\le i\le n} \sum_{j=1}^n |a_{ij}|, \ \text{ and } \ \|A\|_1 := \max_{1\le  j\le n} \sum_{i=1}^n |a_{ij}|.
\]
It is a frequently useful fact that $\|A\|_2 \le \sqrt{\|A\|_1 \|A\|_\infty}$.

The following theorem gives a simple way of obtaining concentration bounds for $f(X)$ under a contractivity condition on the conditional laws:
\begin{thm}\label{dobthm}
Suppose $A = (a_{ij})$ is an $n\times n$ matrix with nonnegative entries and zeros on the diagonal such that for any $i$, and any $x, y \in \Omega^n$, 
\[
d_{TV}(\mu_i(\cdot|\bar{x}^i), \mu_i(\cdot|\bar{y}^i)) \le \sum_{j=1}^n a_{ij} \ii\{x_j\ne y_j\},
\]
where $d_{TV}$ is the total variation distance on the space of probability measures on $\Omega$.
Suppose $f$ satisfies the generalized Lipschitz condition (\ref{lipsch}). If $\|A\|_2 < 1$, we have
\[
\pp\{|f(X)-\ee f(X)| \ge t\} \le 2\exp\biggl(-\frac{(1-\|A\|_2)t^2}{\sum_i c_i^2}\biggr)
\]
for each $t\ge 0$.
\end{thm}
{\bf Remarks.} The matrix $A$ is called ``Dobrushin's interdependence matrix'' and the condition $\|A\|_\infty < 1$ is called Dobrushin's condition (as introduced by Dobrushin \cite{dobrushin70} and extended by Dobrushin \& Shlosman \cite{dobshlos85, dobshlos87}). Dobrushin's condition ensures, among other things, the uniqueness of Gibbs states at high temperature. In a series of important papers, Stroock and Zegarlinski \cite{stroockzeg92a, stroockzeg92b, stroockzeg92c} showed that for spin systems on a lattice, the Dobrushin-Shlosman conditions are equivalent to the validity of logarithmic Sobolev inequalities for the associated Glauber dynamics. Though concentration inequalities follow from log-Sobolev inequalities, explicit constants are not available from this body of work. Our approach gives a direct way of getting explicit concentration bounds from Dobrushin type conditions. For more on Dobrushin's condition and its consequences, see Georgii's book \cite{georgii88}, Chapter 8. For a recent treatise on logarithmic Soboloev inequalities in the context of spin systems and Dobrushin-Shlosman type conditions, see Ledoux \cite{ledoux01a}. 

As mentioned in section \ref{entropy}, Dobrushin type conditions were recently used by Marton \cite{marton03, marton04} to obtain transportation cost inequalities (and hence concentration) under weak dependence. However, Marton's results do not seem to be suited for the setting that we are now working under. Besides, they require us to know some bounds on certain log-Sobolev constants, which are hard to obtain.
\\
\\
{\bf Proof of Theorem \ref{dobthm}.} To prove this theorem, we shall construct a reversible Markov kernel and a suitable coupling, and then directly apply Lemma \ref{invlmm} in conjunction with the tools from Chapter \ref{results}. 

First, define a reversible Markov kernel as follows: At each step, choose a coordinate $I$ uniformly from $\{1,\ldots,n\}$. Then replace the $I^{\mathrm{th}}$ coordinate of the current state of the chain by an element of $\Omega$ chosen according to the conditional distribution of the $I^{\mathrm{th}}$ coordinate given the values of the other coordinates. This is the usual Gibbs sampling technique, and it is well-known and easy to prove that the chain is reversible. This is also known as the ``Glauber dynamics'' in the case of spin systems. 

Now, we describe the coupling. Suppose at any stage, the $X$-chain is at $x$, and the $\xp$-chain is at $y$. Choose a coordinate $I$ uniformly at random. 
By the well-known property of the total variation distance, we can have two $\xx$-valued random variables $W_1^I$ and $W_2^I$ define on some probability space such that $W_1^I \sim \mu_I(\cdot|\bar{x}^I)$, $W_2^I \sim \mu_I(\cdot|\bar{y}^I)$, and
\begin{align*}
\pp\{W_1^I \ne W_2^I\} &= d_{TV}(\mu_I(\cdot|\bar{x}^I), \mu_I(\cdot|\bar{y}^I))\\
&\le \sum_{j=1}^n a_{Ij} \ii\{\bar{x}^I_j \ne \bar{y}^I_j\} \ \ \ \text{ (by assumption).}
\end{align*}
Having obtained $W_1^I$ and $W_2^I$, update the $X$-chain by putting $W_1^I$ as the $I^{\mathrm{th}}$ coordinate (keeping all other coordinates the same) and update the $\xp$-chain the same way using $W_2^I$. This is perhaps the most naive way to couple dependent variables, and is commonly known as ``the greedy coupling''. 

Of course, people know how to make the above construction completely mathematically precise, and so we prefer to gloss over that aspect. But if someone is interested, he can look up the recent paper of Marton \cite{marton04}, where a more complex coupling has been handled in a very mathematically precise manner.

It is clear from the definition that this coupling indeed gives the Gibbs sampler chains as the marginals. To see that property {\bf P} holds, note that for any $i$, the distribution of $W_1^i$ as defined above depends only on the current state of the $X$-chain, irrespective of the status of the $\xp$-chain. Similarly, the distribution of $W_2^i$ depends only on the current state of the $\xp$-chain. Thus, the distribution of $X_{k+1}$ given $(X_k,\xp_k)$ depends only on $X_k$. Since the coupling is Markovian, we can proceed by induction to get a proof of property {\bf P} for this coupling. 

Now assume without loss of generality that $\ee f(X) = 0$. The summability condition (\ref{summable}) will be verified later in the proof. 

Fix $k\ge 0$, and suppose we used the above scheme to move to the $(k+1)^{\mathrm{th}}$ stage. Then for any $i$,
\begin{align*}
\pp\{X_{k+1,i}\ne \xp_{k+1,i} \text{ and } I\ne i|X_k,\xp_k\} &= \biggl(1-\frac{1}{n}\biggr) \ii\{X_{k,i}\ne \xp_{k,i}\},
\end{align*}
and 
\begin{align*}
& \pp\{X_{k+1,i}\ne \xp_{k+1,i} \text{ and } I = i|X_k,\xp_k\} \\
&= \frac{1}{n}\pp\{W_1^i\ne W_2^i\} \le \frac{1}{n}\sum_{j=1}^n a_{ij} \ii\{X_{k,j} \ne \xp_{k,j}\}.
\end{align*}
Adding up, and taking expectation given $(X_0, \xp_0)$ on both sides, we get
\begin{align*}
&\pp\{X_{k+1,i} \ne \xp_{k+1,i}|X_0,\xp_0\} \\
&\le \biggl(1-\frac{1}{n}\biggr) \pp\{X_{k,i} \ne \xp_{k,i}|X_0,\xp_0\} +\frac{1}{n}\sum_{j=1}^n a_{ij} \pp\{X_{k,j}\ne \xp_{k,j}|X_0,\xp_0\}.
\end{align*}
For each $k$, let $\ell_k = \ell_k(X_0,\xp_0)$ be the random vector whose $i^{\mathrm{th}}$ component is $\pp\{X_{k,i} \ne \xp_{k,i}|X_0,\xp_0\}$. Let $B = (1-\frac{1}{n})I + \frac{1}{n}A$, where $I$ denotes the identity matrix. 
The above inequality shows that $\ell_{k+1} \le B\ell_k$ for every $k$, where `$x \le y$' means that $x_i \le y_i$ for each coordinate $i$. Since everything is nonnegative, we can continue by induction to see that
\[
\ell_k \le B^k \ell_0.
\]
Now let $c = (c_1,\ldots, c_n)$ be the vector of constants from (\ref{lipsch}). Note that
\begin{equation}\label{b2}
\|B\|_2 \le 1-\frac{1}{n} + \frac{\|A\|_2}{n} < 1,
\end{equation}
and hence, by the generalized Lipschitz property of $f$, 
\[
\sum_{k=0}^\infty |\ee(f(X_k)-f(\xp_k)|X_0,\xp_0)| \le \sum_{k=0}^\infty \ee(c\cdot\ell_k|X_0,\xp_0) \le \|c\|\|\ell_0\|\sum_{k=0}^\infty \|B\|^k  < \infty.
\]
This establishes condition (\ref{summable}), and so we can now invoke Lemma \ref{invlmm}. Produce $\xp$ by starting from $X$ and taking one step according to the Gibbs sampler kernel. Putting $X_0=X$ and $\xp_0=\xp$, we get
\[
F(X,\xp) = \sum_{k=0}^\infty\ee(f(X_k)-f(\xp_k)|X_0,\xp_0)
\]
by Lemma \ref{invlmm}. Thus, if $v$ is defined as in (\ref{ddef}), then 
\begin{equation}\label{vbd}
v(X)\le \frac{1}{2}\sum_{k=0}^\infty \ee\bigl(|(f(X_0)-f(\xp_0))(f(X_k)-f(\xp_k))|\bigl|X_0).
\end{equation}
Now, by our definition of $\ell_k$ and the property (\ref{lipsch}) of $f$, we have
\[
\ee\bigl(|(f(X_0)-f(\xp_0))(f(X_k)-f(\xp_k))|\bigl|X_0,\xp_0) \le (c\cdot \ell_0)(c\cdot \ell_k).
\]
Note that with $(X_0,\xp_0)=(X,\xp)$, $X_0$ and $\xp_0$ can differ at most at one coordinate. Suppose they differ at coordinate $i$. Then
\[
(c\cdot \ell_0)(c\cdot \ell_k) \le (c\cdot \ell_0)(c\cdot B^k\ell_0) = c_i(c \cdot B^k e_i),
\]
where $e_i$ denotes the vector whose $i^{\mathrm{th}}$ coordinate is $1$ and all other coordinates are $0$. Now, the probability that the $X_0$ and $\xp_0$ differ at the $i^{\mathrm{th}}$ coordinate is $\le 1/n$. Hence,
\begin{align*}
&\ee\bigl(|(f(X_0)-f(\xp_0))(f(X_k)-f(\xp_k))|\bigl|X_0) \\
&\le \frac{1}{n}\sum_{i=1}^n c_i(c\cdot B^ke_i) = \frac{1}{n} c\cdot B^k c \le \frac{\|B\|^k_2\|c\|^2}{n}
\end{align*}
Thus, by (\ref{vbd}) and (\ref{b2}), we get
\[
v(X) \le \frac{\|c\|^2}{2n(1-\|B\|_2)} \le \frac{\|c\|^2}{2(1-\|A\|_2)}.
\]
The proof is now completed by using Theorem \ref{hoeffding}. \hfill $\Box$

\section{Example: Spin systems}\label{spinsystems}
We now specialize to the case $\Omega = [-1,1]$. Let $\nu$ be a probability measure on $[-1,1]$. Suppose our original measure $\mu$ on $[-1,1]^n$ has a density with respect to $\nu^n$, represented in the Boltzmann form: $Z^{-1} e^{-H(x)}$, where $Z$ is the normalizing constant.

The measure $\nu$ will usually be either the normalized Lebesgue measure, or the symmetric distribution on $\{-1,1\}$. The generalized framework will allow us to deal with both discrete and continuous problems.

The following lemma gives us a simple way of applying Theorem \ref{dobthm} in this setting under the assumption that $H$ can be extended to a twice continuously differentiable function on $[-1,1]^n$. 
\begin{lmm}\label{spinlmm}
For each pair $(i,j)$ with $i\ne j$, define
\[
a_{ij} := 4\sup_{x\in [-1,1]^n} \biggl|\frac{\partial^2 H}{\partial x_i \partial x_j}(x)\biggr|
\]
and let $a_{ii} = 0$ for each $i$. Then for each $i$ and $x,y\in [-1,1]^n$, we have
\[
d_{TV}(\mu_i(\cdot|\bar{x}^i), \mu_i(\cdot|\bar{y}^i)) \le \sum_{j=1}^n a_{ij} \ii\{x_j\ne y_j\}.
\]
\end{lmm}
{\bf Remark.} It is obvious that this lemma, combined with Theorem \ref{dobthm}, gives concentration inequalities for most of the familiar spin models at sufficiently high temperature. It covers, for instance, the Ising ferromagnets (in any dimension), the xy-model, and even the Ising spin glasses. This lemma is in fact, just a more easily recognizable version (for probabilists) of Proposition 8.8 in Georgii's book \cite{georgii88}, which gives a similar condition for Gibbs potentials to satisfy Dobrushin's condition, and gives some more examples.
\vskip.1in

\noindent{\bf Example.} Consider the Ising model on a graph $G = (V,E)$, where $V = \{1,\ldots,n\}$ and the maximum degree is $r$. Here $\Omega = \{-1,1\}$ and $\nu$  is the symmetric probability distribution on $\Omega$. The Hamiltonian at inverse temperature $\beta$ and external field $h$ is given by 
\[
H_\beta(x) := -\beta \sum_{(i,j)\in E} x_ix_j - \beta h \sum_i x_i.
\]
We shall consider the concentration of the {\it magnetization}, defined as
\[
m(x) := \frac{1}{n}\sum_{i=1}^n x_i.
\]
The Hamiltionian has a natural extension as a $C^2$ function on $\rr^n$, and a simple computation gives
\[
a_{ij} =
\begin{cases}
4\beta &\text{ if } (i,j)\in E,\\
0 &\text{ if } (i,j)\not \in E.
\end{cases}
\]
Thus, the $L^\infty$ and $L^1$ norms of the matrix $A = (a_{ij})$ are both bounded by $4\beta r$, where recall that $r$ is the maximum degree of $G$. Thus, $\|A\|_2 \le \sqrt{\|A\|_1\|A\|_\infty} \le 4\beta r$.

Again, it is clear that
\[
|m(x)-m(y)| \le \frac{2}{n}\sum_{i=1}^n\ii\{x_i\ne y_i\}
\]
Thus, if $X = (X_1,\ldots,X_n)$ is a spin configuration drawn from the Gibbs measure at inverse temperature $\beta < 1/4r$ and external field $h$, and $m(X)$ is its magnetization, then by Lemma \ref{spinlmm} and Theorem \ref{dobthm}, we have
\[
\pp\{|m(X) -\ee m(X)| \ge t\} \le e^{-\frac{1}{4}n(1-4\beta r) t^2}.
\]
This gives an explicit concentration bound on the magnetization at sufficiently high temperature, though it need not necessarily cover the entire high temperature domain.
\vskip.3in

\noindent {\bf Proof of Lemma \ref{spinlmm}.} The density with respect to $\nu$ of the conditional distribution $\mu_i(\cdot|\bar{x}^i)$ is given by
\[
\rho_i(u|\bar{x}^i) = \frac{e^{-H(u, \bar{x}^i)}}{\int e^{-H(v,\bar{x}^i)} \nu(dv)},
\] 
where $(u, \bar{x}^i)$ denotes the vector obtained by substituting the number $u$ as the $i^{\mathrm{th}}$ coordinate of the vector $x$. 

Now fix $i$ and $j$, where $i\ne j$. Direct computation shows that
\begin{align*}
\frac{\partial}{\partial x_j} \rho_i(u|\bar{x}^i) &=  -\biggl(\frac{\partial H}{\partial x_j}(u, \bar{x}^i) - \int\frac{\partial H}{\partial x_j}(v, \bar{x}^i) \rho_i(v|\bar{x}^i)\nu(dv)\biggr) \rho_i(u|\bar{x}^i)
\end{align*}
and hence
\begin{align*}
\biggl|\frac{\partial}{\partial x_j} \rho_i(u|\bar{x}^i)\biggr| &\le \sup_v \biggl| \frac{\partial H}{\partial x_j}(u, \bar{x}^i) - \frac{\partial H}{\partial x_j}(v, \bar{x}^i)\biggr| \rho_i(u|\bar{x}^i) \\
&\le \frac{a_{ij}}{2} \rho_i(u|\bar{x}^i).
\end{align*}
Thus, for any $A\subseteq [-1,1]$, we have
\begin{align*}
\biggl|\frac{\partial}{\partial x_j} \mu_i(A|\bar{x}^i)\biggr| &= \biggl|\int_A \frac{\partial}{\partial x_j} \rho_i(u|\bar{x}^i) \nu(du)\biggr| \le \frac{a_{ij}}{2}
\end{align*}
Thus, we can conclude that for any $x,y \in [-1,1]^n$, 
\begin{align*}
&d_{TV}(\mu_i(\cdot|\bar{x}^i), \mu_i(\cdot|\bar{y}^i)) \\
&= \sup_A |\mu_i(A|\bar{x}^i) - \mu_i(A|\bar{y}^i)| \le \sum_{j=1}^n a_{ij} \ii\{x_j\ne y_j\}.
\end{align*}
This completes the argument. \hfill $\Box$

\section{Example: Graph colorings}\label{graphcolorings}
Suppose $G= (V,E)$ is a graph with $n$ vertices and maximum degree $r$. A $k$-coloring of $G$ is an assignment of $k$ colors to the vertices of $G$. In other words, it is an element of $\{1,\ldots, k\}^G$. A {\it proper} $k$-coloring is a $k$-coloring of $G$ in which no two adjacent vertices receive the same color. 

Let $X = (X_i, i\in V)$ be a coloring of $G$ chosen uniformly from the set of all proper $k$-colorings. 

A substantial amount of recent effort in theoretical computer science and associated probability theory has been devoted to the study of this random variable $X$ and Markov chains converging in distribution to $X$. Most of it centers around temporal and spatial mixing and decay of correlations. An early observation of Jerrum \cite{jerrum95} and Salas \& Sokal \cite{salassokal97}, followed by vigorous activity in the last few years, have resulted in a spate of sophisticated coupling techniques and improved results. For up-to-date references, one can look at the recent articles \cite{coloring04, coloring04b}.

However, although coupling techniques have revolutionized the analysis of spatial and temporal decay of correlations in graph colorings, it is still extremely difficult to get concentration inequalities for general functions in these settings. In fact, we can see no direct way of getting concentration bounds for as simple a functional as the proportion of vertices having a particular color using currently available techniques.
 
In this section, we shall use our methods to get concentration inequalities for arbitrary (generalized) Lipschitz functions of randomly chosen proper graph colorings. That is, we shall consider $f:\{1,\ldots,k\}^G \ra \rr$ satisfying $|f(x)-f(y)|\le \sum_{i=1}^n c_i\ii\{x_i\ne y_i\}$, and get tail bounds on $f(X)$, where $X$ is our random proper coloring of $G$. 

Both Jerrum \cite{jerrum95} and Salas \& Sokal \cite{salassokal97} use the greedy coupling to get decay of correlations under the condition $k > 2r$. This is sort of a ``high temperature phase'' for proper graph colorings. Since our purpose is only to demonstrate the manner of application of our technique, we shall stick to this primitive approach. 
\begin{prop}\label{coloring}
Suppose $G$ is a finite graph with maximum degree $r$, $X$ is a uniformly chosen proper $k$-coloring of $G$, and $f:\{1,\ldots,k\}^G\ra \rr$ satisfies $|f(x)-f(y)|\le \sum_{i=1}^n c_i\ii\{x_i\ne y_i\}$. If $k > 2r$, then for any $t\ge 0$ we have
\[
\pp\{|f(X)-\ee f(X)| \ge t\} \le 2\exp\biggl(-\frac{\gamma t^2}{\sum_{i=1}^n c_i^2}\biggr),
\] 
where $\gamma = (k-2r)/(k-r)$.
\end{prop}
{\bf Example 1.} Let $f(x)$ be the proportion of vertices receiving color $1$. It's clear that  
\[
|f(x)-f(y)| \le \frac{1}{n}\sum_{i=1}^n\ii\{x_i\ne y_i\}.
\]
Also, by symmetry, $\ee f(X) = 1/k$, where $X$ is drawn uniformly from the set of all proper $k$-colorings. Using the above result, we get
\[
\pp\{|f(X)- k^{-1}|\ge t\} \le 2e^{-n\gamma t^2},
\]
where $\gamma = (k-2r)/(k-r)$ as in the lemma.
\vskip.1in
\noindent {\bf Example 2.} Let $f(x)$ be the number of neighbors of vertex $1$ which receive color $1$. Then we can take $c_i = 1$ if $i$ is a neighbor of $1$ and $c_i=0$ otherwise. Consequently, $\sum_i c_i^2 =  r_1$, where $r_1$ is the degree of vertex $1$. Thus, 
\[
\pp\{|f(X)-\ee f(X)| \ge t\} \le 2e^{-\gamma t^2/r_1}.
\]
\vskip.3in

\noindent {\bf Proof of Proposition \ref{coloring}.} As usual, denote the law of $X_i$ given $\bar{X}^i$ by $\mu_i(\cdot|\bar{x}^i)$. 
It is easy to see from definition that $\mu_i(\cdot|\bar{x}^i)$ is the uniform distribution over the set
\[
\{1,\ldots,k\} \backslash \{x_j: (i,j)\in E\}.
\]
Thus, given two vectors $x$ and $y$ which differ only at a coordinate $j$, where $j$ is adjacent to $i$, it is clear that 
\begin{align*}
d_{TV}(\mu_i(\cdot|\bar{x}^i), \mu_i(\cdot|\bar{y}^i)) &= \frac{1}{2}\sum_{u=1}^k |\pp\{X_i=u|\bar{X}^i = \bar{x}^i\} - \pp\{X_i = u|\bar{X}^i = \bar{y}^i\}|\\
&\le \frac{1}{2}\bigl(\pp\{X_i=y_j|\bar{X}^i=\bar{x}^i\} + \pp\{X_i = x_j |\bar{X}^i = \bar{y}^i\}\bigr)\\
&\le \frac{1}{k-r}.
\end{align*}
Thus, we can take $a_{ij} = 1/(k-r)$ if $(i,j)\in E$ and $a_{ij}=0$ otherwise, and apply Theorem \ref{dobthm} to complete the proof. \hfill $\Box$

\section{Concentration of Haar measures}\label{group}
Let $G$ be a compact topological group. Then there exists a $G$-valued random variable $X$ with the properties that for any $x\in G$, the random variables $xX$, $Xx$ and $X^{-1}$ all have the same distribution as $X$. As mentioned in section \ref{gpreview}, where we review the literature about concentration of measure on groups, the existence and uniqueness of the distribution of $X$ (which is called the ``Haar measure on $G$''), is a classical result (see e.g.\ Rudin \cite{rudin73}, Theorem 5.14). Our goal in this section is to study the concentration of the Haar measure using the convergence properties of certain kinds of Markov chains which have the Haar measure as their stationary distribution. An original application to the concentration of the spectra of sums of random matrices (related to free probability) will be given in the next section.

We shall not repeat the discussion of the existing literature on concentration of Haar measures, which is done in section \ref{gpreview} of the review chapter. It suffices to say that there is not much work. The only general theorem we could detect is a result of Milman \& Schechtman \cite{milmanschechtman86} which is stated in section \ref{gpreview}, where we also discuss available results for special groups like $SO_n$ (the special orthogonal group) and $S_n$ (the symmetric group). 

We should, however, mention that we are not the first to analyze measures on groups using Stein's method. Diaconis \cite{diaconis91} has an application of Stein's method to the analysis of random walks on groups. More recently, Jason Fulman has applied Stein's method to analyze the Haar measure \cite{fulman04} and the Plancherel measure \cite{fulman05} on~$S_n$.
\vskip.1in

\noindent Let us now introduce our setting. Let $Y$ be a $G$-valued random variable having the following properties:
\begin{enumerate}
\item The random variable $Y^{-1}$ has the same distribution as $Y$; that is, the law of $Y$ is {\it symmetric}.
\item For any $x\in G$, $xYx^{-1}$ has the same distribution as $Y$. In other words, the distribution of $Y$ is uniform on each conjugacy class of $G$.
\end{enumerate} 
Any such $Y$ defines a reversible Markov kernel $P$ is a natural way: For any $f:G\ra \rr$ such that $\ee|f(X)|<\infty$, let
\begin{equation}\label{ykernel}
Pf(x) := \ee f(Yx) = \ee f(xx^{-1}Yx) = \ee f(xY).
\end{equation}
The reversibility of this kernel can be proved as follows: Since $yX$ has the same distribution as $X$ for any $y\in G$, therefore $Y$ and $YX$ are independent. Also, $Y^{-1}$ has the same distribution as $Y$. Hence, the pair $(X,Y)$ has the same distribution as $(YX, Y^{-1})$. Consequently, the pairs $(X,YX)$ and $(YX, Y^{-1}YX) = (YX,X)$ also have the same distribution. In other words, if we let $\xp=YX$, the $(X,\xp)$ is an exchangeable pair. Finally, note that $Pf(x)=\ee(f(\xp)|X=x)$. 

We shall attempt to study the concentration of $X$ using the properties of this kernel $P$. A version of this problem was considered by Schmuckenschl\"ager \cite{schmuck98}, but no practical solution was given. We shall now present a theorem that should suffice in many problems, but with ``an extra $\log n$ factor''. As mentioned before, concrete examples will be provided later on.
\begin{thm}\label{gpaction}
Let $G, X, Y$ be as above. Let $f:G \ra \rr$ be a measurable function such that $\ee f(X) = 0$. Let $\|f\|_\infty = \sup_{x \in G} |f(x)|$ and 
\[
\|f\|_Y := \sup_{x\in G} \bigl[\ee(f(x)-f(Yx))^2\bigr]^{1/2}.
\]
Let $Y_1,Y_2,\ldots,$ be i.i.d.\ copies of $Y$. Suppose $a$ and $b$ are two positive constants such that $d_{TV}(Y_1Y_2\cdots Y_k, X) \le ae^{-bk}$ for every $k$, where $d_{TV}$ is the total variation metric. Let $A$ and $B$ be two numbers such that $\|f\|_\infty \le A$ and $\|f\|_Y \le B$. Let
\[
C = \frac{B^2}{2b}\biggl[\biggl(\log \frac{4aA}{B}\biggr)^+ + \frac{b(2-e^{-b})}{1-e^{-b}}\biggr].
\]
Then $\var(f(X))\le C$, and for any $t\ge 0$, $\pp\{|f(X)|\ge t\} \le 2e^{-t^2/2C}$. 
\end{thm}
{\bf Remarks.} The rate of decay of the total variation distance for random walks on groups has been a topic of interest ever since the work of Diaconis \& Shahshahani~\cite{diaconis81}, who introduced group representation tools to overcome the deficiencies of the usual spectral gap approach (which gives a wrong answer for the random walk with transpositions on $S_n$). The property that $Y$ is uniformly distributed on conjugacy classes is particularly suited to this line of analysis. 

Since the seminal paper \cite{diaconis81}, much work has been done; several other methods have been developed by various authors, and the rates have been evaluated in many interesting problems. 
For a recent survey with an extensive bibliography, we refer to Diaconis \cite{diaconis03}. 

To the best of our knowledge, Theorem \ref{gpaction} is the only result as of now, which connects the rates of convergence of random walks on groups --- which is a widely explored area --- with the concentration of Haar measures, which is not so widely explored.

Before we move on to complete the proof of Theorem \ref{gpaction}, it is time for an easy application, which will not give a better-than-existing result, but merely demonstrate how Theorem \ref{gpaction} can be applied. An original application will be given in the next section. 
\vskip.1in
\noindent {\bf Random permutations.}
Let $G=S_n$, the group of all permutations of $n$ elements. Then the Haar measure is just the uniform distribution on $G$. Let $X$ be a uniformly distributed random variables on $G$. We define the distribution of our $Y$ by putting mass $1/n$ on the identity permutation and $2/n^2$ on each transposition of two elements. Since the set of transpositions is closed under conjugation and inversion, $Y$ satisfies the required properties 1 and 2 from the previous section. The kernel $P$ defined by $Y$ as in (\ref{ykernel}) is a familiar object which has been studied by several authors. The study was initiated by Diaconis \& Shahshahani \cite{diaconis81}, who proved the following result:
\begin{thm}\label{diaconis}
\textup{[Diaconis \& Shahshahani \cite{diaconis81}]} Let $Y_1,Y_2,\ldots,$ be i.i.d.\ copies of the random variable $Y$ defined above. Then
$d_{TV}(Y_1\cdots Y_k, X) \le 6n e^{-2k/n}$
for every $k > \frac{1}{2}n\log n$. 
\end{thm}
The above version of the result has been quoted from \cite{diaconis03}, where it is written in a slightly different manner. Note that if we substitute $k=\frac{1}{2}n\log n$, the right hand side becomes $6$, which is $\ge 1$, and hence the condition $k>\frac{1}{2}n\log n$ is redundant. 

Using the Diaconis-Shahshahani result, we get $a = 6n$ and $b= 2/n$ in Theorem \ref{gpaction}, and hence the following result:
\begin{prop}\label{permgen}
Let $G=S_n$, where $n\ge 2$, and let $X,Y$ be as above. Let $f:G \ra \rr$ be a function such that $\ee f(X)=0$. As in Theorem  \ref{gpaction}, let $\|f\|_Y = \max_{\sigma\in G} [\ee(f(\sigma)-f(Y\sigma))^2]^{1/2}$. Let $A$ and $B$ be two numbers such that $\|f\|_\infty \le A$ and $\|f\|_Y \le B$. Let
\[
C = \frac{nB^2}{4}\biggl[\biggl(\log \frac{24nA}{B}\biggr)^+ + \frac{(2/n)(2-e^{-2/n})}{1-e^{-2/n}} \biggr].
\]
Then $\var(f(X))\le C$, and for any $t\ge 0$, we have $\pp\{|f(X)|\ge t\} \le 2e^{-t^2/2C}$.
\end{prop}
For a simple application of this proposition, consider the {\it descent function} on $S_n$, defined as
\[
D(\sigma) = \sum_{i=1}^{n-1} \ii\{\sigma(i) > \sigma(i+1)\}.
\]
The number of descents of a permutation is an interesting quantity from statistical and combinatorial points of view. They have been studied extensively by Foata and Schutzenberger \cite{foata70}, Knuth \cite{knuth73} and in the unpublished notes of Diaconis \& Pitman~\cite{diaconispitman91}. Fulman \cite{fulman04} has applied Stein's method to the study of descents.

Let $X$ be uniformly distributed on $S_n$. Clearly, $\|D(X)-\ee D(X)\|_\infty \le n$. Now, for any $x\in G$ and any transposition $y$, $|D(x)-D(yx)|\le 4$. Putting $A=n$ and $B = 4$ in Theorem \ref{permgen}, and assuming $n\ge 10$, we get $C\le 4n(2\log n + 3.1)$. Thus, for $n \ge 10$, we get $\var(D(X)) \le 4n(2\log n + 3.1)$ and for any $t\ge 0$, 
\[
\pp\{|D(X) -\ee D(X)|\ge t\} \le 2\exp\biggl(-\frac{t^2}{8n(2\log n + 3.1)}\biggr).
\]
The bounds are clearly off by a factor of $\log n$, but that will be a perpetual inconvenience of using Theorem \ref{gpaction}. However, it is interesting to note that Maurey's theorem (\cite{maurey79}, stated as Theorem \ref{maurey} in section \ref{gpreview}) gives the bound $\pp\{|D(X)-\ee D(X)|\ge t\} \le 2e^{-t^2/256n}$ in this problem, which, though technically ``better'' than our bound, is always going to be worse in all practical situations.

That apart, there is a more serious reason why Proposition \ref{permgen} may give better results in  some problems. Maurey's result essentially uses $\|f\|_{\mathrm{Lip}}$ instead of $\|f\|_Y$, where 
$\|f\|_{\mathrm{Lip}} = \max\{ |f(\sigma)-f(\tau\sigma)| : \sigma\in S_n, \ \tau \text{ is a transposition}\}$.
Clearly, $\|f\|_Y \le \|f\|_{\mathrm{Lip}}$, and the difference may be significant in some situations.

Talagrand's result (\cite{talagrand95}, Theorem 5.1), on the other hand, will probably dominate Proposition \ref{permgen} most of the time, though it is not clear whether Talagrand's method can be used to derive a result based on something like $\|f\|_Y$.
 
We now move on to prove the central result of this section.
\vskip.1in
\noindent {\bf Proof of Theorem \ref{gpaction}.} Let $\xp = YX$. As observed before, $(X,\xp)$ is an exchangeable pair. Recall that we defined $Pf(x)= \ee f(Yx)$. By Lemma \ref{inv1}, 
\begin{align}\label{vx}
v(x) &\le \frac{1}{2}\sum_{k=0}^\infty \ee|(f(x)-f(Yx)(P^kf(x)-P^kf(Yx))|,
\end{align}
where $v(x)$ is the usual quantity in our theory, as defined in (\ref{ddef}). The criterion (\ref{sum1}) required for Lemma \ref{inv1} holds, since for any $z\in G$,
\begin{align*}
|P^k f(z)| &= |P^k f(z) - \ee f(X)| = |\ee f(Y_1\cdots Y_k z) - \ee f(Xz)|\\
&\le 2\|f\|_\infty d_{TV}(Y_1\cdots Y_k, X) \le 2\|f\|_\infty a e^{-bk}.
\end{align*}
This observation also gives
\begin{equation}\label{min1}
\begin{split}
&\ee|(f(x)-f(Yx)(P^kf(x)-P^kf(Yx))| \\
&\le 4\|f\|_\infty a e^{-bk} \ee|f(x)-f(Yx)| \le 4\|f\|_\infty a e^{-bk} \|f\|_Y.
\end{split}
\end{equation}
Now recall the assumption 2 that for any $y\in G$, $y^{-1}Yy$ has the same distribution as $Y$. Thus, for any $x,y\in G$,
\[
Pf(yx) = \ee(Yyx) = \ee(y y^{-1}Yyx) = \ee(yYx).
\] 
So, if we let $Y^\prime$ be an independent copy of $Y$, then
\begin{align*}
\ee(Pf(x)-Pf(Yx))^2 &= \ee(\ee(f(Y^\prime x) - f(YY^\prime x)|Y)^2) \\
&\le \ee(f(Y^\prime x) - f(YY^\prime x))^2 \\
&\le \sup_{y^\prime\in G} \ee(f(y^\prime x) - f(Yy^\prime x))^2 = \|f\|_Y^2.
\end{align*}
Thus, $\|Pf\|_Y \le \|f\|_Y$. Continuing by induction, we get $\|P^k f\|_Y \le \|f\|_Y$. Thus,
\begin{align}\label{min2}
&\ee|(f(x)-f(Yx)(P^kf(x)-P^kf(Yx))| \nonumber \\
&\le \bigl(\ee(f(x)-f(Yx))^2\bigr)^{1/2}\bigl(\ee(P^kf(x)-P^kf(Yx))^2\bigr)^{1/2} \nonumber \\
&\le \|f\|_Y \|P^k f\|_Y \le \|f\|_Y^2.
\end{align}
Using (\ref{min1}) and (\ref{min2}) in (\ref{vx}), we get
\begin{align}\label{vx1}
v(x) &\le \frac{1}{2}\sum_{k=0}^\infty \min\{\|f\|_Y^2, 4a \|f\|_\infty\|f\|_Y e^{-bk}\}\nonumber\\
&\le \frac{1}{2}\sum_{k=0}^\infty \min\{B^2, 4a ABe^{-bk}\}\nonumber \\
&= \frac{B^2}{2}\sum_{k=0}^\infty \min\{1, 4aAB^{-1} e^{-bk}\}.
\end{align}
We shall now compute a bound on the above sum. For ease of notation let $\beta = 4aAB^{-1}$, and let 
$\gamma = b^{-1}\log \beta$.
If $\beta < 1$, the sum is just a geometric series which is easy to evaluate. Now assume $\beta\ge 1$. Then $\gamma$ is nonnegative. Now, an easy verification shows that $\beta e^{-b\gamma} = 1$, and $1 \ge \beta e^{-bk}$ if and only if $k \ge \gamma$. Hence,
\begin{align*}
&\sum_{k=0}^\infty \min\{1, \beta e^{-bk}\} \le \gamma + 1+ \sum_{k\ge \gamma} \beta e^{-bk} \\
&\le \gamma +1+ \beta e^{-b\gamma} \sum_{r=0}^\infty e^{-br} = \gamma + 1+\frac{1}{1-e^{-b}}.
\end{align*}
To finish, we substitute this bound in (\ref{vx1}), and use Theorems \ref{weaklaw} and \ref{hoeffding} from~Chapter \ref{results}.\hfill $\Box$

\section{Application to random matrices and free \newline  probability}\label{freeprob}
Let $M$ be an $n \times n$ complex hermitian (i.e.\ self-adjoint) matrix. We shall fix the following notation for the rest of this section.
\begin{itemize}
\item The {\it empirical spectral measure} $M$ is the probability measure on $\rr$, denoted by $\mu_M$, which puts $1/n$ on each eigenvalue of $M$, repeated by multiplicities.
\item The empirical distribution function of $M$, denoted by $F_M$, is the distribution function corresponding to the empirical spectral measure.
\item Let $\Delta$ be a diagonal matrix whose diagonal elements are the eigenvalues of $M$. Then any hermitian matrix which has the same spectrum as $M$ can be written as $U\Delta U^*$, where $U$ is a unitary matrix. Thus, a uniform measure on the set of all such matrices is naturally induced by the Haar measure on the group of all unitary matrices of order $n$, which we shall denote by $\uu_n$. This probability measure will be denoted by $\rho_M$. 
\end{itemize}
A fundamental observation of Voiculescu \cite{voiculescu91} is the following: If $M$ and $N$ are two hermitian matrices of order $n$ with empirical measures $\mu_M$ and $\mu_N$, and $n$ is large, then the empirical measure $\mu_{M+N}$ of $M+N$ is ``approximately determined'' by $\mu_M$ and $\mu_N$ (irrespective of $M$ and $N$) for ``most choices of $M$ and $N$''.
The meaning of the phrases in quotes is made precise by the following result: 
\begin{thm}
\textup{[Voiculescu's result, as stated in Biane \cite{biane03}]} 
For each positive integer $n$, let $A_n$ and $B_n$ be two hermitian matrices, whose eigenvalues are bounded uniformly in $n$. Let $\mu_1$ and $\mu_2$ be two probability measures with compact supports on $\rr$, such that $\mu_{A_n} \ra \mu_1$ and $\mu_{B_n} \ra \mu_2$ weakly as $n \ra \infty$. Then there exists a probability measure, depending only on $\mu_1$ and $\mu_2$, denoted by $\mu_1\boxplus \mu_2$, such that
$\mu_{A_n^\prime+B_n^\prime} \ra \mu_1\boxplus \mu_2$, whenever $A_n^\prime$ and $B_n^\prime$ are random matrices chosen independently with distributions $\rho_{A_n}$ and $\rho_{B_n}$.
\end{thm}
Voiculescu \cite{voiculescu91} proved this striking result in a limiting form, using the method of moments and some concentration results of Szarek \cite{szarek90} and Gromov \& Milman~\cite{gromovmilman83}. The theorem was used by Voiculescu to establish a connection between free probability theory and random matrices, which resulted in an explosion of activity in the free probability literature. Another proof of Voiculescu's observation, using Stieltjes transforms, was given by Pastur and Vasilchuk \cite{pastur00}. A good review of the literature, as well as a good exposition, is given in Biane \cite{biane03}. Another useful reference is the lecture notes by Voiculescu \cite{voiculescu00}.

To the best of our knowledge, all the existing results are limiting statements; quantitative bounds on the concentration of $F_{M+N}$ given $F_M$ and $F_N$ for finite $n$ are not available in the literature. We shall now state and prove one such result, as an application of the machinery developed in the previous section.
\begin{thm}\label{free}
Let $\Delta_1$ and $\Delta_2$ be two $n\times n$ real diagonal matrices. Let $U$ and $V$ be independent Haar distributed random elements of $\uu_n$, the group of all unitary matrices of order $n$. Let
\[
H = U\Delta_1 U^* + V\Delta_2 V^*,
\]
and let $F_H$ be the empirical distribution function of $H$. Then, for every $x\in \rr$,
$\var(F_H(x)) \le \kappa n^{-1}\log n$.
where $\kappa$ is a universal constant \underline{not} depending on $n$, $\Delta_1$, $\Delta_2$ or $x$.
Moreover, we also have the concentration inequality 
\[
\pp\{|F_H(x)-\ee(F_H(x))| \ge t\} \le 2\exp\biggl(-\frac{nt^2}{2\kappa\log n}\biggr)
\]
for every $t\ge 0$, where $\kappa$ is the same as in the variance bound.
\end{thm}
To prove Theorem \ref{free}, we first need to establish a theorem about the concentration of the Haar measure on $\uu_n$. Existing results of the type discussed in section \ref{gpreview} cannot give concentration bounds for $F_H$, since they are based on the Hilber-Schmidt distance which is too crude for such a delicate problem. Instead, we shall try to find the concentration of the Haar measure with respect to the {\it rank distance}, defined as $d(M,N) := \mathrm{rank}(M-N)$. That this is indeed a metric, follows from the subadditivity of rank. The empirical distribution function is well-behaved with respect to this metric, as shown by the following lemma of Bai \cite{bai99}:
\begin{lmm}\label{bailmm}
\textup{[Bai \cite{bai99}, Lemma 2.2]}
Let $M$ and $N$ be two $n\times n$ hermitian matrices, with empirical distribution functions $F_M$ and $F_N$. Then 
\[
\|F_M-F_N\|_\infty \le \frac{1}{n}\mathrm{rank}(M-N).
\]
\end{lmm}
This lemma is an easy consequence of the interlacing inequalities for eigenvalues of hermitian matrices. It seems possible that this already existed in the literature before Bai \cite{bai99}, but we could not find any reference.

We shall use Theorem \ref{gpaction} to find the concentration of the Haar measure on $\uu_n$ with respect to the rank metric. For that purpose, we need a random walk which takes ``small steps'' with respect to that metric. 

Let $G = \uu_n$ and $X$ be a Haar-distributed random variable on $\uu_n$. We define the $Y$ required for generating the random walk for Theorem \ref{gpaction} as follows: Let $Y= I-(1-e^{i\varphi})uu^*$, where $u$ is drawn uniformly from the unit sphere in $\cc^n$, and $\varphi$ is drawn independently from the distribution on $[0, 2\pi)$ with density proportional to $(\sin(\varphi/2))^{n-1}$. Multiplication by $Y$ represents reflection across a randomly chosen subspace.
It is easy to verify that $Y\in \uu_n$. Now, for any $U\in \uu_n$,
\[
UYU^* = I- (1-e^{i\varphi})(Uu)(Uu)^*,
\]
and $Uu$ is again uniformly distributed over the unit sphere in $\cc^n$. Also, $Y^{-1}= Y^* = I - (1-e^{-i\varphi})uu^* = I-(1-e^{i(2\pi - \varphi)})uu^*$ has the same distribution as $Y$, since $2\pi -\varphi$ has the same distribution as $\varphi$. Hence $Y$ satisfies the properties $1$ and $2$ from section \ref{group}. Following a sketch of Diaconis \& Shahshahani \cite{diaconis86}, Ursula Porod \cite{porod96} proved the following result:
\begin{thm}
\textup{[Porod \cite{porod96}]}
Let $X,Y$ be as above. Let $Y_1,Y_2,\ldots,$ be i.i.d.\ copies of $Y$. There exists universal constants $\alpha, \beta, c_0$, such that whenever $n \ge 16$ and $k \ge \frac{1}{2}n\log n + c_0n$, we have
\begin{equation}\label{alphabeta}
d_{TV}(Y_1\cdots Y_k, X) \le \alpha n^{\beta/2} e^{-\beta k/n}.
\end{equation}
\end{thm}
Substituting $k=\frac{1}{2}n\log n + c_0 n$, we get $\alpha e^{-\beta c_0}$ on the right hand side. Thus by suitably increasing $\alpha$ such that $\alpha e^{-\beta c_0} \ge 1$, we can drop the condition that $k\ge \frac{1}{2}n\log n + c_0n$. 
Combining Porod's theorem with Theorem \ref{gpaction}, we get the following result about concentration of the Haar measure on $\uu_n$:
\begin{prop}\label{unitary}
Let $G=\uu_n$ and $X,Y$ be as above, with $n\ge 16$. Let $f:\uu_n \ra\rr$ be a function such that $\ee f(X)=0$. Let $\|f\|_Y = \sup_{U\in \uu_n} [\ee(f(U)-f(YU))^2]^{1/2}$. Let $A$ and $B$ be constants such that $\|f\|_\infty \le A$ and $\|f\|_Y \le B$. Let 
\[
C = \frac{nB^2}{2\beta}\biggl[\biggl(\log \frac{4\alpha n^\beta A}{B}\biggr)^+ + \frac{(\beta/n)(2-e^{-\beta/n})}{1-e^{-\beta/n}}\biggr],
\]
where $\alpha$ and $\beta$ are as in (\ref{alphabeta}).
Then $\var(f(X))\le C$, and for any $t\ge 0$, we have $\pp\{|f(X)|\ge t\} \le 2e^{-t^2/2C}$.
\end{prop}
We are now ready to finish the proof of Theorem \ref{free}:
\vskip.1in
\noindent {\bf Proof of Theorem \ref{free}.} We shall carry on with all the notation we have already defined in this section. 
The matrix $V^*HV = V^*U\Delta_1 U^*V + \Delta_2$ has the same spectrum as $H$. Also, $V^*U$ is again Haar distributed. Hence, we can write, without loss of generality, 
\[
H = X\Delta_1 X^* + \Delta_2,
\]
where $X$ follows the Haar distribution on $\uu_n$.  Now let
\[
H^\prime = (YX) \Delta_1 (YX)^* + \Delta_2.
\]
Recall that $Y = I-(1-e^{i\varphi})uu^*$, where $u$ is drawn from the uniform distribution on the unit sphere in $\cc^n$, and $\varphi$ is drawn independently from the distribution on $[0,2\pi)$ with density proportional to $(\sin(\varphi/2))^{n-1}$. Let $\delta = 1-e^{i\varphi}$. Then 
\begin{align*}
H-H^\prime &= X\Delta_1 X^* - (I-\delta uu^*)X\Delta_1 X^*(I-\bar{\delta} uu^*)\\
&=\delta Huu^* + \bar{\delta} uu^*H - |\delta|^2 uu^*Huu^*.
\end{align*}
The three summands are all of rank $1$. It follows by the subadditivity of rank that $\mathrm{rank}(H-H^\prime) \le 3$. Thus by Lemma \ref{bailmm}, we get 
\begin{equation}\label{fdiff}
\|F_H - F_{H^\prime}\|_\infty\le \frac{3}{n}. 
\end{equation}
Now fix a point $x\in \rr$, and let $f:\uu_n \ra \rr$ be the map which takes $X$ to $F_H(x)$. Then by (\ref{fdiff}), we have
\[
|f(X)-f(YX)| \le \frac{3}{n} \ \ \text{ for all possible values of $X$ and $Y$.}
\]
Thus, $\|f\|_Y \le 3/n$. Also, $\|f\|_\infty \le 1$. Thus, in Proposition \ref{unitary}, we get $C \le \kappa \log n + c$ for some universal constants $\kappa$ and $c$. By choosing $\kappa$ large enough, we can drop the assumption that $n\ge 16$ and also put $c=0$. 
This completes the proof. \hfill $\Box$

\bibliographystyle{plain}

\end{document}